\documentclass[letterpaper, times]{elsarticle}
\makeatletter
\def\ps@pprintTitle{%
 \let\@oddhead\@empty
 \let\@evenhead\@empty
 \def\@oddfoot{\thepage\hfill\footnotesize\itshape\today}
 \let\@evenfoot\@oddfoot}
\makeatother

\usepackage[margin=2.54cm]{geometry}
\usepackage[table]{xcolor}
\usepackage{mathtools}
\usepackage{graphicx}
\usepackage{amssymb}
\usepackage{caption}
\usepackage{subcaption}
\usepackage{tcolorbox}
\usepackage{scrextend}
\usepackage{listings}
\usepackage{bm}
\usepackage{float}
\usepackage[colorlinks]{hyperref}
\usepackage{multirow}
\usepackage[linesnumbered,lined,boxed,ruled,vlined]{algorithm2e}
\usepackage{afterpage}
\usepackage{etoolbox}
\usepackage{dutchcal}
\usepackage{makecell}
\usepackage{cleveref}
\usepackage{enumitem}
\usepackage[makeroom,smaller]{cancel}
\usepackage{stmaryrd}
\usepackage{placeins}
\usepackage[pages=some,placement=top]{background}

\backgroundsetup{
	scale=1,
	color=black,
	opacity=1.0,
	angle=0,
	contents={%
		\includegraphics[width=\paperwidth,height=2cm]{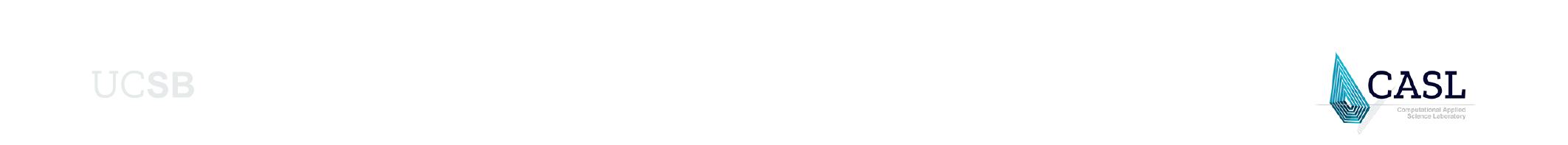}
	}%
}

\hypersetup{
	colorlinks=true,      		
	linkcolor=blue,				
	citecolor=blue,				
	filecolor=black,      		
	urlcolor=red,				
	bookmarks=false,
	pdffitwindow=true,
	pdfpagelayout=SinglePage
}

\definecolor{navy}{RGB}{2, 48, 71}
\definecolor{aqua}{RGB}{33, 158, 188}
\definecolor{cloud0}{RGB}{249, 249, 249}
\definecolor{cloud1}{RGB}{230, 234, 235}
\definecolor{goodgreen}{RGB}{30, 135, 74}
\definecolor{darkcoral}{RGB}{196, 52, 36}

\newcommand{\colorsection}[1]{\sffamily\color{navy}\section{#1}\color{black}\rmfamily}
\newcommand{\colorsubsection}[1]{\sffamily\color{navy}\subsection{#1}\color{black}\rmfamily}
\newcommand{\colorsubsubsection}[1]{\sffamily\color{navy}\subsubsection{#1}\color{black}\rmfamily}

\patchcmd{\abstract}{Abstract}{\sffamily{Abstract}\rmfamily}{}{}
\abstracttitle{\sffamily\textcolor{navy}{Abstract}}

\newcommand{\etal}{et al.\ }

\newcommand{\vv}[1]{ \boldsymbol{\mathbf{#1}} }

\newcommand{\eten}[2]{ #1\times 10^{#2} }

\newcommand{\keyval}[2]{#1{:}\,{#2}}

\let\OLDthebibliography\thebibliography
\renewcommand\thebibliography[1]{
  \OLDthebibliography{#1}
  \setlength{\parskip}{1pt}
  \setlength{\itemsep}{1pt plus 0.3ex}
}

\begin{document}

\title{\textcolor{aqua}{\sffamily\bfseries Error-correcting neural networks for two-dimensional curvature computation in the level-set method}}

\author[1]{\textcolor{gray}{Luis \'{A}ngel} Larios-C\'{a}rdenas\corref{cor1}}
\ead{lal@cs.ucsb.edu}

\author[1,2]{\textcolor{gray}{Fr\'{e}d\'{e}ric} Gibou}
\ead{fgibou@ucsb.edu}

\cortext[cor1]{Corresponding author}

\address[1]{\textcolor{darkgray}{Computer Science Department, University of California, Santa Barbara, CA 93106, USA}}
\address[2]{\textcolor{darkgray}{Mechanical Engineering Department, University of California, Santa Barbara, CA 93106, USA}}

\BgThispage


\begin{abstract}
We present an error-neural-modeling-based strategy for approximating two-dimensional curvature in the level-set method.  Our main contribution is a redesigned hybrid solver [Larios-C\'{a}rdenas and Gibou, \textit{J. Comput. Phys.} (May 2022), \href{https://doi.org/10.1016/j.jcp.2022.111291}{10.1016/j.jcp.2022.111291}] that relies on numerical schemes to enable machine-learning operations on demand.  In particular, our routine features double predicting to harness curvature symmetry invariance in favor of precision and stability.  The core of this solver is a multilayer perceptron trained on circular- and sinusoidal-interface samples.  Its role is to quantify the error in numerical curvature approximations and emit corrected estimates for select grid vertices along the free boundary.  These corrections arise in response to preprocessed context level-set, curvature, and gradient data.  To promote neural capacity, we have adopted sample negative-curvature normalization, reorientation, and reflection-based augmentation.  In the same manner, our system incorporates dimensionality reduction, well-balancedness, and regularization to minimize outlying effects.  Our training approach is likewise scalable across mesh sizes.  For this purpose, we have introduced dimensionless parametrization and probabilistic subsampling during data production.  Together, all these elements have improved the accuracy and efficiency of curvature calculations around under-resolved regions.  In most experiments, our strategy has outperformed the numerical baseline at twice the number of redistancing steps while requiring only a fraction of the cost.
\end{abstract}


\maketitle


\colorsection{Introduction}
\label{sec:Introduction}

Curvature plays a crucial role in free boundary problems (FBP) \cite{Friedman10} for its relation to surface tension in physics \cite{Osher1988, Popinet;NumModelsOfSurfTension;18} and its regularization property in optimization.  Such FBPs are ubiquitous in the physical sciences and engineering, with a broad range of applications in multiphase flows \cite{Sussman;Smereka;Osher:94:A-Level-Set-Approach, Sussman;Fatemi;Smereka;etal:98:An-Improved-Level-Se, Gibou;Chen;Nguyen;etal:07:A-level-set-based-sh, Theillard:2019aa, Losasso;Gibou;Fedkiw:04:Simulating-Water-and, Losasso:2006aa, Gibou:2019aa, Egan;etal;DirNumSimIncompFlowsOctree;2021}, solidification processes \cite{Chen;Min;Gibou:09:A-numerical-scheme-f, Papac;Gibou;Ratsch:10:Efficient-symmetric-, Papac;Helgadottir;Ratsch;etal:13:A-level-set-approach, Mirzadeh;Gibou:14:A-conservative-discr, Theillard;Gibou;Pollock:14:A-Sharp-Computationa}, and biological morphogenesis \cite{Boudon;etal;3DPlantMorphogenesis;2015, Ocko;Heyde;Mahadevan;MorphTermiteMounds;2019, AliasBuenzli20}.  Accurate curvature computations are thus critical to constructing sound FBP models.  For example, estimating curvature correctly \textit{at} the interface is essential for recovering continuous equilibrium solutions in multiphase flows \cite{Popinet;NumModelsOfSurfTension;18}.  Similarly, high precision is necessary to avoid erroneous pressure jumps that may affect breakup and coalescence \cite{Lervag;CalcCurvatureLSM;2014}.

One of the most widely used mathematical tools for solving FBPs is the level-set method \cite{Osher1988, Sethian:99:Level-set-methods-an, Osher;Fedkiw:02:Level-Set-Methods-an, GFO18}.  It belongs to a family of Eulerian formulations, such as the volume-of-fluid (VOF) \cite{Hirt;Nichols:81:Volume-of-Fluid-VOF-} and the phase-field \cite{QB10} representations, that capture and advect interfaces with the help of implicit functions.  Unlike Lagrangian schemes (e.g., front-tracking \cite{Tryggvason;Bunner;Esmaeeli;etal:01:A-Front-Tracking-Met}), implicit frameworks possess the natural ability to handle complex topological changes.  Their main advantage lies in the possibility of solving FBPs with no need for explicit procedures to reconstruct moving fronts.

The level-set method also makes it easy to estimate curvature at any point in the computational domain.  Yet, the lack of smoothness in the level-set field often compromises the accuracy of such approximations. \cite{Chene;Min;Gibou:08:Second-order-accurat, Popinet;NumModelsOfSurfTension;18}.  To address this problem, one typically uses iterative procedures to reinitialize the former into a signed distance function.  These reshaping algorithms can be as costly as the high-order schemes presented in \cite{Chene;Min;Gibou:08:Second-order-accurat} or as cheap as the first-order fast-sweeping methods described in \cite{Zhao:04:A-Fast-Sweeping-Meth, Detrixhe;Gibou;Min:13:A-parallel-fast-swee}.  However, none of these procedures can guarantee all-time reliable results, especially around under-resolved regions and in nonuniform grids.  Here, we consider high-order schemes to regularize the level-set field with only a few redistancing steps.  Our goal is to design a curvature solver working with inexpensive reinitialization to deliver accurate estimates at a low cost and particularly in under-resolved regions.

Further level-set research has been devoted to increasing curvature precision by numerical means.  Macklin, Lowengrub, and Lerv\r{a}g \cite{Macklin;Lowengrub;ImprovedCurvatureAppTumorGrowth;2006, Lervag;CalcCurvatureLSM;2014}, for instance, have proposed geometry-aware schemes based on interface reconstruction (IR) that depends on a smoothness quality metric.  The latter resembles a classifier that enables IR and curvature estimation with much smaller stencils or through interpolation from a recomputed level-set field.  In general, these approaches yield second-order accuracy but defeat the purpose of the implicit formulation.  Also, interface parametrization in three dimensions may become prohibitive, making this method unsuitable for highly deforming surfaces.  Our proposal is instead a hybrid solver.  Rather than solely relying on conventional means, we have instrumented an algorithm augmented with machine learning (ML) to improve numerical curvature only when necessary.

Our work is a revised version of the preliminary research presented in \cite{LALariosFGibou;LSCurvatureML;2021} and later extended in \cite{Larios;Gibou;HybridCurvature;2021}.  These approaches emerged as a level-set adaptation from Qi and coauthors' ML developments to approximate curvature in VOF technologies \cite{CurvatureML19}.  Qi and colleagues' idea was simple: to use shallow neural networks \cite{A18, Mehta19} to compute curvature at interface cells as a function of nine-point-stencil volume fractions.  First, they trained their models with samples extracted from circles of varying sizes.  Then, these were tested on static and moving boundaries, leading to reasonably good results.  Such a clever use of data-driven methods to find alternative views on FBPs has spurred a lot of joint ML and computational science research.  Patel \etal \cite{VOFCurvature3DML19}, for example, have redesigned Qi and coauthors' framework for three-dimensional interfaces.  Their systematic approach and optimized networks, in particular, have proven effective in several free-boundary experiments, including standard bubble simulations.  Likewise, Despr\'{e}s and Jourdren \cite{DespresJourdren;MLDesignOfVOF;20} have developed a family of VOF-ML schemes adapted to bi-material compressible Euler calculations on Cartesian grids.  The core of their system is a finite-volume flux function embodied in dense-layered neural models trained with lines, arcs, and corners.  Another piece of contemporary work is due to Ataei \etal \cite{NPLIC20}.  They have proposed a neural piece-wise linear IR method as a way of circumventing the VOF difficulties to locate moving fronts.  Yet, all this progress has not been exclusive to the VOF community.  Buhendwa \etal \cite{Buhendwa;Bezgin;Adams;IRinLSwithML;2021}, for instance, have devised a data-driven mechanism to estimate volume fractions and apertures in the level-set method.  Their ML solution has worked well for under-resolved regions while preserving the numerical convergence around well-resolved sectors.  Also, Fran\c{c}a and Oishi \cite{Franca;Oishi;MLCurvatureFrontTracking;2022} have recently introduced a neural strategy motivated by \cite{CurvatureML19} for calculating curvature in front-tracking.  Their multilayer perceptrons \cite{A18, Mehta19} simply consume normal- and tangential-vector angles to emit curvature estimates at the marker points.  As Fran\c{c}a and Oishi have shown, these networks can predict statistically accurate curvatures to replicate coalescence, bouncing, and separation behaviors in free-surface flows.

In a prior study \cite{LALariosFGibou;LSCurvatureML;2021}, we introduced level-set curvature neural networks analogous to \cite{CurvatureML19}.  Instead of volume fractions, these models ingested nine-point-stencil level-set values and produced the dimensionless curvature for interface nodes\footnote{We use the terms \textit{node}, \textit{grid point}, and \textit{vertex} interchangeably.}.  Eventually, further analysis of this network-only strategy led to the hybrid inference system presented in \cite{Larios;Gibou;HybridCurvature;2021}.  This new approach combined numerical schemes with multilayer perceptrons trained on circular and sinusoidal contours to improve curvature estimations as needed.  In addition, these models were trained only on half of the curvature spectrum, thus leading to higher capacity and fewer degrees of freedom in the observed input patterns.  The work in \cite{LALariosFGibou;LSCurvatureML;2021}, however, left a few questions unanswered.  The most important ones referred to overfitting and the lack of learning data-set scalability for highly resolved grids.  In this manuscript, we attempt to address these items by looking at the curvature problem from the error-correction perspective of \cite{Larios;Gibou;ECNetSemiLagrangian;2021}.  The latter is a level-set adaptation from the work of Pathak \etal \cite{Pathak;etal;MLToAugCoarseGridCFD;2020}, who developed a PDE-ML strategy to increase coarse turbulent-flow simulations' accuracy.  In \cite{Larios;Gibou;ECNetSemiLagrangian;2021}, we engineered a deep-learning-augmented semi-Lagrangian scheme inspired by image super-resolution methodologies \cite{Dong;Loy;He;SuperResolution;2014}.  The hearth of our algorithm was a feedforward model that compensated for numerical diffusion in response to context velocity, level-set, curvature, and gradient information.  To optimize such a network, first, we built data sets with samples collected from coarse and fine meshes advected simultaneously.  Then, we fitted the statistical model to the inter-grid level-set discrepancies and used it to minimize artificial viscosity and area loss.

The present study thus elaborates on an improved hybrid strategy \cite{Larios;Gibou;HybridCurvature;2021} that extends the error-correction notion of \cite{Larios;Gibou;ECNetSemiLagrangian;2021}.  Our main contribution is a redesigned curvature solver that relies on numerical schemes to enable ML operations on demand.  In particular, our procedure features double predicting to harness curvature symmetry invariance in favor of precision and stability.  As in \cite{Larios;Gibou;HybridCurvature;2021}, the core of this solver is a multilayer perceptron trained on circular- and sinusoidal-interface samples.  Its role is to quantify the error in numerical curvature approximations and emit corrected estimates for select grid vertices along the free boundary.  These corrections arise after processing context information, such as level-set, curvature, and gradient data.  To promote neural capacity, we have adopted sample negative-curvature normalization, reorientation, and reflection-based augmentation.  In the same manner, our system incorporates dimensionality reduction, well-balancedness, and regularization to minimize outlying effects.  Our network-training approach is likewise scalable across mesh sizes.  For this purpose, we have introduced dimensionless parametrization and probabilistic subsampling during data production.  Together, all these elements have improved the accuracy and efficiency of curvature calculations around under-resolved regions.  In fact, in most experiments, our proposed strategy has outperformed the numerical baseline at twice the number of redistancing steps while requiring only a fraction of the cost.

The rest of this manuscript is organized as follows.  First, we briefly review the level-set method.  After that, we give a thorough description of our strategy in \Cref{sec:Methodology}, including inference, data-set generation, and preprocessing algorithms, alongside relevant technical details.  \Cref{sec:Results} then provides a series of experiments evaluating the curvature hybrid solver's efficacy and efficiency.  Finally, we conclude our presentation in \Cref{sec:Conclusions} with a summary of our findings and a few pointers for future work.


\FloatBarrier
\colorsection{The level-set method}
\label{sec:TheLevelSetMethod}

The level-set method, introduced by Osher and Sethian \cite{Osher1988}, is an implicit formulation for capturing and transporting surfaces undergoing complex topological changes.  This framework represents the interface as the zero-isocontour $\Gamma(t) \doteq \{\vv{x}: \phi(\vv{x},t) = 0\}$ of a scalar, Lipschitz continuous relation $\phi(\vv{x},t): \mathbb{R}^{n+1} \mapsto \mathbb{R}$ known as the \textit{level-set function}\footnote{We will limit ourselves to $n = 2$ dimensions in this manuscript.  Also, we will omit the explicit time dependence of $\phi$ and $\Gamma$ for compactness.}.  The interface splits the computational domain $\Omega \subseteq \mathbb{R}^n$ into two non-overlapping regions, denoted by $\Omega^- \doteq \{\vv{x}: \phi(\vv{x}) < 0\}$ and $\Omega^+ \doteq \{\vv{x}: \phi(\vv{x}) > 0\}$.  Further, if $\vv{u}(\vv{x}): \mathbb{R}^n \mapsto \mathbb{R}^n$ is a well-defined velocity field, we can evolve the moving front and $\phi(\vv{x})$ by solving the Hamilton--Jacobi \textit{level-set equation}

\begin{equation}
\phi_t(\vv{x}) + \vv{u}(\vv{x})\cdot\nabla\phi(\vv{x}) = 0.
\label{eq:LevelSetEquation}
\end{equation}

\begin{figure}[t]
	\centering
	\includegraphics[width=\textwidth]{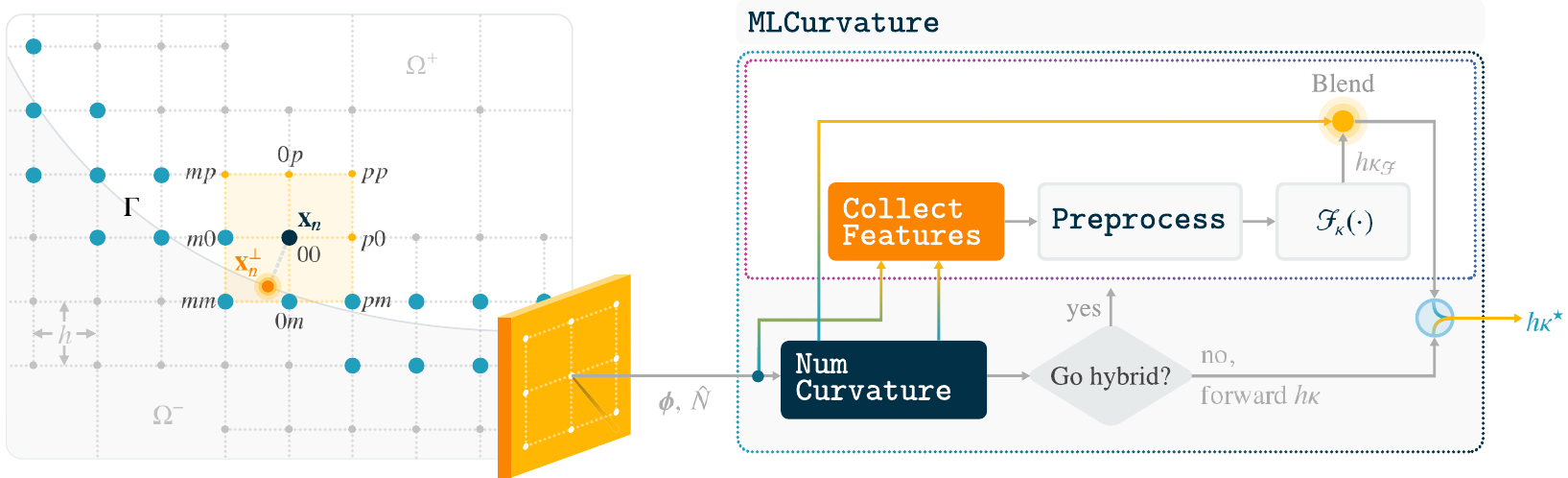}
	\caption{Overview of the hybrid inference system.  The left panel illustrates an adaptive mesh with grid points adjacent to the interface (in blue) suitable for our ML-augmented algorithm.  Among these, we highlight a node at $\vv{x}_\mathcal{n}$, its nine-point stencil, and its normal projection $\vv{x}_\mathcal{n}^\perp$ onto $\Gamma$, where curvature should be calculated.  We have also labeled each of the stencil vertices with $ij$, where $i,j = \{m,0,p\}$.  On the right, we show the {\tt MLCurvature()} routine, which couples the numerical procedure {\tt NumCurvature()} with the neural model $\mathcal{F}_\kappa(\cdot)$ to yield an improved (dimensionless) curvature estimation $h\kappa^\star$.  The full description of {\tt MLCurvature()} appears in \Cref{alg:MLCurvature}.  The latter uses as starting point the standard $h\kappa$ computation (see \cref{eq:Curvature}) to enable the hybrid portion of the routine (top right) on demand.  The synergy between the neural and the numerical components helps increase curvature accuracy beyond what one could attain with either unit alone.  (Color online.)}
	\label{fig:Overview}
\end{figure}

In addition, provided that $\phi(\vv{x})$ remains sufficiently smooth after advection, we can compute the normal vectors and curvature values at \textit{any} $\vv{x} \in \Omega$ by estimating

\begin{subequations}
\begin{align}
\hat{\vv{n}}(\vv{x}) &= \frac{\nabla\phi(\vv{x})}{\|\nabla\phi(\vv{x})\|_2} \quad \mathrm{and} \label{eq:Normal}\\
      \kappa(\vv{x}) &= \nabla\cdot\frac{\nabla\phi(\vv{x})}{\|\nabla\phi(\vv{x})\|_2} = \left. \frac{\phi_x^2 \phi_{yy} - 2\phi_x\phi_y\phi_{xy} + \phi_y^2 \phi_{xx}}{\left( \phi_x^2 + \phi_y^2 \right)^{3/2}}\right|_{\vv{x}} \label{eq:Curvature}
\end{align}
\label{eq:NormalAndCurvature}
\end{subequations}
with first- and second-order accurate finite differences \cite{Min;Gibou:07:A-second-order-accur, Chene;Min;Gibou:08:Second-order-accurat}.  The left panel in \cref{fig:Overview} illustrates the standard nine-point stencil used in the discrete approximation of $\kappa(\vv{x}_\mathcal{n})$ at node $\mathcal{n}$.

In general, there is an infinite number of relations whose zero level sets describe the same surface as $\Gamma$.  However, one opts for a signed distance level-set function in practice because it helps to simplify computations (e.g., $\|\nabla\phi(\vv{x})\|_2 = 1$), reduce artificial mass loss, and produce robust numerical results \cite{Sussman;Smereka;Osher:94:A-Level-Set-Approach}.  Also, previous studies \cite{LALariosFGibou;LSCurvatureML;2021, Larios;Gibou;HybridCurvature;2021} have shown that the signed distance property is beneficial for improving the accuracy and performance of neural networks that estimate curvature.  Thus, since the numerical solution of \cref{eq:LevelSetEquation} does not guarantee that $\phi(\vv{x})$ will remain a signed distance function, it is typical to recover such a feature by solving the pseudo-time transient \textit{reinitialization equation} \cite{Sussman;Smereka;Osher:94:A-Level-Set-Approach}

\begin{equation}
\phi_\tau(\vv{x}) + \texttt{sgn}\left(\phi^0(\vv{x})\right)(\|\nabla\phi(\vv{x})\|_2 - 1) = 0.
\label{eq:Reinitialization}
\end{equation}
Here, $\tau$ represents fictitious time, $\phi^0(\vv{x})$ is the level-set function before redistancing, and $\texttt{sgn}(\cdot)$ is a smoothed-out signum function.  Traditionally, one approaches \cref{eq:LevelSetEquation,eq:Reinitialization} with a TVD Runge--Kutta scheme in time \cite{Shu;Osher:89:Efficient-Implementa} and a high-order HJ-WENO discretization in space (e.g., see \cite{Jiang;Peng:00:Weighted-ENO-Schemes}).  As for \cref{eq:Reinitialization}, in particular, it is customary to stop the iterative process after $\nu$ steps, depending on the application and how close $\phi^0(\cdot)$ is to a signed distance function.  Min and Gibou \cite{Min;Gibou:07:A-second-order-accur} and Mirzadeh \etal \cite{Mirzadeh;etal:16:Parallel-level-set} have redesigned the advection and reinitialization procedures for adaptive Cartesian grids, which we have used as the basis for the present work.  For a thorough exposition of the level-set method, we refer the reader to the classic texts by Sethian \cite{Sethian:99:Level-set-methods-an} and Osher and Fedkiw \cite{Osher;Fedkiw:02:Level-Set-Methods-an}.  Similarly, the reader may consult \cite{GFO18} for a recent review of the level-set technologies and related state-of-the-art algorithms.


\FloatBarrier
\colorsection{Methodology}
\label{sec:Methodology}

Here, we describe our hybrid solver based on neural corrections.   Our approach improves the frameworks in \cite{Larios;Gibou;HybridCurvature;2021} and \cite{LALariosFGibou;LSCurvatureML;2021} by incorporating the error-correcting notion of \cite{Pathak;etal;MLToAugCoarseGridCFD;2020} and the ML enhancements of \cite{Larios;Gibou;ECNetSemiLagrangian;2021}.  The key components that make up the proposed curvature solver appear in \cref{fig:Overview}.

The main contribution of this study is the {\tt MLCurvature()} routine depicted in the right panel of \cref{fig:Overview}.  Its purpose is to manipulate level-set- and gradient-field information to estimate the dimensionless curvature \textit{at} the interface with higher precision than conventional means.  The heart of the hybrid system is the error-modeling multilayer perceptron $\mathcal{F}_\kappa(\cdot)$.  For a mesh size $h$, $\mathcal{F}_\kappa(\cdot)$ adjusts the numerical $h\kappa$ approximation interpolated at $\Gamma$ and computed in {\tt NumCurvature()} with finite-difference schemes.  Analogous to \cite{Larios;Gibou;HybridCurvature;2021}, we use $h\kappa$ to distinguish between well- and under-resolved regions and enable neural corrections only when deemed necessary.  However, rather than simply thresholding, we blend $h\kappa$\footnote{For us, $h\kappa$ denotes the \textit{numerical} approximation to dimensionless curvature at the interface, $h\kappa_\mathcal{F}$ is the \textit{neurally corrected} estimation, $h\kappa^\star$ is {\tt MLCurvature()}'s output, and $h\kappa^*$ represents the exact value.} and the corrected estimation, $h\kappa_\mathcal{F}$, to smooth the transition between these solutions when $h\kappa \approx 0$.  Also, compared to \cite{Larios;Gibou;HybridCurvature;2021}, we have redesigned the sample assembling process.  In particular, we have leveraged curvature invariance to reorient and reflect samples, which has led to more compact learning sets and more accurate approximations at the inference stage.  

Next, we discuss the {\tt MLCurvature()} function, emphasizing the ML component.  Then, we conclude this section with a technical description of the tools employed for building, optimizing, and deploying {\tt MLCurvature()} and $\mathcal{F}_\kappa(\cdot)$.


\colorsubsection{A curvature hybrid solver enhanced by error-correcting neural networks}
\label{subsec:CurvatureHybridSolverAndECNets}

Despite the straightforward formulation in \cref{eq:NormalAndCurvature}, the level-set method lacks inbuilt mechanisms to compute curvature accurately \textit{at} the interface.  This deficiency has been linked to the inability of traditional level-set schemes to recover equilibrium solutions in surface tension models \cite{Popinet;NumModelsOfSurfTension;18}.  To relieve such a limitation, one usually evaluates \cref{eq:Normal,eq:Curvature} in a subset of nodes around $\Gamma$.  Then, one interpolates $\kappa$ at their normal projections onto the free boundary (e.g., see \cite{Chen;Min;Gibou:09:A-numerical-scheme-f} and Algorithm \href{https://www.sciencedirect.com/science/article/pii/S002199911630242X\#fg0140}{5} in \cite{Mirzadeh;etal:16:Parallel-level-set} for curvature-driven implementations of the Stefan problem).  These steps, however, do not always guarantee satisfactory results in under-resolved regions and when $\phi(\vv{x})$ is insufficiently smooth.  In this study, our goal is to investigate whether ML (embodied in the {\tt MLCurvature()} module) can produce on-the-fly corrections to numerical $h\kappa$ interpolations at $\Gamma$.

The starting point within {\tt MLCurvature()} is the {\tt NumCurvature()} procedure.  {\tt NumCurvature()} is the typical function that uses nodal level-set values, $\vv{\phi}$, and unit normal vectors, $\hat{N}$, to yield (mean) curvature, $\vv{\kappa}$\footnote{For consistency, we represent one-element nodal variables as $M$-vectors in lowercase bold faces (e.g., $\vv{\phi}$) and variables with $d > 1$ values per node as $d$-by-$M$ matrices in caps (e.g., $\hat{N}$).  $M$ is the number of vertices, e.g., all independent nodes that a {\tt p4est} \cite{Burstedde;Wilcox;Ghattas:11:p4est:-Scalable-Algo} macromesh $\mathcal{G}$ is aware of.}.  Here, we are only interested in correcting $h\kappa$ for nodes with at least one neighbor across $\Gamma$.  Now, let $\mathcal{n}$ be one such interface node at $\vv{x}_\mathcal{n}$, as shown in the left panel of \cref{fig:Overview}.  First, the {\tt NumCurvature()} subroutine uses $\mathcal{n}$'s nine-point-stencil information to evaluate $\kappa$ with \cref{eq:Curvature}.  Then, we resort to bilinear interpolation \cite{Strain1999, Min;Gibou:07:A-second-order-accur} to estimate $h\kappa$ at

\begin{equation}
\vv{x}_\mathcal{n}^\Gamma = \vv{x}_\mathcal{n} - \phi(\vv{x}_\mathcal{n})\frac{\nabla\phi(\vv{x}_\mathcal{n})}{\|\nabla\phi(\vv{x}_\mathcal{n})\|_2},
\label{eq:ProjectionOntoGamma}
\end{equation}
where $\vv{x}_\mathcal{n}^\Gamma$ is a rough approximation of $\vv{x}_\mathcal{n}$'s perpendicular projection, $\vv{x}_\mathcal{n}^\perp$, onto $\Gamma$.  Furthermore, if $h\kappa^*$ represents the exact dimensionless curvature at $\vv{x}_\mathcal{n}^\perp$, we can characterize its deviation from $h\kappa$ with

\begin{equation}
h\kappa^* = h\kappa + \bar{\varepsilon},
\label{eq:HKError}
\end{equation}
where $\bar{\varepsilon}$ is the numerical error.

\begin{figure}[t]
	\centering
	\includegraphics[width=15.2cm]{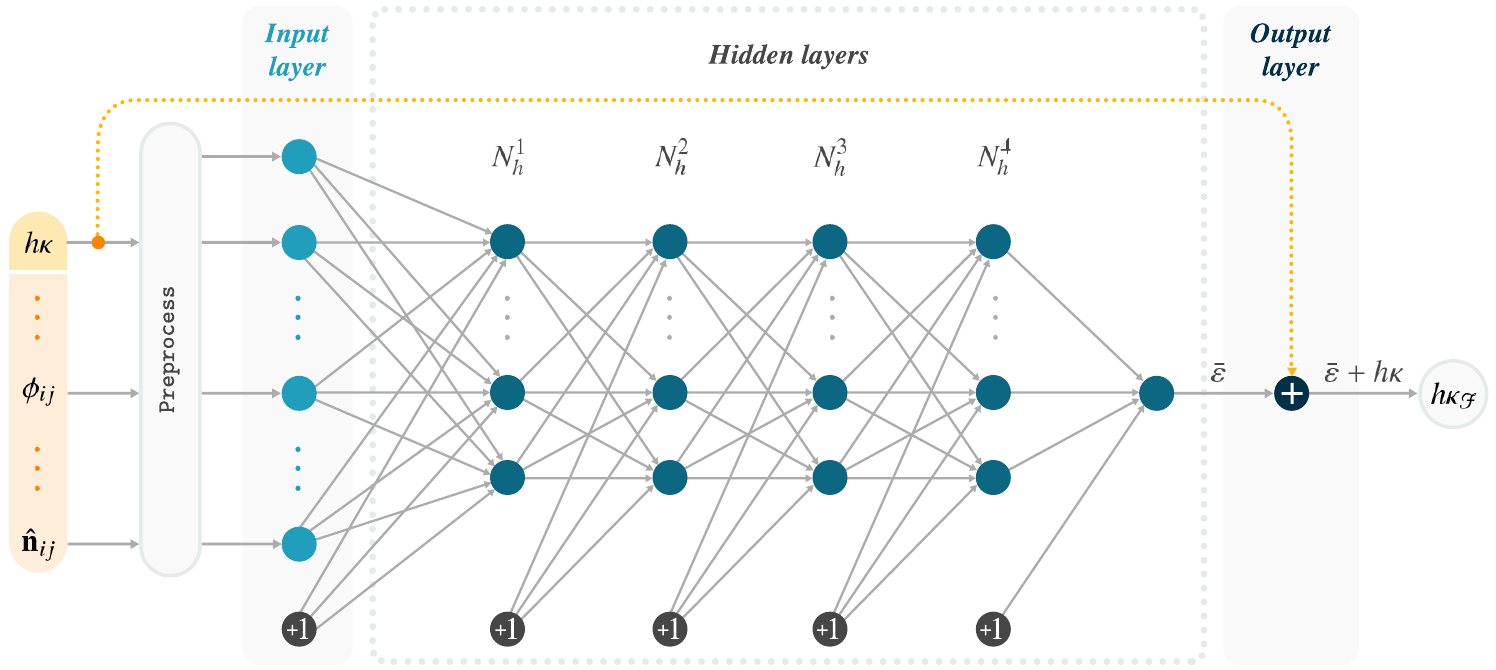}
	\caption{The error-correcting neural network, $\mathcal{F}_\kappa(\cdot)$, employed in the {\tt MLCurvature()} module of \cref{fig:Overview} and \Cref{alg:MLCurvature}.  We also show the {\tt Preprocess()} subroutine (see \Cref{alg:Preprocess}) that transforms the level-set values $\phi_{ij}$, unit normal vectors $\hat{\vv{n}}_{ij}$, and numerical $h\kappa$ bilinearly interpolated at the closest point on $\Gamma$.  $\mathcal{F}_\kappa(\cdot)$ outputs $h\kappa_\mathcal{F}$, which is subject to linear blending with $h\kappa$ if the latter is close to zero.  (Color online.)}
	\label{fig:ECNet}
\end{figure}

In this study, we propose a neural function, $\mathcal{F}_\kappa(\cdot)$, to quantify $\bar{\varepsilon}$ in \cref{eq:HKError}.  $\mathcal{F}_\kappa(\cdot)$ is a multilayer perceptron resembling the error-correcting network introduced by \cite{Larios;Gibou;ECNetSemiLagrangian;2021} for semi-Lagrangian transport.  \Cref{fig:ECNet} outlines the structure of our neural network.  Unlike \cite{LALariosFGibou;LSCurvatureML;2021, Larios;Gibou;HybridCurvature;2021} in the level-set method and \cite{CurvatureML19, VOFCurvature3DML19} in the VOF framework, $\mathcal{F}_\kappa(\cdot)$ digests more statistical information (e.g., gradient data) than just level-set or volume-fraction values.  In addition, besides its $\bar{\varepsilon}$-reconstructing feedforward architecture, our model features a skip connection that carries the $h\kappa$ input to a non-trainable unit that assembles $h\kappa_\mathcal{F} \doteq \bar{\varepsilon} + h\kappa$.  Although it is impossible to recover the exact $\bar{\varepsilon}$, we expect $|h\kappa^* - h\kappa_\mathcal{F}|$ to be smaller than $\bar{\varepsilon}$.  \Cref{fig:ECNet} also shows the {\tt Preprocess()} component that transforms incoming data into a suitable form for the downstream estimator.  As proven in \cite{LALariosFGibou;LSCurvatureML;2021, Larios;Gibou;HybridCurvature;2021, Larios;Gibou;ECNetSemiLagrangian;2021}, preprocessing is beneficial for learning algorithms because it facilitates convergence and reduces topological complexity \cite{scikit-learn11}.  We discuss the {\tt Preprocess()} subroutine later in \Cref{subsec:TechnicalAspects}.

\begin{figure}[t]
	\centering
	\begin{subfigure}[b]{5cm}
		\includegraphics[width=\textwidth]{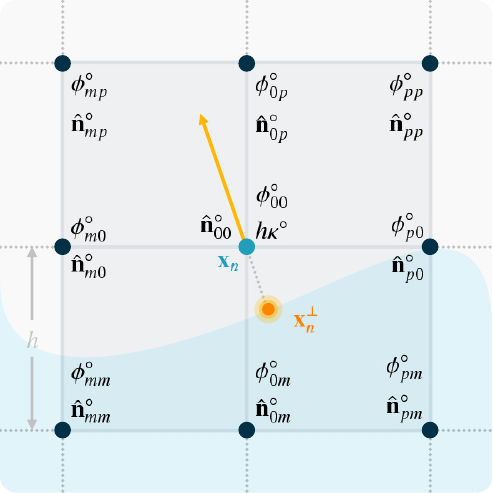}
		\caption{\footnotesize Original data packet $\mathcal{p}^\circ$}
		\label{fig:DataPacket.Original}
	\end{subfigure}
	~
	\begin{subfigure}[b]{5cm}
		\includegraphics[width=\textwidth]{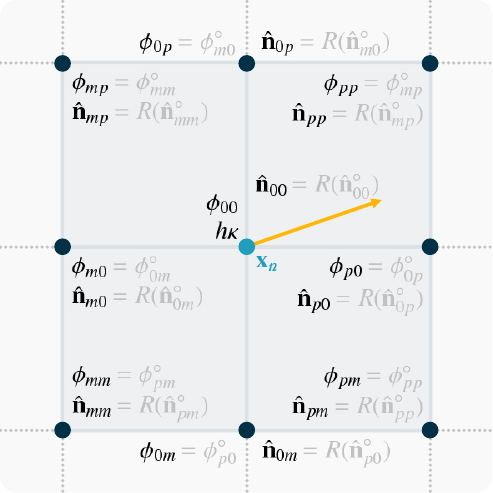}
		\caption{\footnotesize Reoriented data packet $\mathcal{p}$ (standard form)}
		\label{fig:DataPacket.Reoriented}
	\end{subfigure}
	~
	\begin{subfigure}[b]{5cm}
		\includegraphics[width=\textwidth]{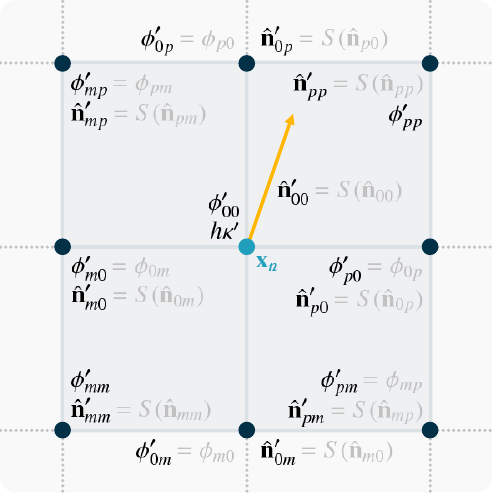}
		\caption{\footnotesize Augmented data packet $\mathcal{p}'$}
		\label{fig:DataPacket.Augmented}
	\end{subfigure}
	\caption{A data packet associated with interface node $\mathcal{n}$ located at $\vv{x}_\mathcal{n}$ with an $h$-uniform nine-point stencil at the maximum level of refinement.  The original sampled data appears in (a), denoted with a $\circ$ superscript.  The ML engine digests only numerical values font-colored in black.  For reference, we picture $\vv{x}_\mathcal{n}^\perp$, which we estimate with $\vv{x}_\mathcal{n}^\Gamma$ to approximate $h\kappa^\circ$ via bilinear interpolation.  After reorientation, we obtain the data packet in standard form in (b), where $R$ is a rotation transformation (assuming that $h\kappa^\circ \leqslant -h\kappa_{\min}^{low}$).  Invariant data to rotation transfers directly from (a) to (b), and we omit their equalities (e.g., $\phi_{00} = \phi_{00}^\circ$).  (c) exemplifies data augmentation by reflection about the line $y = x + \beta$, with $\beta \in \mathbb{R}$, going through $\vv{x}_\mathcal{n}$.  A permutation matrix $S$ swaps the unit normal vector components in (b).  Invariant data to reflection transfers directly from (b) to (c) (e.g., $\phi_{00}' = \phi_{00}$).  Augmented features are labeled with prime superscripts.  (Color online.)}
	\label{fig:DataPacket}
\end{figure}

The input vector to $\mathcal{F}_\kappa(\cdot)$ is the transformed version of data directly involved in \cref{eq:Curvature} for node $\mathcal{n}$.  This information includes level-set values, unit gradient vector components, and the numerical $h\kappa$ estimation.  Formally, we group all these entries into a \textit{data packet} defined by

\begin{equation}
\mathcal{p} = \left(\begin{array}{rl}
	        \phi_{ij}\!: & \textrm{nine-point stencil level-set values} \\
	\hat{\vv{n}}_{ij}\!: & \textrm{nine-point stencil unit normal vectors} \\
	          h\kappa\!: & \textrm{dimensionless curvature interpolated at } \vv{x}_\mathcal{n}^\Gamma
\end{array}\right) \in \mathbb{R}^{28},
\label{eq:DataPacket}
\end{equation}
where $\psi_{ij} = \psi(x_0 + ih, y_0 + jh)$ is a vertex feature with $i,j \in \{\keyval{m}{-1},\; \keyval{0}{0},\; \keyval{p}{+1}\}$\footnote{This is a dictionary mapping coordinate codes to numerical values; for example, $m$ maps to $-1$.}, and $\vv{x}_\mathcal{n} = (x_0, y_0)$ are $\mathcal{n}$'s coordinates.  In practice, we sample most of the information with the {\tt CollectFeatures()} function in \cref{fig:Overview}.  Also, notice that we require $h$-uniform nine-point stencils, especially when working with adaptive Cartesian grids \cite{Min;Gibou:07:A-second-order-accur, Mirzadeh;etal:16:Parallel-level-set}.  The schematic in \cref{fig:DataPacket.Original} illustrates a data packet in its spatial context.  We have introduced the $\circ$ superscript to distinguish it from the \textit{standard form} supplied to the {\tt Preprocess()} subroutine.

In prior research, we showed that inference and learning could substantially improve if one narrows down the curvature problem to only the negative spectrum.  More specifically, in \cite{Larios;Gibou;HybridCurvature;2021}, we first employed the numerical $h\kappa$ estimation as a convexity indicator.  Then, we flipped the sign of the level-set stencil and the inferred curvature accordingly.  In \cite{Larios;Gibou;ECNetSemiLagrangian;2021}, we also exploited curvature rotation invariance to minimize feature variations by \textit{reorienting} data packets into canonical forms.  This technique is a convenient normalization procedure in computer vision \cite{Turk;Pentland;Eigenfaces;1991, Parker;CS170A;2016} and proved effective for our image super-resolution approach.  Here, we have adopted these methodologies to reduce degrees of freedom and facilitate neural network design and data-set composition.  Thus, to normalize a data packet (as seen in \cref{fig:DataPacket.Original}), we first negate $\mathcal{p}^\circ$'s components if $h\kappa^\circ \geqslant h\kappa_{\min}^{low}$, for some threshold $h\kappa_{\min}^{low} > 0$.  Next, we reorient $\mathcal{p}^\circ$ into its standard form $\mathcal{p}$ by rotating the stencil until the angle between the horizontal and $\mathcal{p}.\hat{\vv{n}}_{00}$ lies between $0$ and $\pi/2$.  \Cref{fig:DataPacket.Reoriented} portrays the resulting data packet after reorienting $\mathcal{p}^\circ$, assuming that $h\kappa^\circ \leqslant -h\kappa_{\min}^{low}$.  Reorientation does not change the inherent level-set values in $\mathcal{p}$; however, it affects unit normal vectors since it entails a rigid-body transformation by a rotation $R(\theta)$, where $\theta = \pm k\pi/2$ and $k \in \{0, 1, 2\}$.


\begin{algorithm}[!t]
\small
\SetAlgoLined
\SetKwFunction{numcurvature}{NumCurvature}
\SetKwFunction{interpolate}{Interpolate}
\SetKwFunction{collectfeatures}{CollectFeatures}
\SetKwFunction{preprocess}{Preprocess}
\SetKwFunction{sign}{Sign}

\KwIn{$\mathcal{n}$: node object; $\mathcal{F}_\kappa(\cdot)$: error-correcting neural network; $\vv{\phi}$: nodal level-set values; $\hat{N}$: nodal unit normal vectors; $h$: mesh size; $h\kappa_{\min}^{low}$: minimum $|h\kappa|$ to enable neural inference; $h\kappa_{\min}^{up}$: $|h\kappa|$'s upper bound for blending numerical with neurally corrected approximation.}
\KwResult{$h\kappa^\star$: dimensionless curvature at $\vv{x}_\mathcal{n}^\Gamma$.}
\BlankLine

\tcp{Numerical computation}
$K \leftarrow$ \numcurvature{\normalfont $\mathcal{n}.\texttt{stencil}$, $\vv{\phi}$, $\hat{N}$}\tcp*[r]{$\kappa$ for $\mathcal{n}$'s stencil using \cref{eq:Curvature}}
$\vv{x}_\mathcal{n}^\Gamma \leftarrow \mathcal{n}.\vv{x} - \vv{\phi}[\mathcal{n}]\hat{N}[\mathcal{n}]$\tcp*[r]{See \cref{eq:ProjectionOntoGamma}}
$h\kappa \leftarrow$ $h \cdot$\interpolate{\normalfont $\mathcal{n}.\texttt{stencil}$, $K$, $\vv{x}_\mathcal{n}^\Gamma$}\;
\BlankLine
			
\tcp{Selectively enabling neural correction}
\eIf{$|h\kappa| \geqslant h\kappa_{\min}^{low}$}{
	$\mathcal{p} \leftarrow$ \collectfeatures{\normalfont $\mathcal{n}.\texttt{stencil}$, $\vv{\phi}$, $\hat{N}$}\tcp*[r]{Populate $\mathcal{p}$'s level-set values and normals}
	$\mathcal{p}.h\kappa \leftarrow h\kappa$\tcp*[r]{Add numerical $h\kappa$ too (see \cref{fig:DataPacket.Original})}
	\BlankLine
	
	\tcp{Produce two samples for $\mathcal{n}$}
	transform $\mathcal{p}$ so that $\mathcal{p}.h\kappa$ is negative\;
	reorient $\mathcal{p}$ so that the angle of $\mathcal{p}.\hat{\vv{n}}_{00}$ lies between $0$ and $\pi/2$\tcp*[r]{See \cref{fig:DataPacket.Reoriented}}
	$h\kappa_\mathcal{F} \leftarrow \mathcal{F}_\kappa([$\preprocess{$\mathcal{p}$, $h$}, $\mathcal{p}.h\kappa])$\tcp*[r]{First try}
	\BlankLine
	
	let $\mathcal{p}'$ be the reflected data packet about the line $y = x + \beta$ going through $\mathcal{n}.\vv{x}$\tcp*[r]{See \cref{fig:DataPacket.Augmented}}
	$h\kappa_\mathcal{F}' \leftarrow \mathcal{F}_\kappa([$\preprocess{$\mathcal{p}'$, $h$}, $\mathcal{p}'.h\kappa])$\tcp*[r]{Second try}
	\BlankLine
	
	$h\bar{\kappa}_\mathcal{F} \leftarrow \frac{1}{2}(h\kappa_\mathcal{F} + h\kappa_\mathcal{F}')$\tcp*[r]{Average neural prediction}
	\BlankLine
	
	\tcp{Linearly blending neural and numerical estimations near zero}
	\If{$|h\kappa| \leqslant h\kappa_{\min}^{up}$}{
		$\lambda \leftarrow (h\kappa_{\min}^{up} - |h\kappa|)/(h\kappa_{\min}^{up} - h\kappa_{\min}^{low})$\;
		$h\bar{\kappa}_\mathcal{F} = (1 - \lambda)h\bar{\kappa}_\mathcal{F} + \lambda( -|h\kappa|)$\;
	}
	\BlankLine
	
	$h\kappa^\star \leftarrow$\sign{$h\kappa$}$ \cdot\,|h\bar{\kappa}_\mathcal{F}|$\tcp*[r]{Restore sign}
}{
	$h\kappa^\star \leftarrow h\kappa$\;
}
\BlankLine

\Return $h\kappa^\star$\;

\caption{\small $h\kappa^\star \leftarrow$ {\tt MLCurvature(}$\mathcal{n}$, $\mathcal{F}_\kappa(\cdot)$, $\vv{\phi}$, $\hat{N}$, $h$, $h\kappa_{\min}^{low}$, $h\kappa_{\min}^{up}${\tt )}: Compute the dimensionless curvature for the interface node $\mathcal{n}$ using the standard scheme with error correction provided by $\mathcal{F}_\kappa(\cdot)$.}
\label{alg:MLCurvature}
\end{algorithm}

We close our system's overview with a concise description {\tt MLCurvature()} in \Cref{alg:MLCurvature}.  We provide the procedure for a single node, but one can be adapt it for batches of vertices in a grid $\mathcal{G}$.  The main ingredient in \Cref{alg:MLCurvature} is $\mathcal{F}_\kappa(\cdot)$, which one must optimize for a prescribed mesh size $h$.  Later, in \Cref{subsec:Training}, we describe the methodology for training $\mathcal{F}_\kappa(\cdot)$.  Besides the (reinitialized) nodal level-set values and unit gradient vectors, the {\tt MLCurvature()} routine requires $h\kappa_{\min}^{low}$ and $h\kappa_{\min}^{up}$.  These parameters define the minimum $|h\kappa|$ to trigger ML evaluations and the upper bound to transition smoothly from $h\kappa$ to $h\kappa_\mathcal{F}$.  $h\kappa_{\min}^{low}$, in particular, is critical to restrain data-set size, too, as discussed in \Cref{subsubsec:CircularInterfaceDataSetConstruction,subsubsec:SinusoidalInterfaceDataSetConstruction}.  Experimentally, we have found that $h\kappa_{\min}^{low} = 0.004$ and $h\kappa_{\min}^{up} = 0.007$ yield satisfactory results for any of the assessed grid resolutions in \Cref{sec:Results}.

As noted above, \Cref{alg:MLCurvature} is suitable for grid points satisfying at least one of the four conditions: $\phi(\mathcal{n}.x,\, \mathcal{n}.y)\cdot\phi(\mathcal{n}.x \pm h,\, \mathcal{n}.y) \leqslant 0$ or $\phi(\mathcal{n}.x,\, \mathcal{n}.y)\cdot\phi(\mathcal{n}.x,\, \mathcal{n}.y \pm h) \leqslant 0$, where $\mathcal{n}.\vv{x} = (\mathcal{n}.x,\, \mathcal{n}.y) = \vv{x}_\mathcal{n}$ are the vertex locations in \cref{fig:Overview}.  For each such node $\mathcal{n}$, we first estimate $h\kappa$ at $\vv{x}_\mathcal{n}^\Gamma$.  Then, we determine whether we should evaluate $\mathcal{F}_\kappa(\cdot)$ to improve $h\kappa$ in a relatively steep region.  This decision appears in \cref{fig:Overview} as the right-end multiplexor.  If $|h\kappa| \geqslant h\kappa_{\min}^{low}$, we continue to populate and normalize $\mathcal{n}$'s data packet $\mathcal{p}$ to comply with the characterization in \cref{fig:DataPacket.Reoriented}.  Using the features in $\mathcal{p}$, we next construct two samples: one with the normalized information and another with the augmented data packet $\mathcal{p}'$.  The example in \cref{fig:DataPacket.Augmented} illustrates the reflected data packet about the line $y = x + \beta$ for the reoriented $\mathcal{p}$ in \cref{fig:DataPacket.Reoriented}.  This form of data augmentation exploits curvature reflection invariance and has been used successfully in ML-guided IR \cite{Buhendwa;Bezgin;Adams;IRinLSwithML;2021} and semi-Lagrangian advection \cite{Larios;Gibou;ECNetSemiLagrangian;2021}.  In our case, this symmetry-based invariance allows us to improve accuracy by averaging $h\kappa_\mathcal{F}$ and $h\kappa_\mathcal{F}'$ into $h\bar{\kappa}_\mathcal{F}$ for the corresponding

\begin{equation*}
[\texttt{Preprocess(}\mathcal{p},\, h\texttt{)},\, \mathcal{p}.h\kappa] \quad \textrm{and} \quad
[\texttt{Preprocess(}\mathcal{p}',\, h\texttt{)},\, \mathcal{p}'.h\kappa]
\end{equation*}
Both samples are two-part input tensors conforming to $\mathcal{F}_\kappa(\cdot)$'s expected format in \cref{fig:ECNet}.

Finally, we blend $-|h\kappa|$ and $h\bar{\kappa}_\mathcal{F}$ linearly if $h\kappa_{\min}^{low} \leqslant |h\kappa| \leqslant h\kappa_{\min}^{up}$.  By doing so, we rid our solver from sharp curvature transitions and minimize the effects of the neural divergence documented in Section 4.2 of \cite{Larios;Gibou;HybridCurvature;2021}.  The last step in \Cref{alg:MLCurvature} restores the curvature sign.  Then, we return the improved dimensionless curvature to the calling function as $h\kappa^\star$.


\colorsubsection{Training}
\label{subsec:Training}

\begin{figure}[t]
	\centering
	\begin{subfigure}[b]{7.5cm}
		\includegraphics[width=\textwidth]{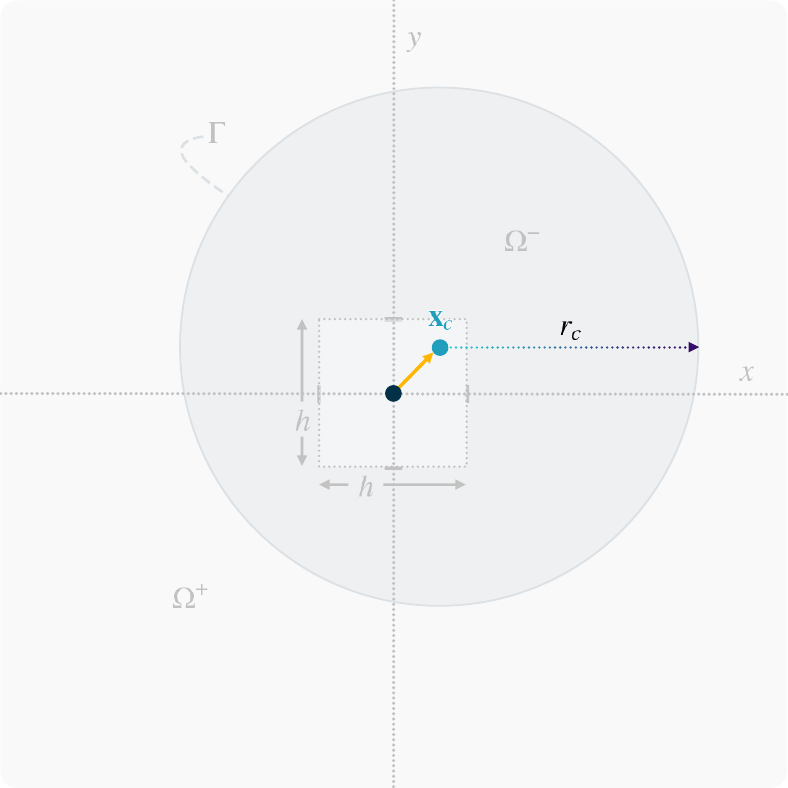}
		\caption{\footnotesize Circular interface}
		\label{fig:TrainingInterfaces.Circular}
	\end{subfigure}
	~
	\begin{subfigure}[b]{7.5cm}
		\includegraphics[width=\textwidth]{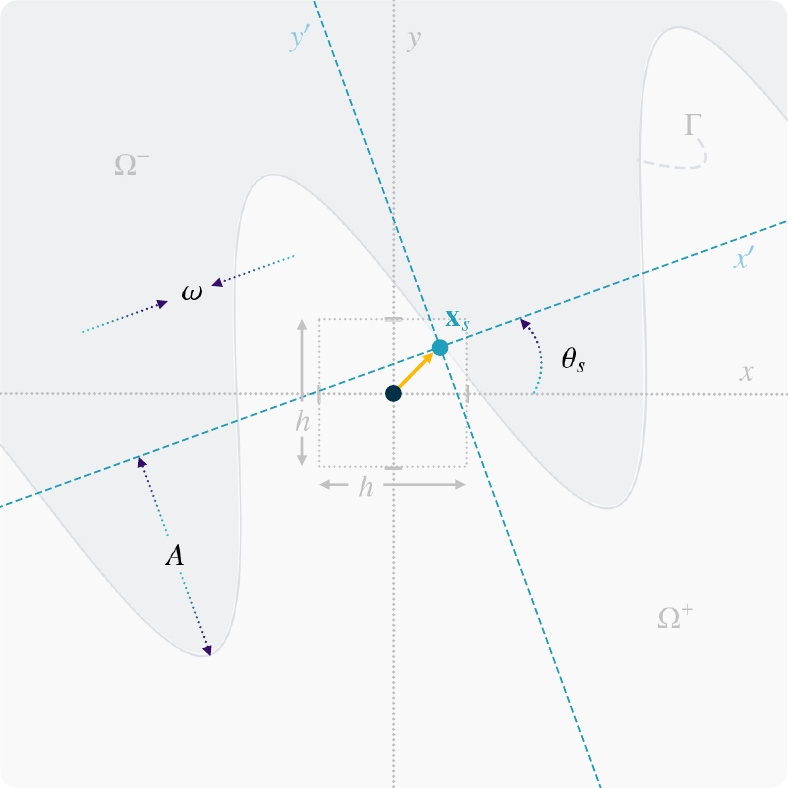}
		\caption{\footnotesize Sinusoidal interface}
		\label{fig:TrainingInterfaces.Sinusoidal}
	\end{subfigure}
	\caption{Training interfaces.  The innermost $h$-by-$h$ region defines the range on $x$ and $y$ where we randomly place the circle's center, $\vv{x}_c$, or the origin, $\vv{x}_s$, of the sine wave's canonical frame.  The parameters we vary to diversify sampling appear in black and blue.  (Color online.)}
	\label{fig:TrainingInterfaces}
\end{figure}

To train $\mathcal{F}_\kappa(\cdot)$, we assemble a data set $\mathcal{D}$ with circular- and sinusoidal-interface samples.  Our methodology is analogous to \cite{Larios;Gibou;HybridCurvature;2021}, with some critical changes that enable scalability across grid resolutions.  One such difference is the omission of exact-signed-distance-function sampling, which helps cut $|\mathcal{D}|$ by half.  The second difference is the introduction of $h\kappa_{\max}^*$ and $h\kappa_{\min}^*$\footnote{One should choose $h\kappa_{\min}^* \leqslant h\kappa_{\min}^{low}$, where the latter is the lower bound from \Cref{alg:MLCurvature}.}, which restrain $\mathcal{F}_\kappa(\cdot)$'s scope and makes data-set construction independent of the mesh size.  \Cref{fig:TrainingInterfaces} shows the two kinds of interfaces and the shape parameters considered in this study.  Next, we describe the sampling process to build $\mathcal{D}$ by joining the tuples extracted from circumferences ($\mathcal{D}_c$) and sine waves ($\mathcal{D}_s$).

\colorsubsubsection{Circular-interface data-set construction}
\label{subsubsec:CircularInterfaceDataSetConstruction}

The procedure for building $\mathcal{D}_c$ entails sampling from circular-interface level-set functions in a randomized fashion.  As in \cite{CurvatureML19, LALariosFGibou;LSCurvatureML;2021, Larios;Gibou;HybridCurvature;2021}, we have chosen circular interfaces for their simplicity and because $\kappa$ is just the reciprocal of the radius.  In our approach, we space radii based on a uniform curvature distribution.  Also, we ensure that all (discrete) target $h\kappa^*$ values have a fair frequency in $\mathcal{D}_c$.  Our strategy is to provide each class (i.e., $h\kappa^*$) with a sufficient probability of being accounted for during training.  And, to do so, we prescribe the number of data packets to produce for any $h\kappa^*$ by using a trapezoidal density function.  Our heuristic reduces the over-representation of shallow-curvature samples and provides more visibility to small-circle patterns.

\Cref{alg:GenerateCircularDataSet} summarizes the steps to assemble $\mathcal{D}_c$.  Its basis is the circular-interface (non-signed distance) level-set function defined by

\begin{equation}
\phi_c(\vv{x}) = \|\vv{x} - \vv{x}_c\|_2^2 - r_c^2,
\label{eq:CircularLevelSetFunction}
\end{equation}
where $\vv{x}_c$ and $r_c$ are the center and radius shown in \cref{fig:TrainingInterfaces.Circular}.  The goal of \Cref{alg:GenerateCircularDataSet} is to vary $\vv{x}_c$ (randomly) and $r_c$ (deterministically) to widen the range of patterns in $\mathcal{D}_c$.


\begin{algorithm}[!t]
\small
\SetAlgoLined
\SetKwFunction{linspace}{Linspace}
\SetKwFunction{circlelevelset}{CircleLevelSet}
\SetKwFunction{generategrid}{GenerateGrid}
\SetKwFunction{evaluate}{Evaluate}
\SetKwFunction{reinitialize}{Reinitialize}
\SetKwFunction{computenormals}{ComputeNormals}
\SetKwFunction{numcurvature}{NumCurvature}
\SetKwFunction{getnodesnexttogamma}{GetNodesNextToGamma}
\SetKwFunction{collectfeatures}{CollectFeatures}
\SetKwFunction{randomsamples}{RandomSamples}

\KwIn{$\eta$: maximum level of refinement per unit-square quadtree; $h\kappa_{\min}^*$ and $h\kappa_{\max}^*$: minimum and maximum target $|h\kappa^*|$; {\tt CPH}: number of circles per $h$; {\tt SPH2}: number of samples per $h^2$; $\nu$: number of iterations for level-set reinitialization; {\tt KeepEveryX}: subsampling constant.}
\KwResult{$\mathcal{D}_c$: data set of circular-interface samples.}
\BlankLine

$h \leftarrow 2^{-\eta}$\tcp*[r]{Mesh size}
$\kappa_{\min}^* \leftarrow h\kappa_{\min}^*/h; \quad \kappa_{\max}^* \leftarrow h\kappa_{\max}^*/h; \quad 
r_{\min} \leftarrow 1/\kappa_{\max}^*; \quad r_{\max} \leftarrow 1/\kappa_{\min}^*$\tcp*[r]{Target curvature and radius bounds}
{\tt NC} $\leftarrow \left\lceil \texttt{CPH} \cdot \left(\frac{r_{\max} - r_{\min}}{h} + 1\right) \right\rceil$\tcp*[r]{Number of circles}
\BlankLine

$\bar{r} = \frac{1}{2}(r_{\min} + r_{\max})$\tcp*[r]{Mean radius modulates sampling for each $r \in [r_{\min}, r_{\max}]$}
$\texttt{AvgSPR} \leftarrow \frac{1}{\texttt{KeepEveryX}}\left\lceil \texttt{SPH2} \cdot \frac{\pi}{h^2}\left(\bar{r}^2 - (\bar{r} - h)^2\right) \right\rceil$\;
$\texttt{SPR} \leftarrow$ \linspace{\normalfont $0.75\cdot\texttt{AvgSPR}$, $1.25\cdot\texttt{AvgSPR}$, $\texttt{NC}$}\tcp*[r]{Number of samples per radius}
\BlankLine
			
\tcp{Data generation loop}
$\mathcal{D}_c \leftarrow \varnothing$\;
$\texttt{TgtK} \leftarrow$ \linspace{\normalfont $\kappa_{\max}^*$, $\kappa_{\min}^*$, $\texttt{NC}$}\tcp*[r]{Equally spaced $\kappa^* \in [\kappa_{\min}^*, \kappa_{\max}^*]$}
\For{\normalfont $\texttt{c} \leftarrow 0, 1, \dots, \texttt{NC} - 1$}{
	$\kappa^* \leftarrow \texttt{TgtK}[\texttt{c}]; \quad r_c \leftarrow 1/\kappa^*$\tcp*[r]{Target curvature and associated radius}
	\BlankLine
	
	\tcp{Collect samples for radius $r_c$ into the set $\mathcal{S}$}
	$\mathcal{S} \leftarrow \varnothing$\;
	\While{\normalfont $|\mathcal{S}| < 2\cdot\texttt{SPR}[\texttt{c}]$}{
		$\vv{x}_c \sim \mathcal{U}(-h/2, +h/2)$\tcp*[r]{Random center around the origin}
		$\phi_c(\cdot) \leftarrow$ \circlelevelset{$\vv{x}_c$, $r_c$}\tcp*[r]{Create level-set function using \cref{eq:CircularLevelSetFunction}}
		let $\mathcal{B}_\Omega$ be the domain bounds ensuring that $\|\vv{x}-\vv{x}_c\|_\infty \leqslant r_c + 4h$ holds for any $\vv{x} \in \Omega$\;
		$\mathcal{G} \leftarrow$ \generategrid{\normalfont $\phi_c(\cdot)$, $\eta$, $\mathcal{B}_\Omega$, $2h\sqrt{2}$}\tcp*[r]{Discretizing $\Omega$ with an adaptive mesh}
		\BlankLine
		
		$\vv{\phi} \leftarrow$ \evaluate{$\mathcal{G}$, $\phi_c(\cdot)$}\tcp*[r]{Nodal level-set values}
		$\vv{\phi} \leftarrow$ \reinitialize{$\vv{\phi}$, $\nu$}\tcp*[r]{Redistancing with $\nu$ iterations (see \cref{eq:Reinitialization})}
		$\hat{N} \leftarrow$ \computenormals{$\mathcal{G}$, $\vv{\phi}$}\tcp*[r]{Nodal unit normal vectors and curvature (see \cref{eq:NormalAndCurvature})}
		$\vv{\kappa} \leftarrow$ \numcurvature{$\mathcal{G}$, $\vv{\phi}$, $\hat{N}$}\;
		\BlankLine
		
		$\mathcal{N} \leftarrow$ \getnodesnexttogamma{$\mathcal{G}$, $\vv{\phi}$}\;
		\ForEach{node $\mathcal{n} \in \mathcal{N}$}{
			\lIf{\normalfont $\mathcal{U}(0,1) \geqslant 1/\texttt{KeepEveryX}$}{skip node $\mathcal{n}$}
			\BlankLine
			
			$\mathcal{p} \leftarrow$ \collectfeatures{\normalfont $\mathcal{n}.\texttt{stencil}$, $\vv{\phi}$, $\hat{N}$}\tcp*[r]{Populate $\mathcal{p}$ (see \cref{fig:DataPacket.Original})}
			$\vv{x}_\mathcal{n}^\Gamma \leftarrow \mathcal{n}.\vv{x} - \vv{\phi}[\mathcal{n}]\hat{N}[\mathcal{n}] \quad \mathcal{p}.h\kappa \leftarrow h \cdot$\interpolate{$\mathcal{G}$, $\vv{\kappa}$, $\vv{x}_\mathcal{n}^\Gamma$}\tcp*[r]{See \cref{eq:ProjectionOntoGamma}}
			\BlankLine
			
			let $\vv{\xi}$ be the learning tuple $\left(\mathcal{p}, h\kappa^*\right)$ with inputs $\mathcal{p}$ and expected output $h\kappa^* := h\cdot\kappa^*$\;
			transform $\vv{\xi}$ so that $\mathcal{p}.h\kappa^*$ is negative\;
			rotate $\vv{\xi}$ so that the angle of $\mathcal{p}.\hat{\vv{n}}_{00}$ lies between 0 and $\pi/2$\tcp*[r]{See \cref{fig:DataPacket.Reoriented}}
			add $\vv{\xi}$ to $\mathcal{S}$\;
			\BlankLine
			
			let $\vv{\xi}'$ be the reflected tuple about the line $y=x+\beta$ going through $\mathcal{n}.\vv{x}$\tcp*[r]{See \cref{fig:DataPacket.Augmented}}
			add $\vv{\xi}'$ to $\mathcal{S}$\;
		}
	}	
	$\mathcal{D}_c \leftarrow \mathcal{D}_c\, \cup$ \randomsamples{\normalfont $\mathcal{S}$, $2\cdot\texttt{SPR}[\texttt{c}]$}\tcp*[r]{Grow $\mathcal{D}_c$ by drawing $\texttt{SPR}[\texttt{c}]$ samples from $\mathcal{S}$ if needed}
}
\BlankLine

\Return $\mathcal{D}_c$\;

\caption{\small $\mathcal{D}_c \leftarrow$ {\tt GenerateCircularDataSet(}$\eta$, $h\kappa_{\min}^*$, $h\kappa_{\max}^*$, {\tt CPH}, {\tt SPH2}, $\nu$, {\tt KeepEveryX)}: Generate a randomized data set with samples from circular-interface level-set functions.}
\label{alg:GenerateCircularDataSet}
\end{algorithm}

\begin{figure}[t]
	\centering
	\includegraphics[width=15cm]{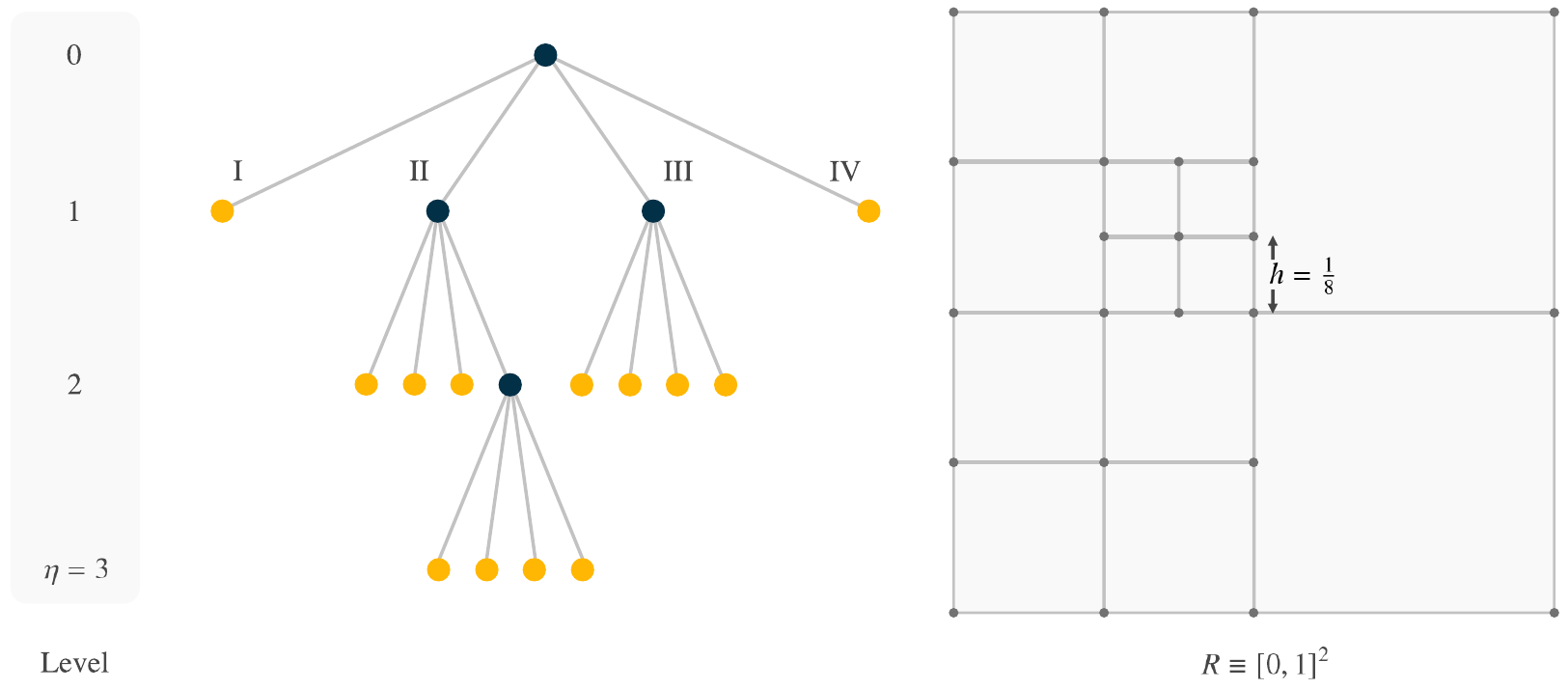}
	\caption{A unit-square quadtree with maximum level of refinement $\eta = 3$ and its corresponding subdivision (adapted from \cite{BKOS00}).  Quadrants (i.e., cells) are enumerated from I to IV in the standard Cartesian order.  Leaves appear in yellow in the center data structure.  (Color online.)}
	\label{fig:Quadtree}
\end{figure}

The formal parameters in \Cref{alg:GenerateCircularDataSet} include $h\kappa_{\min}^*$, $h\kappa_{\max}^*$, and the highest level of refinement to instantiate $\mathcal{G}$.  In addition, one must supply the number of steps for redistancing \cref{eq:CircularLevelSetFunction} and other constants that modulate randomization and interface dispersal.  In our case, $h\kappa_{\min}^* = 0.004 = h\kappa_{\min}^{low}$ and $h\kappa_{\max}^* = 2/3$, where $\kappa_{\max}^*$ is the maximal curvature associated with $r_c = 1.5h$.  Also, we have chosen adaptive Cartesian grids to discretize the computational domain \cite{Min;Gibou:07:A-second-order-accur, Mirzadeh;etal:16:Parallel-level-set}.  Under these schemes, $\eta \in \mathbb{N}^+$ represents the maximum number of subdivisions for a unit-square quadtree \cite{BKOS00}, so that $h = 2^{-\eta}$.  Quadtrees are rooted partitioning structures composed of cells that cover a rectangular region $R \subseteq \Omega$.  Each of these cells either has four children or is a leaf \cite{Strain1999}.  The quadtree building process thus comprises recursive cellular subdivisions governed by Min's extended Whitney decomposition criterion \cite{Min2004}.  In particular, we mark cell $C$ for refinement if the inequality

\begin{equation*}
\min_{\mathcal{v} \in \textrm{vertices}(C)} |\phi(\mathcal{v}.\vv{x})| \leqslant \textrm{Lip}(\phi)\cdot \texttt{diag}(C)
\end{equation*}
holds, where $\mathcal{v}$ is a vertex, $\textrm{Lip}(\phi)$ is the Lipschitz constant (equivalent to $1.2$), and $\texttt{diag}(C)$ is $C$'s diagonal.  Quadtrees are convenient spatial structures for solving FBPs efficiently in the level-set method \cite{Mirzadeh;etal:16:Parallel-level-set}.  The reason is that they help increase resolution only close to $\Gamma$, where accuracy is needed the most \cite{Strain1999}.  As an example, \cref{fig:Quadtree} shows a four-level unit-square quadtree, where $\eta = 3$ and $h = 2^{-3} = \frac{1}{8}$.

The first statements in \Cref{alg:GenerateCircularDataSet} set up the minimum and maximum radii and curvature values.  Then, the procedure defines how many distinct radii, $\texttt{NC}$, we will evaluate.  $\texttt{NC}$ is proportional to the $h$-normalized distance from $r_{\min}$ to $r_{\max}$, where the proportionality constant is the user-defined $\texttt{CPH} \in \mathbb{N}^+$---the number of circles per $h$.  For us, $\texttt{CPH} = 2$ has worked well across grid resolutions.  Next, the routine uses a heuristic to estimate the number of samples, $\texttt{AvgSPR}$, for the mean radius $\bar{r}$.  $\texttt{AvgSPR}$ is proportional to the $h^2$-normalized area difference of the circles with radii $\bar{r}$ and $\bar{r} - h$.  And, $\texttt{AvgSPR}$'s constant of proportionality is $\texttt{SPH2} \in \mathbb{N}^+$, which is the number of samples per $h^2$.  Similarly, we dampen $\texttt{AvgSPR}$ by a factor of $\texttt{KeepEveryX} \in \mathbb{N}^+$, where the latter acts as a probabilistic filter to minimize stencil overlaps later during sample extraction.  In this study, $\texttt{SPH2} = 5$ and $\texttt{KeepEveryX} = 4$ have produced good results.

Next, \Cref{alg:GenerateCircularDataSet} uses $\texttt{AvgSPR}$ to determine the sample distribution, $\texttt{SPR}$, for all unique radii.  \texttt{SPR} resembles a (discrete) trapezoidal density function with $\texttt{NC}$ steps, where the images of $r_{\min}$ (i.e., $\kappa_{\max}^*$) and $r_{\max}$ (i.e., $\kappa_{\min}^*$) are $\frac{3}{4}\texttt{AvgSPR}$ and $\frac{5}{4}\texttt{AvgSPR}$, respectively.  After setting \texttt{SPR}, we first spread curvature uniformly with $\texttt{NC}$ values from $\kappa_{\max}^*$ to $\kappa_{\min}^*$ into $\texttt{TgtK}$.  Then, we enter the data production loop, where we accumulate learning tuples for each $r_c \in [r_{\min}, r_{\max}]$ into $\mathcal{S}$ until $|\mathcal{S}|$ is twice $r_c$'s image in $\texttt{SPR}$.  Note that the last factor of 2 is in place to account for the data augmentation scheme introduced in \Cref{alg:MLCurvature}.

Once inside the sampling loop, we draw a random center $\vv{x}_c$ within the $h$-by-$h$ region around the origin to define $\phi_c(\cdot)$ (see \cref{eq:CircularLevelSetFunction}).  Then, we discretize a sufficiently large $\Omega$ with a grid $\mathcal{G}$ by using quadtrees leading to the desired $h$.  Also, we enforce a uniform band of half-width $2h\sqrt{2}$ around $\Gamma$ to ensure valid interface stencils.  Afterward, we use $\mathcal{G}$ to evaluate $\phi_c(\cdot)$ into $\vv{\phi}$ and, subsequently, solve \cref{eq:Reinitialization} with $\nu = 10$ iterations.  Next, we compute the nodal unit normal vectors and curvatures into $\hat{N}$ and $\vv{\kappa}$ to prepare for sampling in the following step.

The innermost loop in \Cref{alg:GenerateCircularDataSet} iterates over each interface node $\mathcal{n}$ and adds its samples into $\mathcal{S}$ with a probability of $1/\texttt{KeepEveryX}$.  For every sampled vertex, we populate its data-packet features in \cref{eq:DataPacket} and construct two learning tuples:

\begin{equation}
\vv{\xi} = \left(\mathcal{p}, h\kappa^*\right) \quad \textrm{and} \quad 
\vv{\xi}' = \left(\mathcal{p}', h\kappa^*\right),
\label{eq:TrainSamples}
\end{equation}
where $h\kappa^*$ is the target dimensionless curvature and $-h\kappa_{\max}^* \leqslant h\kappa^* \leqslant -h\kappa_{\min}^*$.  Similar to \Cref{alg:MLCurvature}, $\mathcal{p}$ is a reoriented data packet (e.g., \cref{fig:DataPacket.Reoriented}), while $\mathcal{p}'$ is its reflected version (e.g.,  \cref{fig:DataPacket.Augmented}).  Upon collecting $|\mathcal{S}| \geqslant 2\cdot\texttt{SPR}[\texttt{c}]$ samples for the $\texttt{c}$th radius, we remove the excess randomly and conclude by adding $\texttt{SPR}[\texttt{c}]$ learning pairs into $\mathcal{D}_c$.


\colorsubsubsection{Sinusoidal-interface data-set construction}
\label{subsubsec:SinusoidalInterfaceDataSetConstruction}

Sinusoidal curves are the second class of elementary interfaces considered when assembling $\mathcal{D}$.  As in \cite{Larios;Gibou;HybridCurvature;2021}, we construct the corresponding level-set functions using sine waves of the form

\begin{equation}
f_s(t) = A\sin(\omega t),
\label{eq:SineWave}
\end{equation}
where $A$ and $\omega$ are the desired amplitude and frequency.  \Cref{eq:SineWave} gives rise to the level-set function

\begin{equation}
\phi_s(x, y) = \left\{\begin{array}{ll}
	-d_s(x, y) & \textrm{if } f_s(x) < y, \\
	 0 & \textrm{if } f_s(x) = y, \\
	+d_s(x, y) & \textrm{if } y < f_s(x),
\end{array}\right.
\label{eq:SineWaveLevelSetFunction}
\end{equation}
where $d_s(x, y)$ is the shortest distance from $(x, y) \in \Omega$ to $f_s(t)$.  The simple expression in \cref{eq:SineWave} also yields a straightforward relation,

\begin{equation}
\kappa_s(t) = -\frac{A\omega^2 \sin(\omega t)}{\left( 1 + A^2 \omega^2 \cos^2(\omega t) \right)^{3/2}},
\label{eq:SineWaveCurvature}
\end{equation}
for computing curvature \cite{Swokowski88} along $f_s(t)$.

In general, varying just $A$ and $\omega$ could diversify the learning samples in the associated data set $\mathcal{D}_s$.  However, further pattern variability is possible if we introduce affine transformations \cite{CGUsingOpenGL01}.  Their purpose is to alter $f_s(t)$'s canonical frame $\mathcal{C_s}$, as depicted in \cref{fig:TrainingInterfaces.Sinusoidal}.  For example, if $T(\vv{x}_s)$ and $R(\theta_s)$ denote a translation and a rotation applied to $\mathcal{C_s}$, $ _0M_s = T(\vv{x}_s)R(\theta_s)$ maps points and vectors in homogeneous coordinates\footnote{In homogenous coordinates, $[x,y,1]^T$ is a point, while $[x,y,0]^T$ is a vector.  $T(\cdot)$ and $R(\cdot)$ are thus $3$-by-$3$ matrices in $\mathbb{R}^2$.} from $\mathcal{C}_s$ into the world coordinate system $\mathcal{C}_0$.  Here, $\vv{x}_s = (x_s, y_s)$ is a prescribed shift from the origin, and $\theta_s$ is $\mathcal{C}_s$'s tilt with respect to $\mathcal{C}_0$'s horizontal.  Similarly, ${}_sM_0 ={}_0M_s^{-1} = R^T(\theta_s)T(-\vv{x}_s)$ allows us to express a point or vector in $\mathcal{C}_0$ in terms of $\mathcal{C}_s$.  Consequently, by tweaking $\vv{x}_s$, $\theta_s$, $A$, and $\omega$, we can span a wide range of level-set configurations without forgoing the simplicity of \cref{eq:SineWaveCurvature} to compute $h\kappa^*$.  More specifically, if $\vv{x} \in \Omega$ is given in the representation of $\mathcal{C}_0$, we can calculate its level-set value as $\phi_s( \cancel{\llbracket}_sM_0 \llbracket\vv{x}\rrbracket\cancel{\rrbracket})$.  Likewise, if $\vv{x}_\mathcal{n}^\perp$ is an interface node's projection onto $\Gamma$, we can retrieve $h\kappa^*$ by evaluating $\kappa_s(t^*)$, where $t^*$ is the first component of $\llbracket {}_sM_0\vv{x}_\mathcal{n}^\perp\rrbracket$.  In these cases, $\llbracket\cdot\rrbracket$ expands its argument into homogenous coordinates, and $\cancel{\llbracket}\cdot\cancel{\rrbracket}$ removes the homogenization suffix.

Furthermore, solving \cref{eq:SineWaveLevelSetFunction,eq:SineWaveCurvature} implies recovering $\tilde{t} = \arg\min_t \left\Vert\tilde{\vv{x}} - [t,\, f_s(t)]^T\right\Vert_2$ for any (transformed) grid point $\tilde{\vv{x}} = \cancel{\llbracket}_sM_0 \llbracket\vv{x}\rrbracket\cancel{\rrbracket}$.  Here, we compute such a $\tilde{t}$ using bisection followed by Newton--Raphson's root-finding \cite{Heath;SciComput;2018}.  Then, we reinitialize the nodal $\phi_s(\cdot)$ values as a way of perturbing the otherwise exact signed distance field.  By doing so, we make sinusoidal level-set patterns compatible with those in $\mathcal{D}_c$.  Also, it helps us mimic the settings of real-world simulations, where reinitialization is the employed mechanism to estimate distances to $\Gamma$.


\begin{algorithm}[!t]
\small
\SetAlgoLined
\SetKwFunction{linspace}{Linspace}
\SetKwFunction{sinelevelset}{SineLevelSet}
\SetKwFunction{generategrid}{GenerateGrid}
\SetKwFunction{evaluate}{Evaluate}
\SetKwFunction{reinitialize}{Reinitialize}
\SetKwFunction{computenormals}{ComputeNormals}
\SetKwFunction{numcurvature}{NumCurvature}
\SetKwFunction{getnodesnexttogamma}{GetNodesNextToGamma}
\SetKwFunction{easeoff}{EaseOff}
\SetKwFunction{collectfeatures}{CollectFeatures}
\SetKwFunction{randomsamples}{RandomSamples}

\KwIn{$\eta$: maximum level of refinement per unit-square quadtree; $h\kappa_{\min}^*$ and $h\kappa_{\max}^*$: minimum and maximum target $|h\kappa^*|$; {\tt EaseOffMidMaxHKPr}: easing-off probability threshold to always keep samples based on $|h\kappa^*|$; {\tt NA}: number of amplitudes; {\tt NT}: number of discrete rotation angles; $\nu$: number of iterations for level-set reinitialization.}
\KwResult{$\mathcal{D}_s$: data set of sinusoidal-interface samples.}
\BlankLine

$h \leftarrow 2^{-\eta}$\tcp*[r]{Mesh size}
$\kappa_{\min}^* \leftarrow h\kappa_{\min}^*/h; \quad \kappa_{\max}^* \leftarrow h\kappa_{\max}^*/h; \quad A_{\min} \leftarrow 4/\kappa_{\max}^*; \quad A_{\max} \leftarrow 1/\kappa_{\min}^*$\tcp*[r]{Curvature and amplitude bounds}
$r_{sam} \leftarrow 2A_{\max}$\tcp*[r]{Sampling radius}
$h\kappa_{\max}^{low} \leftarrow \frac{1}{2}h\kappa_{\max}^*; \quad h\kappa_{\max}^{up} \leftarrow h\kappa_{\max}^*$\tcp*[r]{Bounds for maximum curvature at the crests}
$h\bar{\kappa}_{\max}^* \leftarrow \frac{1}{2}\left(h\kappa_{\max}^{low} + h\kappa_{\max}^{up}\right)$\tcp*[r]{Mean target maximum $h\kappa^*$}
\BlankLine
			
\tcp{Data generation loop; first defining the $A$ and $\omega$ parameters}
$\mathcal{D}_s \leftarrow \varnothing$\;
\ForEach{amplitude $A \in$ \linspace{\normalfont $A_{\min}$, $A_{\max}$, $\texttt{NA}$}}{
	$\omega_{\min} \leftarrow \sqrt{h\kappa_{\max}^{low}/(h \cdot A)}; \quad \omega_{\max} \leftarrow \sqrt{h\kappa_{\max}^{up}/(h \cdot A)}$\tcp*[r]{Frequency range to ensure $|\kappa_{\max}^*| \in \left[\kappa_{\max}^{low}, \kappa_{\max}^{up}\right]$}
	$\omega_{d} \leftarrow \frac{\pi}{2}\left(\omega_{\min}^{-1} - \omega_{\max}^{-1}\right)$\tcp*[r]{Distance between first positive crests with $\omega_{\min}$ and $\omega_{\max}$ freqs.}
	$\texttt{NF} \leftarrow \left\lceil \omega_d / h \right\rceil + 1$\tcp*[r]{Number of frequency steps}
	\BlankLine
	
	\ForEach{frequency $\omega \in$ \linspace{\normalfont $\omega_{\min}$, $\omega_{\max}$, $\texttt{NF}$}}{
		\BlankLine
		
		$\mathcal{S} \leftarrow \varnothing$\tcp*[r]{Collect samples for interface $A\sin(\omega t)$ into $\mathcal{S}$}
		\ForEach{angle $\theta_s \in$ \linspace{\normalfont $-\pi/2$, $+\pi/2$, $\texttt{NT}$}, excluding right end point}{
			$\vv{x}_s \sim \mathcal{U}(-h/2, +h/2)$\tcp*[r]{Random displacement around the origin}
			$\phi_s(\cdot) \leftarrow$ \sinelevelset{$A$, $\omega$, $\vv{x}_s$, $\theta_s$}\tcp*[r]{Create level-set function using \cref{eq:SineWaveLevelSetFunction}}
			let $\mathcal{B}_\Omega$ be the domain bounds ensuring that $\|\vv{x}-\vv{x}_s\|_\infty \leqslant r_{sam} + 4h$ holds for any $\vv{x} \in \Omega$\;
			$\mathcal{G} \leftarrow$ \generategrid{\normalfont $\phi_s(\cdot)$, $\eta$, $\mathcal{B}_\Omega$, $4h\sqrt{2}$}\tcp*[r]{Discretizing $\Omega$ with an adaptive mesh}
			\BlankLine
			
			$\vv{\phi} \leftarrow$ \evaluate{$\mathcal{G}$, $\phi_s(\cdot)$}\tcp*[r]{Nodal level-set values}
			$\vv{\phi} \leftarrow$ \reinitialize{$\vv{\phi}$, $\nu$}\tcp*[r]{Redistancing with $\nu$ iterations (see \cref{eq:Reinitialization})}
			$\hat{N} \leftarrow$ \computenormals{$\mathcal{G}$, $\vv{\phi}$}\tcp*[r]{Nodal unit normal vectors and curvature (see \cref{eq:NormalAndCurvature})}
			$\vv{\kappa} \leftarrow$ \numcurvature{$\mathcal{G}$, $\vv{\phi}$, $\hat{N}$}\;
			\BlankLine
			
			$\mathcal{N} \leftarrow$ \getnodesnexttogamma{$\mathcal{G}$, $\vv{\phi}$}\;
			\ForEach{node $\mathcal{n} \in \mathcal{N}$}{
				compute the target $h\kappa^*$ at the nearest point, $\vv{x}_\mathcal{n}^\perp$, to $\mathcal{n}.\vv{x}$ on the sinusoidal interface\tcp*[r]{See \cref{eq:SineWaveCurvature}}
				\lIf{\normalfont $\|\mathcal{n}.\vv{x} - \vv{x}_s\|_2^2 > r_{sam}^2$ \textbf{or} $|h\kappa^*| < h\kappa_{\min}^*$}{skip node $\mathcal{n}$}
				\BlankLine

				\tcp{Prob. sampling: $\mathrm{Pr}(|h\kappa^*|\!\geqslant\!h\bar{\kappa}_{\max}^*) = \texttt{EaseOffMidMaxHKPr}$ and $\mathrm{Pr}(|h\kappa^*|\!=\!h\kappa_{\min}^*) = 0.01$}				
				$q \leftarrow \min\left( 1, (|h\kappa^*| - h\kappa_{\min}^*)/(h\bar{\kappa}_{\max}^* - h\kappa_{\min}^*) \right)$\;
				\If{$\mathcal{U}(0,1) \leqslant$ \easeoff{\normalfont $0.01$, $\texttt{EaseOffMidMaxHKPr}$, $q$}}{
					$\mathcal{p} \leftarrow$ \collectfeatures{\normalfont $\mathcal{n}.\texttt{stencil}$, $\vv{\phi}$, $\hat{N}$}\tcp*[r]{Populate $\mathcal{p}$ (see \cref{fig:DataPacket.Original})}
					$\vv{x}_\mathcal{n}^\Gamma \leftarrow \mathcal{n}.\vv{x} - \vv{\phi}[\mathcal{n}]\hat{N}[\mathcal{n}]; \quad \mathcal{p}.h\kappa \leftarrow h \cdot$\interpolate{$\mathcal{G}$, $\vv{\kappa}$, $\vv{x}_\mathcal{n}^\Gamma$}\tcp*[r]{See \cref{eq:ProjectionOntoGamma}}
					\BlankLine
			
					let $\vv{\xi}$ be the learning tuple $\left(\mathcal{p}, h\kappa^*\right)$ with inputs $\mathcal{p}$ and expected output $h\kappa^*$\;
					transform $\vv{\xi}$ so that $\mathcal{p}.h\kappa^*$ is negative\;
					rotate $\vv{\xi}$ so that the angle of $\mathcal{p}.\hat{\vv{n}}_{00}$ lies between 0 and $\pi/2$\tcp*[r]{See \cref{fig:DataPacket.Reoriented}}
					add $\vv{\xi}$ to $\mathcal{S}$\;
					\BlankLine
			
					let $\vv{\xi}'$ be the reflected tuple about the line $y=x+\beta$ going through $\mathcal{n}.\vv{x}$\tcp*[r]{See \cref{fig:DataPacket.Augmented}}
					add $\vv{\xi}'$ to $\mathcal{S}$\;
				}
			}
		}
		$\mathcal{D}_s \leftarrow \mathcal{D}_s\, \cup \mathcal{S}$\;
	}
}
\BlankLine

\Return $\mathcal{D}_s$\;

\caption{\small $\mathcal{D}_s \leftarrow$ {\tt GenerateSinusoidalDataSet(}$\eta$, $h\kappa_{\min}^*$, $h\kappa_{\max}^*$, {\tt EaseOffMidMaxHKPr}, {\tt NA}, {\tt NT}, $\nu${\tt )}: Generate a randomized data set with samples from sinusoidal-interface level-set functions.}
\label{alg:GenerateSinusoidalDataSet}
\end{algorithm}

\Cref{alg:GenerateSinusoidalDataSet} describes the steps to construct $\mathcal{D}_s$.  Essentially, the combinatorial procedure varies $A$, $\omega$, and $\theta_s$ to collect samples from rotated and randomly shifted sine waves.  Like \Cref{alg:GenerateCircularDataSet}, its formal parameters include the maximum level of refinement, the target dimensionless-curvature bounds, and the number of iterations to reinitialize \cref{eq:SineWaveLevelSetFunction}.  In addition, the routine requires the $\texttt{EaseOffMidMaxHKPr}$ probability threshold and other constants that control the number of amplitudes ($\texttt{NA}$) and tilts ($\texttt{NT}$).  Again, we have chosen $h\kappa_{\min}^* = 0.004 = h\kappa_{\min}^{low}$, $h\kappa_{\max}^* = 2/3$, and $\nu = 10$.  Similarly, $\texttt{NA} = 34$ and $\texttt{NT} = 38$ have worked well for all $\eta$.

The first few statements in \Cref{alg:GenerateSinusoidalDataSet} initialize the extremum curvatures and amplitudes and the mesh size.  Also, it specifies $r_{sam}$, which denotes the radius of the sampling circular region $R_c \subset \Omega$ centered at $\vv{x}_s$.  Unlike \cite{Larios;Gibou;HybridCurvature;2021}, $R_c$ helps prevent data explosion since it shrinks and grows proportionally to $\eta$.  Afterward, we define the closed interval $[h\kappa_{\max}^{low},\, h\kappa_{\max}^{up}]$, which controls steepness along the crests.  Next, we compute the interval's midpoint $h\bar{\kappa}_{\max}^*$ in anticipation of the upcoming probabilistic sampling.

\Cref{alg:GenerateSinusoidalDataSet} then enters the data generation loop, where the sine wave amplitude varies uniformly from $A_{\min}$ to $A_{\max}$.  For each such amplitude $A$, we first compute the frequencies associated with $\kappa_{\max}^{low}$ and $\kappa_{\max}^{up}$.  And, using these values, we calculate the $t$-parameter distance, $\omega_d$, between the first positive crests with frequencies $\omega_{\min}$ and $\omega_{\max}$.  Once $\omega_d$ is available, we determine $\texttt{NF}$, which denotes the number of discrete frequency steps in the following loop.

For every amplitude-frequency pair, we spawn affine-transformed $f_s(t)$ interfaces and gather their samples into $\mathcal{S}$.  To do so, we iteratively select $\texttt{NT}$ equally spaced $\theta_s$ angles, and, for each $\theta_s$, we draw a random displacement $\vv{x}_s \in [-h/2,\, +h/2]^2$ (see \cref{fig:TrainingInterfaces.Sinusoidal}).  The subsequent tasks are then parallel to the statements in \Cref{alg:GenerateCircularDataSet}.  These include: instantiating $\phi_s(\cdot)$, discretizing $\Omega$ with an adaptive grid $\mathcal{G}$, enforcing a $4h\sqrt{2}$-half-width uniform band around $\Gamma$, evaluating $\phi_s(\cdot)$ into the nodal vector $\vv{\phi}$, reinitializing \cref{eq:SineWaveLevelSetFunction}, and calculating the nodal unit normal vectors and curvatures into $\hat{N}$ and $\vv{\kappa}$.  Lastly, \Cref{alg:GenerateSinusoidalDataSet} gathers the interface nodes in $\mathcal{N}$ before proceeding to the sampling stage.

As we iterate over each $\mathcal{n} \in \mathcal{N}$, we skip grid points outside of $R_c$ or whose expected $|h\kappa^*|$ falls below $h\kappa_{\min}^*$.  Next, we apply probabilistic filtering based on $|h\kappa^*|$ and an ease-off density function \cite{ComputerAnimation08}.  The latter is designed to ensure that $\textrm{Pr}(|h\kappa^*|\!=\!h\kappa_{\min}^*) = 0.01$ and $\textrm{Pr}(|h\kappa^*|\!\geqslant\!h\bar{\kappa}_{\max}^*) = \texttt{EaseOffMidMaxHKPr}$.  In our case, setting $\texttt{EaseOffMidMaxHKPr} = 0.4$ has produced reasonably large sinusoidal data sets.  Finally, we populate the data packets for the selected vertices and insert their $\vv{\xi}$ and $\vv{\xi}'$ samples into $\mathcal{S}$ (see \cref{eq:TrainSamples}).  Eventually, we grow $\mathcal{D}_s$ with the contents of $\mathcal{S}$ after depleting all tilts for each frequency.  Then, we return the resulting $\mathcal{D}_s$ to the calling function upon completing the last iteration for $A_{\max}$.

Despite our efforts to contain shallow-curvature-sample over-representation in \Cref{alg:GenerateSinusoidalDataSet}, $\mathcal{D}_s$'s $h\kappa^*$ distribution is often heavily skewed towards 0.  To remedy this situation, we have resorted to histogram-based post-balancing.  In this technique, we first distribute $\mathcal{D}_s$ into $100$ intervals.  Then, we randomly subsample overpopulated bins until the condition $|b| \leqslant \tfrac{2}{3}m_\mathcal{H}$ holds for each bin $b$ in histogram $\mathcal{H}$, where $m_\mathcal{H}$ is $\mathcal{H}$'s median.  After balancing $\mathcal{D}_s$, we combine it with $\mathcal{D}_c$ into $\mathcal{D}$ and use it to train $\mathcal{F}_\kappa(\cdot)$, as described in the following section.


\colorsubsection{Technical aspects}
\label{subsec:TechnicalAspects}

In this study, we have jump-started $\mathcal{F}_\kappa(\cdot)$'s optimization with our explored configurations in  \cite{LALariosFGibou;LSCurvatureML;2021, Larios;Gibou;HybridCurvature;2021}.  Here, however, we introduce a different loss function, regularization, and dimensionality reduction to speed up training and improve generalization. As seen in \Cref{sec:Results}, these upgrades have decreased outlying artifacts while increasing curvature accuracy around poorly resolved regions.


\begin{algorithm}[!t]
\small
\SetAlgoLined
\SetKwFunction{getstats}{GetStats}

\KwIn{$\mathcal{p}$: data packet to preprocess; $h$: mesh size.}
\KwResult{$\vv{p}$: transformed data packet into a vector form.}
\BlankLine

\tcp{In-place $h$-normalization of level-set values}
Divide each $\mathcal{p}.\phi_{ij}$ by $h$, for $i, j \in \{m, 0, p\}$\;
\BlankLine

\tcp{In-place z-scoring}
$\mathcal{Q} \leftarrow$ \getstats{}\tcp*[r]{Retrieve training stats}
let $\vv{q}$ be the vector of features in $\mathcal{p}$\;
let $\vv{\mu}$ and $\vv{\sigma}$ be the vectors with per-feature mean and std. deviation in $\mathcal{Q}.${\tt StdScaler}\;
$\vv{z} \leftarrow (\vv{q} - \vv{\mu}) \oslash \vv{\sigma}$\footnotemark\;
\BlankLine

\tcp{Dimensionality reduction and whitening}
let $V$ and $\vv{\sigma}$ be the components' matrix and explained std. deviations' vector in $\mathcal{Q}.${\tt PCA}\;
$\vv{p} \leftarrow (V^T\vv{z}) \oslash \vv{\sigma}$\tcp*[r]{PCA transformation}
\BlankLine

\Return $\vv{p}$\;

\caption{\small $\vv{p} \leftarrow$ {\tt Preprocess(}$\mathcal{p}$, $h${\tt )}: Preprocess a data packet $\mathcal{p}$ using training statistics.}
\label{alg:Preprocess}
\end{algorithm}
\footnotetext{Element-wise division.}

First, we remark that \Crefrange{alg:MLCurvature}{alg:GenerateSinusoidalDataSet} have been realized within our C++ implementation of the adaptive-grid parallel-level-set framework of \cite{Mirzadeh;etal:16:Parallel-level-set}.  In addition, we have incorporated Boost's functionalities \cite{Boost;2019} to solve the minimization problem stated in \Cref{subsubsec:SinusoidalInterfaceDataSetConstruction}.  Later, we will discuss other technical features to materialize the hybrid solver in \Cref{alg:MLCurvature}.  

To train $\mathcal{F}_\kappa(\cdot)$, we use TensorFlow \cite{Tensorflow15} and Keras \cite{Keras15} in Python.  Also, we have partitioned $\mathcal{D}$ into training, testing, and validation subsets to guide model selection and prevent overfitting \cite{A18, Mehta19}.  Our partitioning approach emulates a tenfold cross-validation scheme based on the {\tt StratifiedKFold} class from SciKit-Learn \cite{scikit-learn11}.  {\tt StratifiedKFold}, in particular, allows us to assemble learning subsets with the same class distribution as $\mathcal{D}$.  To create such labels, we begin by invoking Panda's {\tt cut()} function \cite{Pandas2010}.  The latter has helped us group samples into one hundred categories according to their $h\kappa^*$ values.  Next, by calling {\tt StratifiedKFold}'s {\tt split()} method three times with $K=10$, we have populated the training, testing, and validation subsets with 70\%, 10\%, and 10\% of the tuples in $\mathcal{D}$.  As in \cite{Larios;Gibou;ECNetSemiLagrangian;2021}, we discard the remaining 10\% while avoiding an exhaustive $K$-fold cross-validation for efficiency.  In the traditional $K$-fold cross-validation method, one usually divides the labeled data into $K$ equal segments.  Then, one of the $K$ segments is reserved for testing, and we use the remaining $(K-1)$ portions for training.  This process is subsequently repeated $K$ times by choosing each segment as the testing subset and reporting their average accuracy at the end \cite{A18}.  Such an approach is often too expensive to train neural networks with large data sets like ours.  For this reason, we have preferred running a couple of training instances with the same settings and random kernel initializations, given that $|\mathcal{D}| \approx \eten{1.4}{6}$ for all $\eta$.

The second critical aspect to discuss is the {\tt Preprocess()} module.  \Cref{alg:Preprocess} lists the preprocessing operations performed on reoriented data packets.  The goal of the {\tt Preprocess()} subroutine is to use training statistics to transform the raw attributes in \cref{eq:DataPacket} into an amenable vector representation for our neural estimator \cite{scikit-learn11}.  To this end, the procedure first $h$-normalizes the stencil level-set values.  Then, it employs the feature-wise means and variances stored in $\mathcal{Q}$ to standardize\footnote{Standardization or \textit{z-scoring} transforms a vector $\vv{\psi}$ into $\vv{\psi}'$, which has mean $0$ and variance $1$.} the vectorized $\mathcal{p}$ into $\vv{z} \in \mathbb{R}^{28}$.  Lastly, \Cref{alg:Preprocess} performs a principal component analysis (PCA) and whitening transformation on $\vv{z}$ to produce the output vector $\vv{p} \in \mathbb{R}^{m_\iota}$, where $m_\iota \in \{18, 20\}$ (see \cref{tbl:results.train.config}).  For us, PCA\footnote{If $D \in \mathbb{R}^{n \times p}$ is a data set, and $M \in \mathbb{R}^{n \times p}$ contains $D$'s column-wise mean $p$-vector $\vv{\mu}$ stacked $n$ times, then $C = \frac{1}{n-1}(D-M)^T(D-M)$ is $D$'s $p$-by-$p$ covariance matrix.  Further, if $\vv{x} = \vv{d} - \vv{\mu}$ is a centered, feature $p$-vector, and $C = U\Sigma V^T$ is $C$'s SVD decomposition, then $\vv{x}_k = V_k^T\vv{x}$ is the PCA-transformed $k$-vector, where $V_k$ contains $V$'s first $k$ columns.  Lastly, \textit{whitening} involves normalizing $\vv{x}_k$ so that $\vv{x}_k' = \Sigma_k^{-1/2} \vv{x}_k$, where $\Sigma_k$ is $\Sigma$'s $k$-order leading principal sub-matrix.} entails a change of basis, where the new axes are better data descriptors defined by the first $m_\iota$ singular vectors of the training subset's correlation matrix\footnote{The correlation matrix is equivalent to the covariance matrix of the standardized data \cite{Parker;CS170A;2016}.}.  Likewise, we have introduced the supplementary whitening step to provide equal importance to the different features in $\vv{p}$ \cite{A18}.  In this regard, prior research has proven that PCA and whitening are beneficial for learning convergence \cite{LeCun;EfficientBackProp;98, LALariosFGibou;LSCurvatureML;2021, Larios;Gibou;HybridCurvature;2021, Larios;Gibou;ECNetSemiLagrangian;2021}.  Also, data scientists have shown that combining PCA with z-scoring prevents variance from directly reflecting the scale of the data \cite{Parker;CS170A;2016}.  Motivated by such qualities, here we have resorted to SciKit-Learn's {\tt PCA} class to find $\mathcal{D}$'s principal components.  This {\tt PCA} object alongside $\mathcal{Q}$'s {\tt StdScaler} is precisely one of the {\tt JSON}-formatted exports we have made available to facilitate future deployments.

The building blocks available in TensorFlow and Keras have made it easy to design and train our error quantifier.  As portrayed in \cref{fig:ECNet}, our model digests a two-part input tensor with $m_\iota + 1$ components, where $m_\iota$ is resolution-dependent.  Each of the first four hidden layers contains $N_h^i$ ReLU neurons with incoming trainable connections, for $i = 1,2,3,4$.  Similarly, the weights in the last hidden-layer unit are updatable, but its activation function is linear to yield an appropriate $\bar{\varepsilon}$ term.  At last, we have included a single neuron that merely adds $h\kappa$ and $\bar{\varepsilon}$ to emit the error-corrected $h\kappa_\mathcal{F}$ output.

\begin{figure}[!t]
	\centering
	\begin{subfigure}[b]{0.49\textwidth}
		\includegraphics[width=\textwidth]{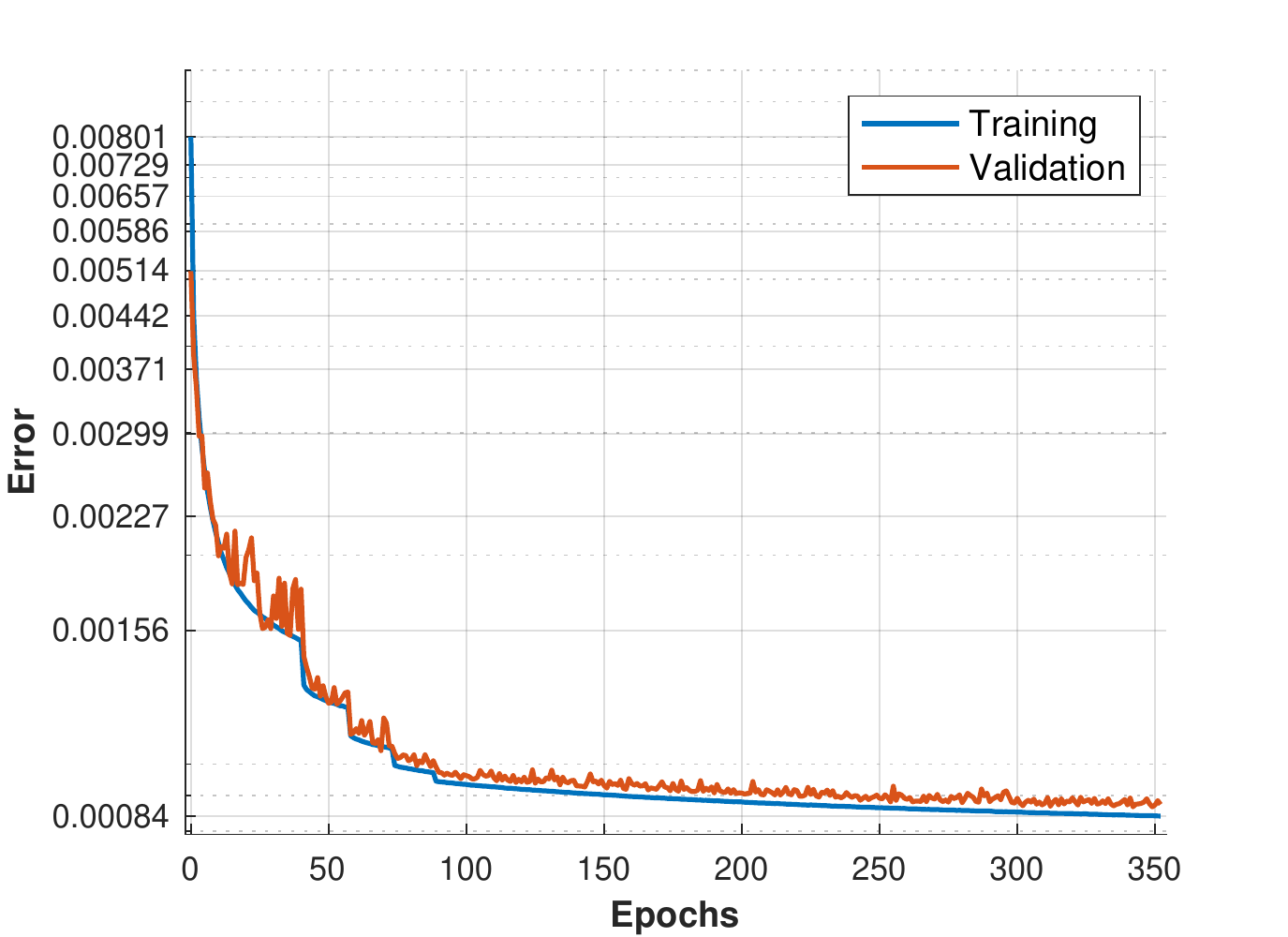}
        \caption{\footnotesize RMSE}
        \label{fig:learningcurves.rmse}
    \end{subfigure}
    ~
	\begin{subfigure}[b]{0.49\textwidth}
		\includegraphics[width=\textwidth]{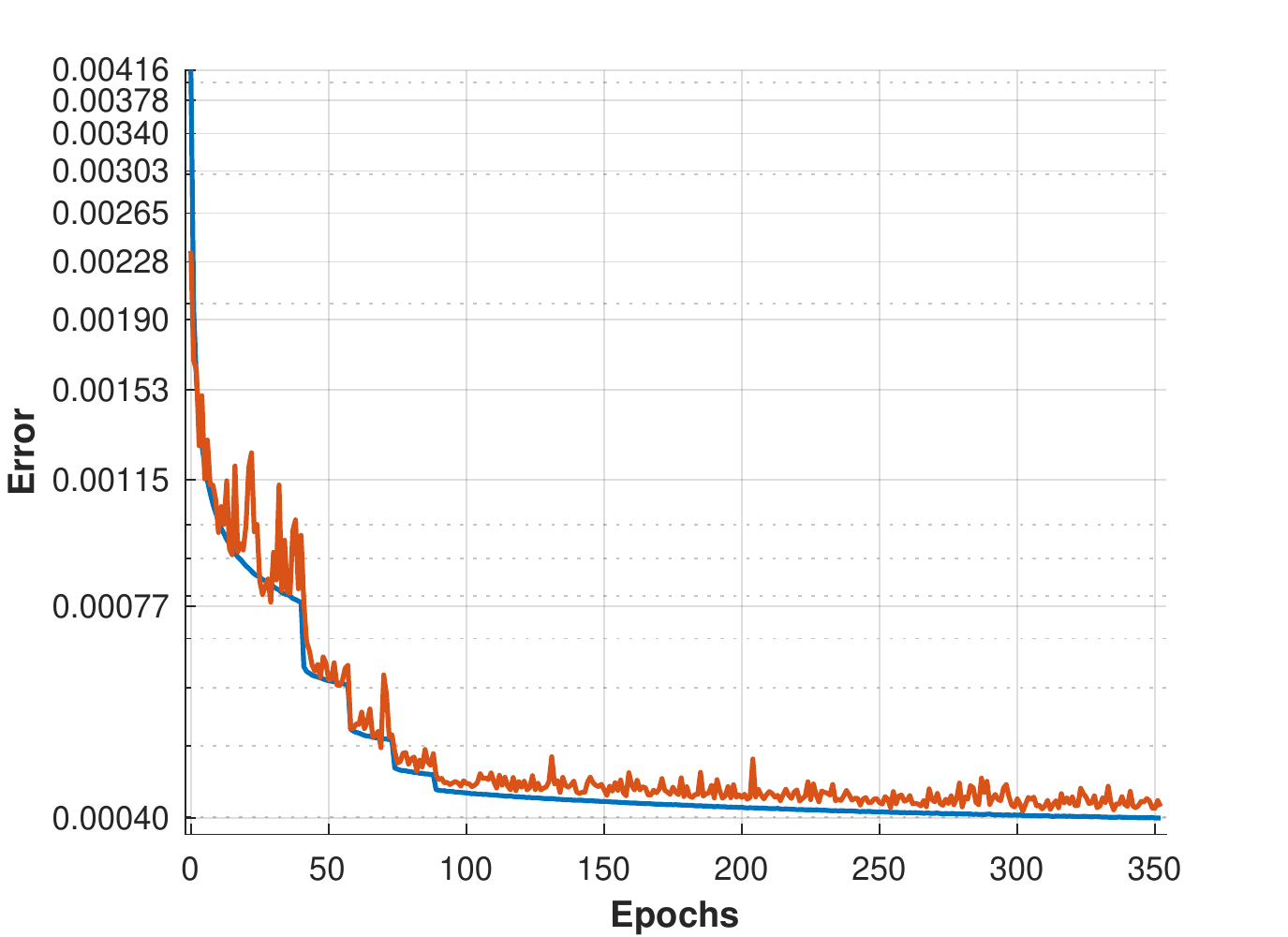}
		\caption{\footnotesize MAE}
		\label{fig:learningcurves.mae}
	\end{subfigure}
	   
	\caption{Example (a) RMSE and (b) MAE training and validation curves for $\eta = 7$.  Optimization results and architectural parameters appear in \Cref{subsec:SteepCurvatureFlowerShapedInterfaces}.  (Color online.)}
	\label{fig:learningcurves}
\end{figure} 

The routine for fitting $\mathcal{F}_\kappa(\cdot)$ to $\mathcal{D}$ incorporates most of the elements introduced in \cite{Larios;Gibou;ECNetSemiLagrangian;2021}.  More precisely, we have adopted the Adam optimizer and the root mean squared error loss (RMSE) to update $\mathcal{F}_\kappa(\cdot)$'s weights.  Every such update has occurred after observing 64-sample batches, with one thousand epochs as a limit.  Besides the RMSE, we monitor the mean absolute error (MAE) and the maximum absolute error (MaxAE) for validation and model selection.  In particular, we use the validation subset's MAE to stabilize backpropagation by halving the learning rate from $\eten{1.5}{-4}$ down to $\eten{1}{-5}$ whenever the error does not improve for fifteen epochs.  Also, we have integrated {\tt L2} kernel regularization and early stopping to promote network resilience.  To realize early stopping, we use a callback function that triggers after fifty epochs of learning stagnation.  Likewise, the layer-wise construction approach in Keras has allowed us to specify {\tt L2} factors without fiddling directly with the loss function.  Lastly, \cref{fig:learningcurves} illustrates a pair of typical training and validation curves.  As seen in these charts, {\tt L2} regularization favors generalization by closing the gap between the learning curves while working as a tunable ``forgetting'' mechanism \cite{A18}.

We close this section with a few comments about implementing \Cref{alg:MLCurvature,alg:Preprocess} in C++.  Analogous to \cite{Larios;Gibou;ECNetSemiLagrangian;2021}, our curvature solver leverages several third-party tools to decode, import, and speed up neural inference.  Among these tools, we have included Lohmann's {\tt json} library \cite{json;2021} and a couple of methods from Hermann's {\tt frugally-deep} framework \cite{frugally-deep;2021}.  These functionalities have helped us read in $\mathcal{F}_\kappa(\cdot)$'s weights (up to the last hidden layer in \cref{fig:ECNet}) and load the training statistics in $\mathcal{Q}$.  In addition, we have carried out neural inference as groups of 32-bit {\tt sgemm}\footnote{A \textbf{s}ingle-precision \textbf{ge}neral \textbf{m}atrix-\textbf{m}atrix multiplication of the form $C = \alpha AB+\beta C$, where $\alpha = 1$ and $\beta = 0$.} operations in OpenBLAS \cite{openblas;2021}.  The latter has provided an acceleration factor above sixty\footnote{Using C++14 compiling optimization enabled via {\tt -O2 -O3 -march=native}.} when comparing batch versus single predictions in {\tt frugally-deep}.


\FloatBarrier
\colorsection{Results}
\label{sec:Results}

Next, we assess our hybrid solver's ability to estimate curvature on irregular interfaces at various resolutions.  To this end, we have selected the bivariate function

\begin{equation}
\phi_{rose}(\vv{x}) = \|\vv{x}\|_2 - a \cos\left( p\theta \right) - b,
\label{eq:PolarRoseLevelSetFunction}
\end{equation} 
whose zero-isocontour delineates the flower-shaped interface

\begin{equation}
\gamma(\theta) = a\cos\left( p\theta \right) + b,
\label{eq:PolarRoseInterface}
\end{equation}
where $a$, $b$, and $p$ are shape parameters, and $\theta$ is $\vv{x}$'s angle with respect to the horizontal.  This equation is often referred to as the \textit{polar rose} with $p = 2k + 1$ arms or petals, for $k \in \mathbb{N}^+$.  Likewise, we can control how fast curvature changes along $\Gamma$ by varying $a$ and $b$.  In particular, one can use \cref{eq:PolarRoseInterface} to compute

\begin{equation*}
\kappa_{rose}(\theta) = \frac{\gamma^2(\theta) + 2\left[\gamma'(\theta)\right]^2 - \gamma(\theta)\gamma''(\theta)}{\left( \gamma^2(\theta) + \left[\gamma'(\theta)\right]^2 \right)^{3/2}},
\end{equation*}
where $\gamma'$ and $\gamma''$ are the first and second derivatives of $\gamma(\theta)$.

In the following, we have chosen $p = 5$ and appropriate $a$ and $b$ values to validate \Cref{alg:MLCurvature}.  \Cref{fig:results.steep.charts} portrays several examples of the type of interfaces we are considering.


\colorsubsection{Steep-curvature flower-shaped interfaces}
\label{subsec:SteepCurvatureFlowerShapedInterfaces}

\begin{table}[!t]
	\centering
	\small
	\bgroup
	\def\arraystretch{1.1}%
	\begin{tabular}{cccccc}
		\hline
		\rowcolor{cloud1}
		$\eta$ & $m_\iota$ & $N_h^i$ & {\tt L2} factor & $|\mathcal{D}|$ & Parameters \\
		\hline \hline
		6      & 20 & 130 & $\eten{5}{-6}$ & 1,464,020 & 53,951 \\
		\hline
		7      & 18 & 130 & $\eten{5}{-6}$ & 1,473,359 & 53,691 \\
		\hline
		8      & 18 & 120 & $\eten{5}{-6}$ & 1,472,775 & 45,961 \\
		\hline
		9      & 18 & 130 & $\eten{5}{-6}$ & 1,464,407 & 53,691 \\
		\hline
		10     & 18 & 130 & $\eten{1}{-5}$ & 1,468,513 & 53,691 \\
		\hline
		11     & 18 & 120 & $\eten{7}{-6}$ & 1,462,741 & 45,961 \\
		\hline
	\end{tabular}
	\egroup
	\caption{Architectural data for $6 \leqslant \eta \leqslant 11$.}
	\label{tbl:results.train.config}
\end{table}

\begin{table}[!t]
	\centering
	\small
	\bgroup
	\def\arraystretch{1.1}%
	\begin{tabular}{clccc}
		\hline
		\rowcolor{cloud1}
		$\eta$              & Method   & MAE                   & MaxAE                 & RMSE \\
		\hline \hline
		\multirow{2}{*}{6}  & $\mathcal{F}_\kappa(\cdot)$ & $\eten{3.660305}{-4}$ & $\eten{1.397054}{-2}$ & $\eten{5.713122}{-4}$ \\
		                    & Baseline & $\eten{1.975861}{-2}$ & $\eten{1.542492}{-1}$ & $\eten{3.095603}{-2}$ \\
		\hline
		\multirow{2}{*}{7}  & $\mathcal{F}_\kappa(\cdot)$ & $\eten{3.923205}{-4}$ & $\eten{1.541521}{-2}$ & $\eten{6.269966}{-4}$ \\
		                    & Baseline & $\eten{1.977165}{-2}$ & $\eten{1.541942}{-1}$ & $\eten{3.090897}{-2}$ \\
		\hline
		\multirow{2}{*}{8}  & $\mathcal{F}_\kappa(\cdot)$ & $\eten{3.991765}{-4}$ & $\eten{1.369361}{-2}$ & $\eten{6.332734}{-4}$ \\
		                    & Baseline & $\eten{1.979684}{-2}$ & $\eten{1.549115}{-1}$ & $\eten{3.094043}{-2}$ \\
		\hline
		\multirow{2}{*}{9}  & $\mathcal{F}_\kappa(\cdot)$ & $\eten{3.952077}{-4}$ & $\eten{1.295839}{-2}$ & $\eten{6.336660}{-4}$ \\
		                    & Baseline & $\eten{1.983267}{-2}$ & $\eten{1.524857}{-1}$ & $\eten{3.097143}{-2}$ \\
		\hline
		\multirow{2}{*}{10} & $\mathcal{F}_\kappa(\cdot)$ & $\eten{4.427094}{-4}$ & $\eten{1.342434}{-2}$ & $\eten{7.062269}{-4}$ \\
		                    & Baseline & $\eten{1.982013}{-2}$ & $\eten{1.549863}{-1}$ & $\eten{3.095504}{-2}$ \\
		\hline
		\multirow{2}{*}{11} & $\mathcal{F}_\kappa(\cdot)$ & $\eten{4.293997}{-4}$ & $\eten{1.456045}{-2}$ & $\eten{6.879887}{-4}$ \\
		                    & Baseline & $\eten{1.983361}{-2}$ & $\eten{1.551221}{-1}$ & $\eten{3.098202}{-2}$ \\
		\hline
	\end{tabular}
	\egroup
	\caption{Training $h\kappa$ and $h\kappa_\mathcal{F}$ error statistics for $6 \leqslant \eta \leqslant 11$ over $\mathcal{D}$.}
	\label{tbl:results.train.stats}
\end{table}

Our first assessment tests how accurately \Cref{alg:MLCurvature} can estimate curvature on interfaces with poorly resolved regions.  At the same time, we show the flexibility of the schemes listed in \Cref{sec:Methodology} to train and deploy error-correcting models, irrespective of the mesh size.  For this purpose, we have optimized $\mathcal{F}_\kappa(\cdot)$ for $\nu = 10$ and $h=2^{-\eta}$, where $6 \leqslant \eta \leqslant 11$.  The interested reader may find these networks, the preprocessing objects, and usage documentation at \url{https://github.com/UCSB-CASL/Curvature_ECNet_2D}.  This repository includes {\tt frugally-deep} and {\tt HDF5} model versions alongside the PCA and standard scalers as {\tt pickle} and {\tt JSON} exports.  In addition, we have made $\mathcal{F}_\kappa(\cdot)$ available as a custom {\tt JSON} file with {\tt Base64}-encoded weights, excluding the output layer from \cref{fig:ECNet}.  The latter encoding helps us preserve TensorFlow's 32-bit precision and ensures compatibility with OpenBLAS (see \Cref{subsec:TechnicalAspects}).

\begin{figure}[!t]
	\centering
	\begin{subfigure}[!t]{0.32\textwidth}
		\includegraphics[width=\textwidth]{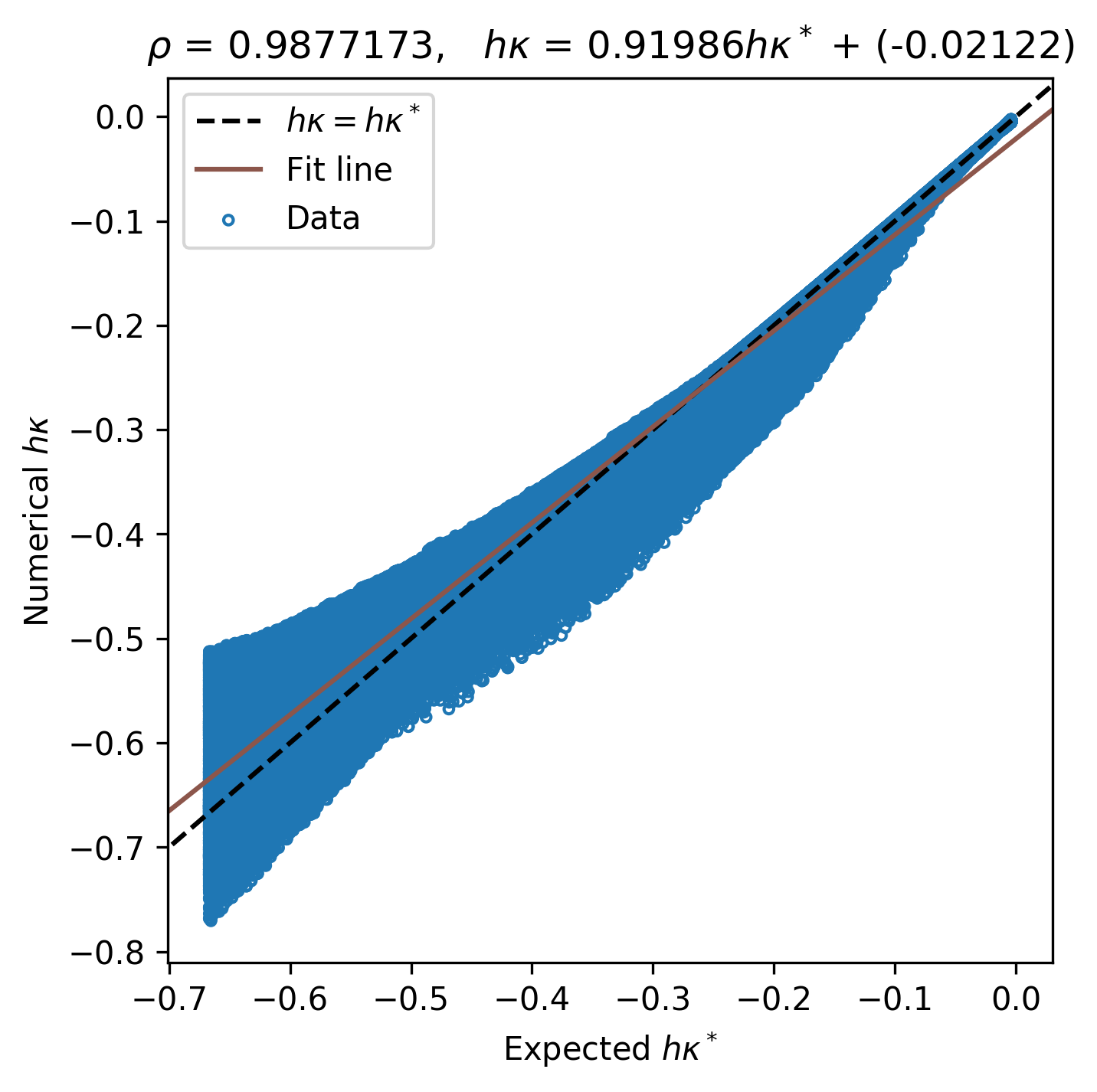}
	\end{subfigure}
	~
	\begin{subfigure}[!t]{0.32\textwidth}
		\includegraphics[width=\textwidth]{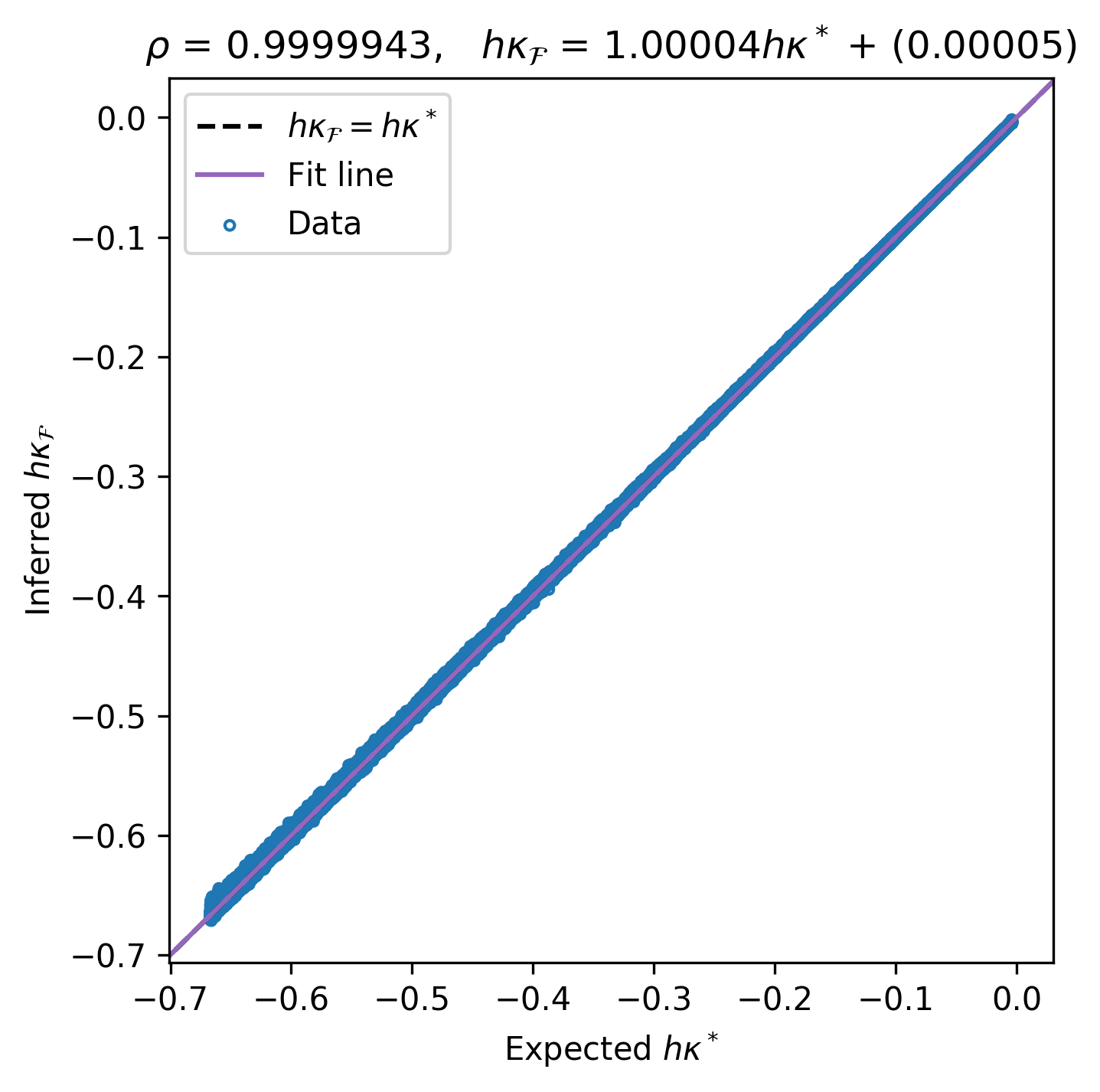}
	\end{subfigure}
	~
	\begin{subfigure}[!t]{0.32\textwidth}
		\includegraphics[width=\textwidth]{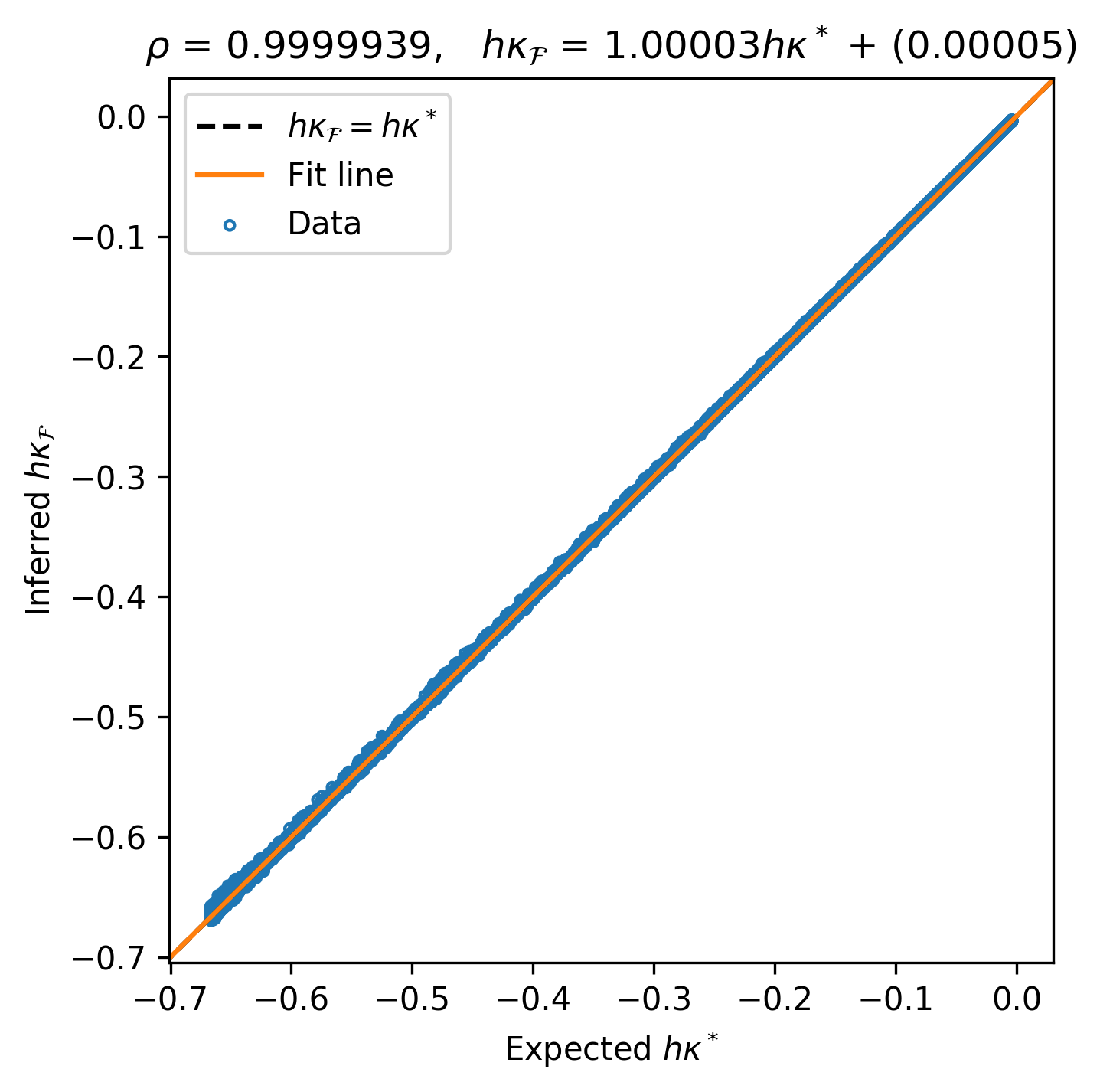}
	\end{subfigure}
	\\
	\begin{subfigure}[!t]{\textwidth}
		\caption{\footnotesize $h = 2^{-7}$}
		\label{fig:results.train.charts.7}
    \end{subfigure}
	\\
    
	\begin{subfigure}[!t]{0.32\textwidth}
		\includegraphics[width=\textwidth]{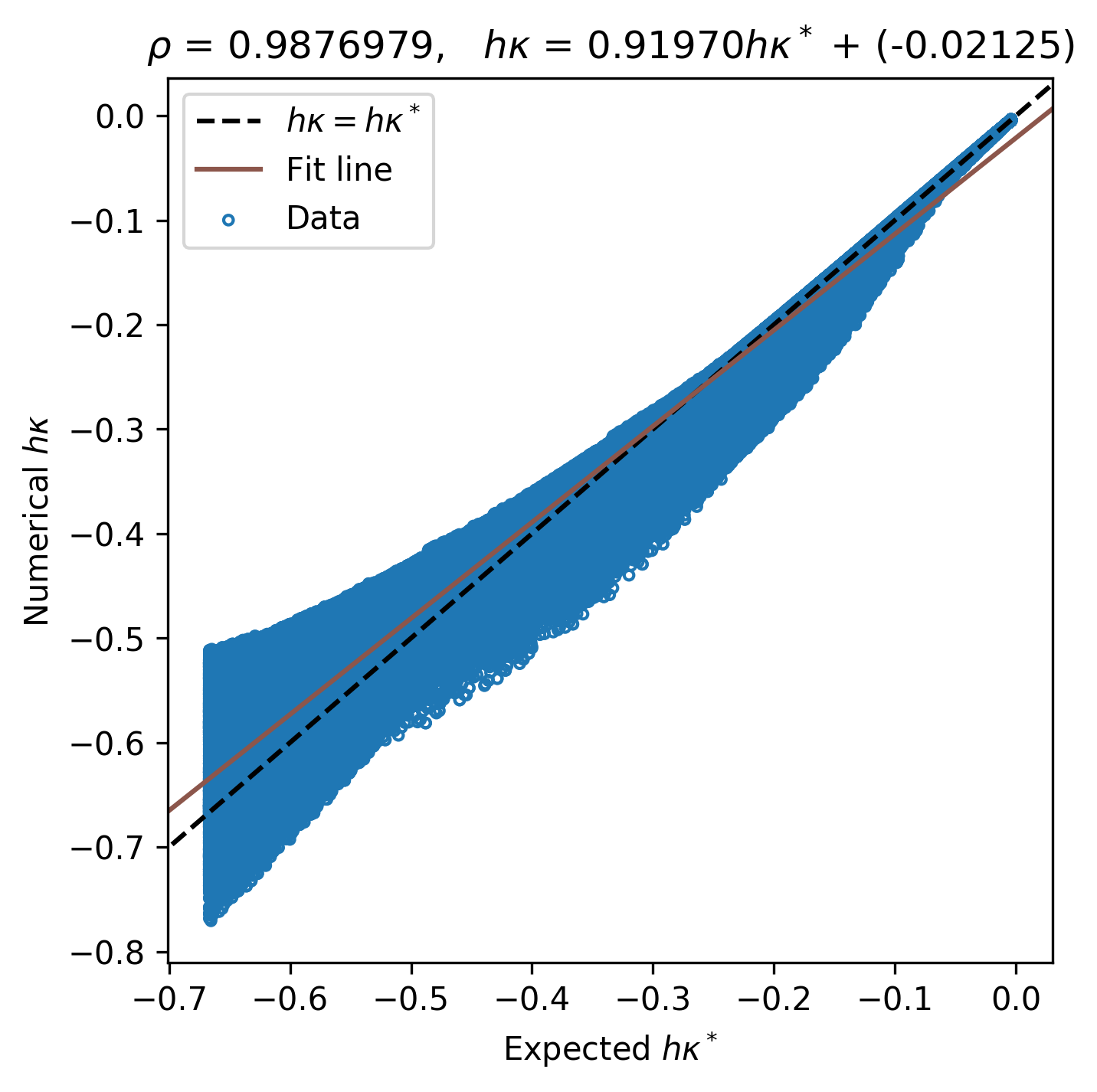}
	\end{subfigure}
	~
	\begin{subfigure}[!t]{0.32\textwidth}
		\includegraphics[width=\textwidth]{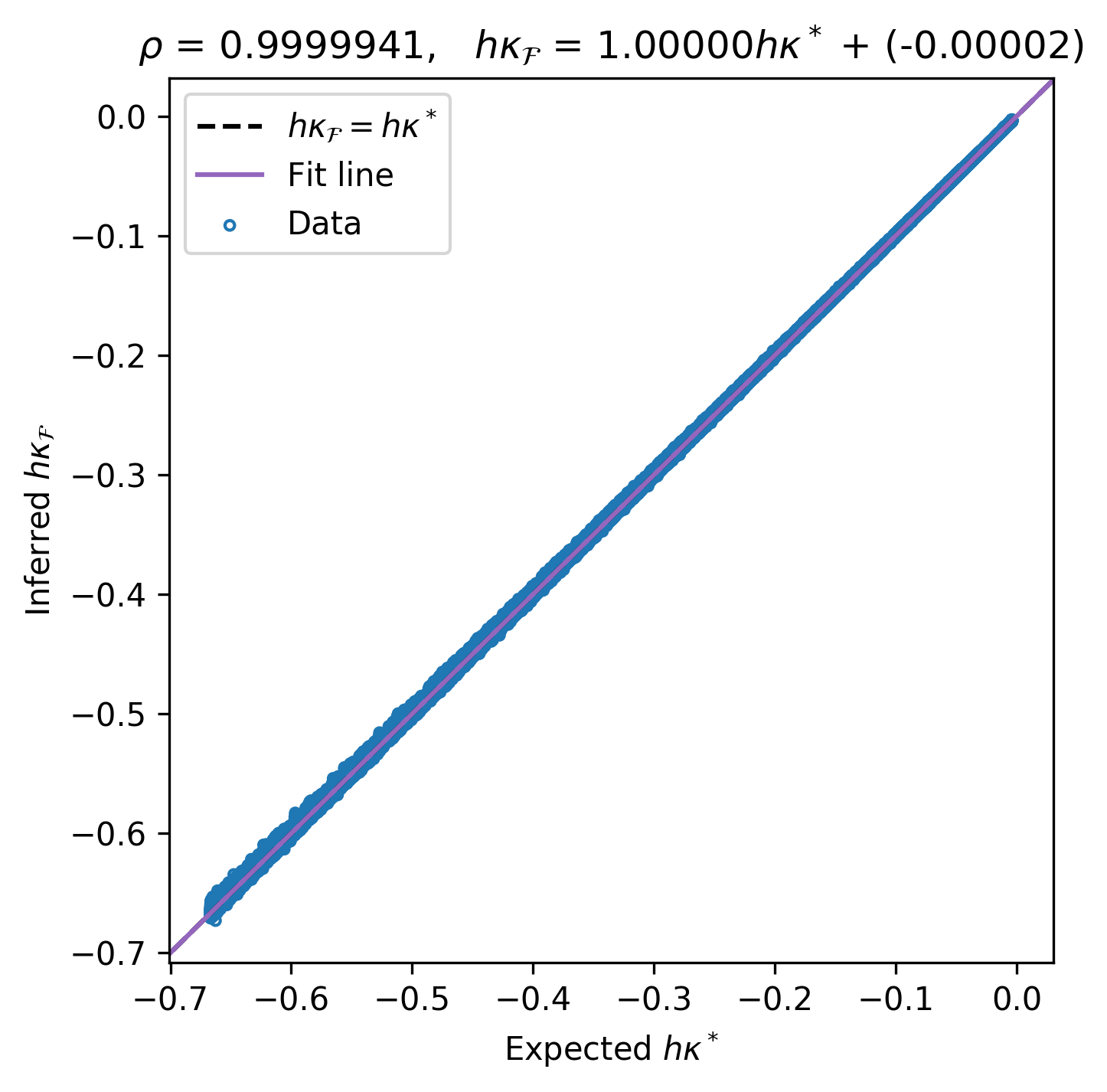}
	\end{subfigure}
	~
	\begin{subfigure}[!t]{0.32\textwidth}
		\includegraphics[width=\textwidth]{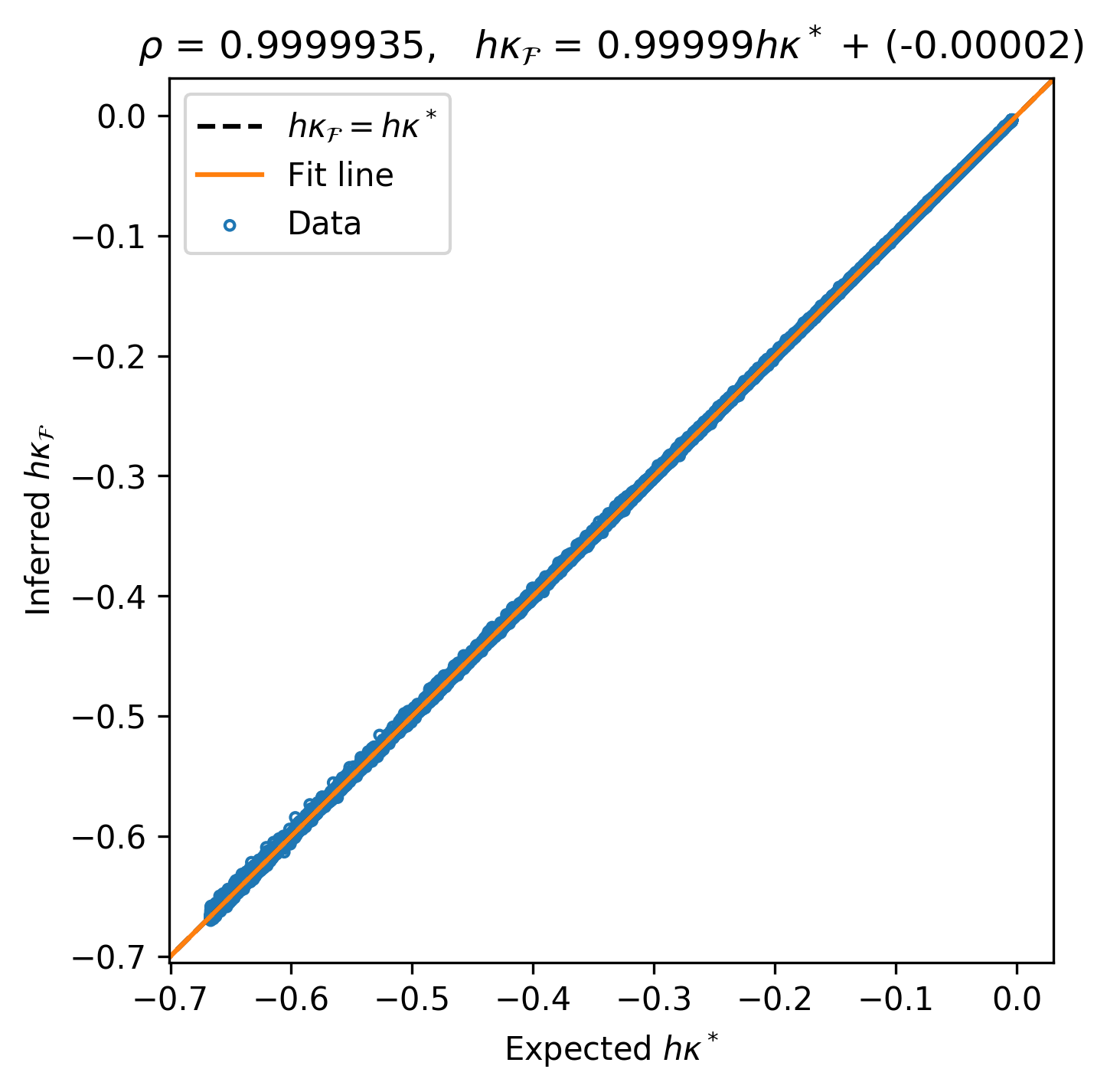}
    \end{subfigure}
	\\
	\begin{subfigure}[!t]{\textwidth}
		\caption{\footnotesize $h = 2^{-8}$}
		\label{fig:results.train.charts.8}
	\end{subfigure}
	\\
    
	\begin{subfigure}[!t]{0.32\textwidth}
		\includegraphics[width=\textwidth]{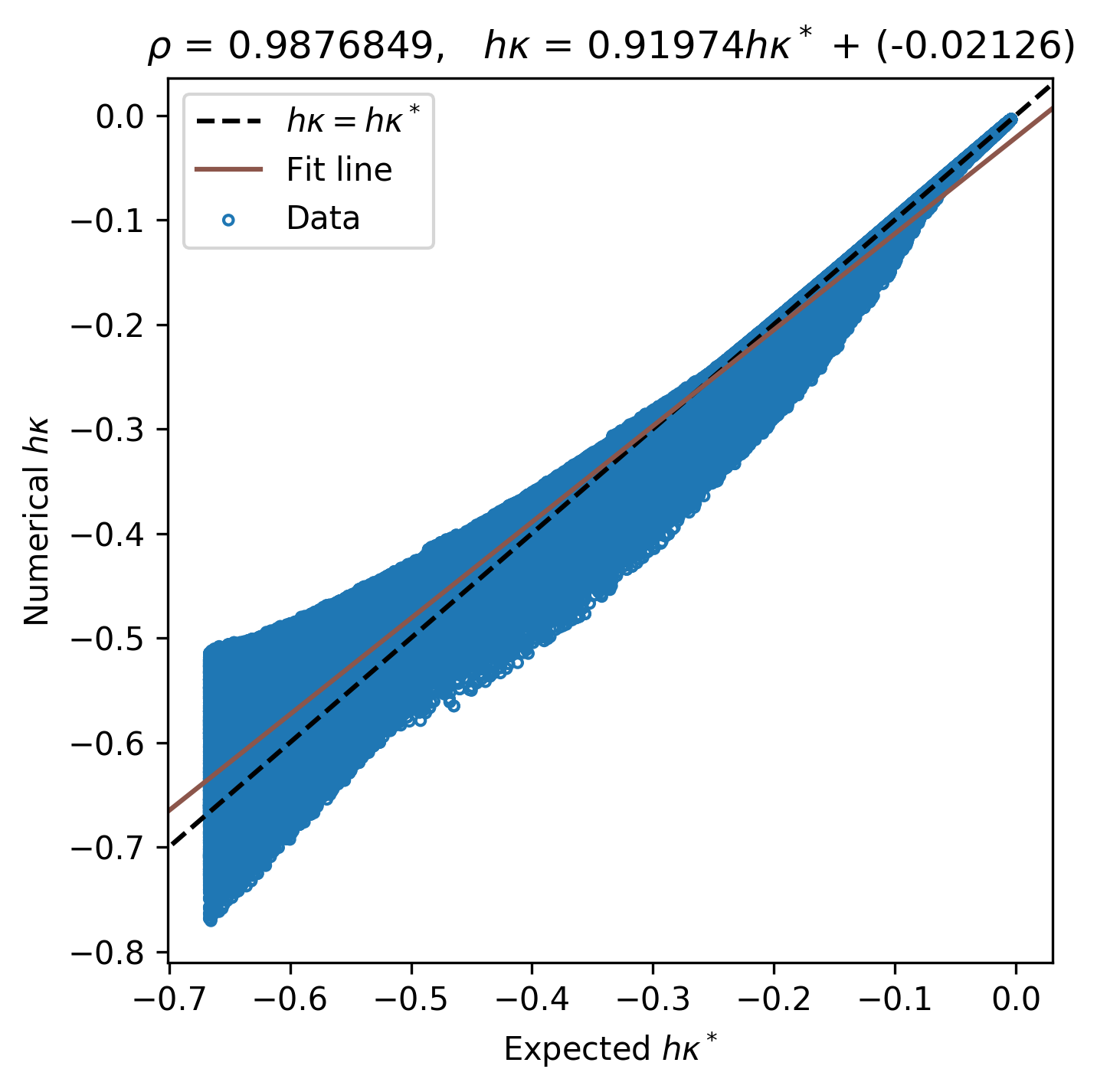}
    \end{subfigure}
	~
	\begin{subfigure}[!t]{0.32\textwidth}
		\includegraphics[width=\textwidth]{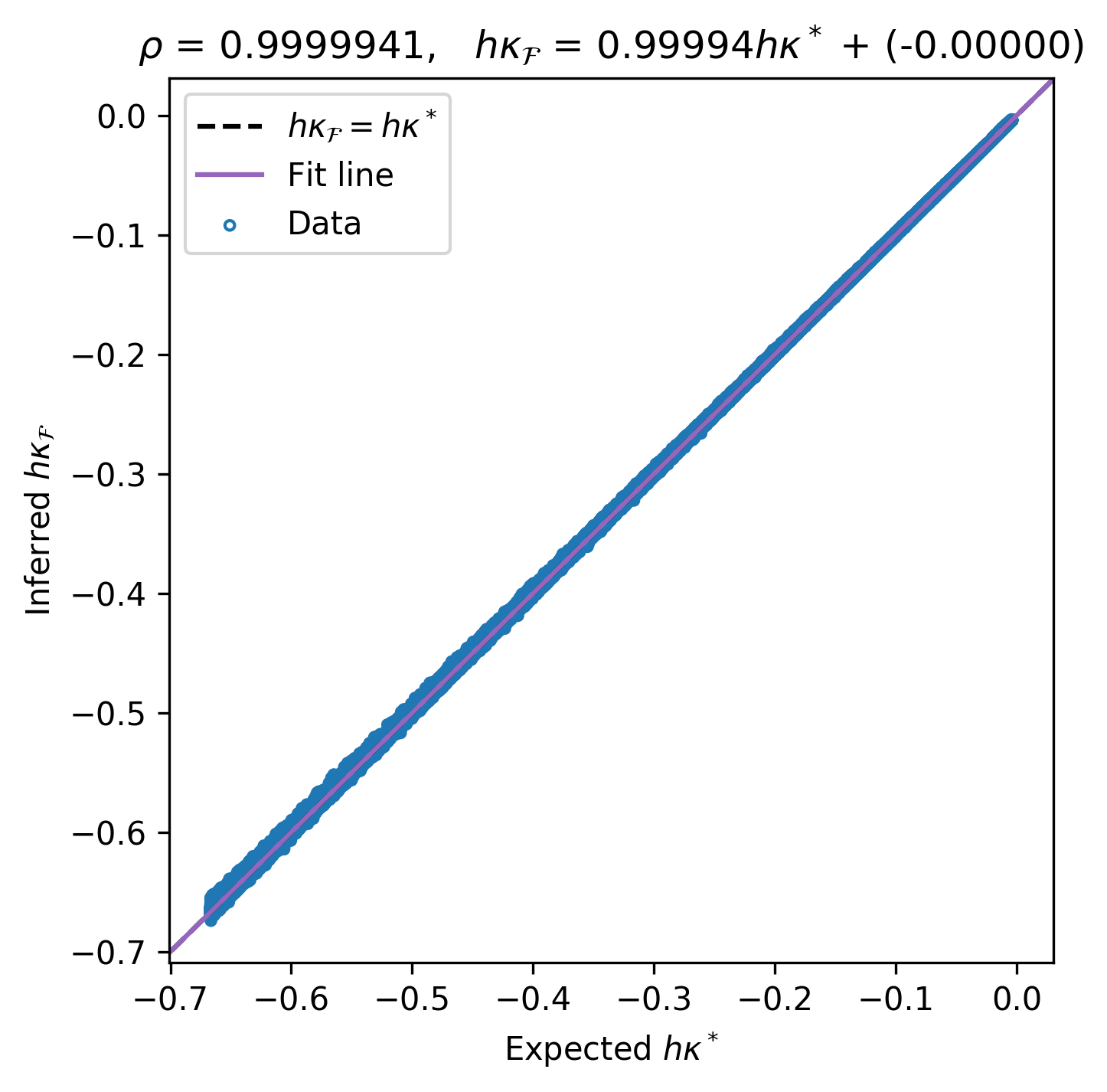}
    \end{subfigure}
	~
	\begin{subfigure}[!t]{0.32\textwidth}
		\includegraphics[width=\textwidth]{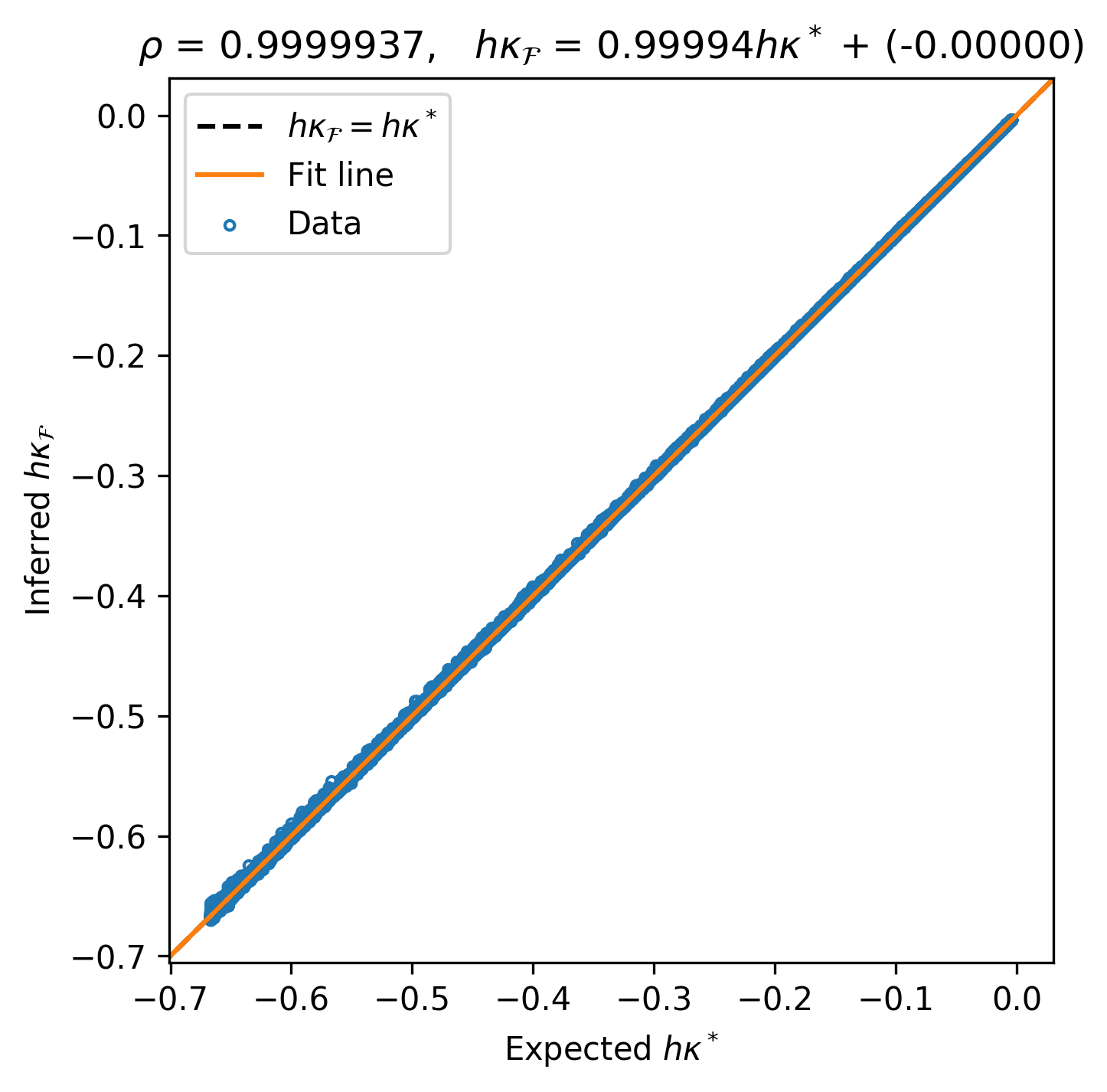}
    \end{subfigure}
    \\
	\begin{subfigure}[!t]{\textwidth}
		\caption{\footnotesize $h = 2^{-9}$}
		\label{fig:results.train.charts.9}
	\end{subfigure}
	\caption{Example correlation plots showing the training outcomes for $\eta = 7$ on row (a), $\eta = 8$ on row (b), and $\eta = 9$ on row (c).  The left column illustrates the quality of the numerical fit (i.e., the $h\kappa$ baseline), while the middle and right columns portray the $h\kappa_\mathcal{F}$ estimations in the entire $\mathcal{D}$ and the testing subset against the true $h\kappa^*$ values.  Each chart bears its $\rho$ correlation factor and the equation of the best regression line through the data.  (Color online.)}
	\label{fig:results.train.charts}
\end{figure}

\Cref{tbl:results.train.config} provides the architectural features for each $\eta$, such as the input-layer size, the number of neurons per hidden layer, the layer-wise {\tt L2} factor, and the total number of tunable parameters (see \cref{fig:ECNet}).  Similarly, \cref{tbl:results.train.stats} overviews the training outputs, including the tracked $h\kappa$- and $h\kappa_\mathcal{F}$-error metrics evaluated over the entire $\mathcal{D}$.  The most important fact is that ML optimization reduces the numerical error by one order of magnitude in the $L^\infty$ norm.  

\Cref{fig:results.train.charts} supplements the statistics in \cref{tbl:results.train.stats} with a few correlation charts.  For simplicity, it only provides the fit quality for $\eta = 7,8,9$ since the other grid resolutions lead essentially to the same results.  In particular, these charts contrast $h\kappa_\mathcal{F}$ and $h\kappa$ with the exact $h\kappa^*$ values.  As one reads through these plots, it is clear that our neural networks can improve curvature calculations by producing suitable corrections in response to preprocessed context information.

\begin{table}[!t]
	\centering
	\small
	\bgroup
	\def\arraystretch{1.1}%
	\begin{tabular}{ccclcrcrcr}
		\hline
		\rowcolor{cloud1}
		$\eta$ & $a$ & $b$ & Method & MAE & \makecell{Reduct.\\factor} & MaxAE & \makecell{Reduct.\\factor} & Time (sec.) & Cost (\%) \\
		\hline \hline
		\multirow{3}{*}{6}  & \multirow{3}{*}{0.085} & \multirow{3}{*}{0.300} 
		                    & {\tt MLCurvature()} & $\eten{2.27828}{-3}$ & -
		                                          & $\eten{1.67204}{-2}$ & -
		                                          & 0.010069              & - \\ 
		                & & & Baseline ($\nu=10$) & $\eten{9.31478}{-3}$ & \textcolor{goodgreen}{4.09} 
		                                          & $\eten{1.37868}{-1}$ & \textcolor{goodgreen}{8.25}
		                                          & 0.007960              & +26.5 \\
		                & & & Baseline ($\nu=20$) & $\eten{9.17273}{-3}$ & \textcolor{goodgreen}{4.03} 
		                                          & $\eten{1.37366}{-1}$ & \textcolor{goodgreen}{8.22}
		                                          & 0.014384              & -30.0 \\
		\hline
		\multirow{3}{*}{7}  & \multirow{3}{*}{0.120} & \multirow{3}{*}{0.305} 
		                    & {\tt MLCurvature()} & $\eten{7.37148}{-4}$ & -
		                                          & $\eten{1.36763}{-2}$ & -
		                                          & 0.020830              & - \\
		                & & & Baseline ($\nu=10$) & $\eten{3.46929}{-3}$ & \textcolor{goodgreen}{4.71}
		                                          & $\eten{1.34981}{-1}$ & \textcolor{goodgreen}{9.87}
		                                          & 0.016341              & +27.5 \\
		                & & & Baseline ($\nu=20$) & $\eten{3.16221}{-3}$ & \textcolor{goodgreen}{4.29}
		                                          & $\eten{1.36403}{-1}$ & \textcolor{goodgreen}{9.97}
		                                          & 0.031805              & -34.5 \\
		\hline
		\multirow{3}{*}{8}  & \multirow{3}{*}{0.170} & \multirow{3}{*}{0.330} 
		                    & {\tt MLCurvature()} & $\eten{4.59714}{-4}$ & -
		                                          & $\eten{2.27442}{-2}$ & -
		                                          & 0.050951              & - \\
		                & & & Baseline ($\nu=10$) & $\eten{1.43790}{-3}$ & \textcolor{goodgreen}{3.13} 
		                                          & $\eten{9.83639}{-2}$ & \textcolor{goodgreen}{4.32}
		                                          & 0.041263              & +23.5 \\
		                & & & Baseline ($\nu=20$) & $\eten{1.12260}{-3}$ & \textcolor{goodgreen}{2.44}
		                                          & $\eten{9.93208}{-2}$ & \textcolor{goodgreen}{4.37}
		                                          & 0.076359              & -33.3 \\
		\hline
		\multirow{3}{*}{9}  & \multirow{3}{*}{0.225} & \multirow{3}{*}{0.355} 
		                    & {\tt MLCurvature()} & $\eten{3.91830}{-4}$ & -
		                                          & $\eten{9.21459}{-3}$ & -
		                                          & 0.110962              & - \\
		                & & & Baseline ($\nu=10$) & $\eten{8.23963}{-4}$ & \textcolor{goodgreen}{2.10}
		                                          & $\eten{1.05026}{-1}$ & \textcolor{goodgreen}{11.40}
		                                          & 0.093926              & +18.1 \\
		                & & & Baseline ($\nu=20$) & $\eten{4.97385}{-4}$ & \textcolor{goodgreen}{1.27}
		                                          & $\eten{1.05255}{-1}$ & \textcolor{goodgreen}{11.42}
		                                          & 0.180637              & -38.6 \\
		\hline
		\multirow{3}{*}{10} & \multirow{3}{*}{0.258} & \multirow{3}{*}{0.356} 
		                    & {\tt MLCurvature()} & $\eten{4.01385}{-4}$ & -
		                                          & $\eten{1.20870}{-2}$ & -
		                                          & 0.229452              & - \\
		                & & & Baseline ($\nu=10$) & $\eten{5.84885}{-4}$ & \textcolor{goodgreen}{1.46}
		                                          & $\eten{9.31230}{-2}$ & \textcolor{goodgreen}{7.70}
		                                          & 0.207452              & +10.6 \\
		                & & & Baseline ($\nu=20$) & $\eten{2.05082}{-4}$ & \textcolor{darkcoral}{0.51}
		                                          & $\eten{9.39947}{-2}$ & \textcolor{goodgreen}{7.78}
		                                          & 0.411581              & -44.3 \\
		\hline
		\multirow{3}{*}{11} & \multirow{3}{*}{0.274} & \multirow{3}{*}{0.345} 
		                    & {\tt MLCurvature()} & $\eten{4.60019}{-4}$ & -
		                                          & $\eten{9.06714}{-3}$ & -
		                                          & 0.510308              & - \\
		                & & & Baseline ($\nu=10$) & $\eten{5.44016}{-4}$ & \textcolor{goodgreen}{1.18}
		                                          & $\eten{8.50492}{-2}$ & \textcolor{goodgreen}{9.38}
		                                          & 0.473656              & +7.7 \\
		                & & & Baseline ($\nu=20$) & $\eten{1.03791}{-4}$ & \textcolor{darkcoral}{0.23} 
		                                          & $\eten{8.43886}{-2}$ & \textcolor{goodgreen}{9.31}
		                                          & 0.917335              & -44.4 \\
		\hline
	\end{tabular}
	\egroup
	\caption{Shape parameters for steep-curvature flower-shaped interfaces and the $h\kappa$ and $h\kappa^\star$ error and performance statistics for $6 \leqslant \eta \leqslant 11$ and $\nu = 10, 20$.  Error-reduction factors in favor of {\tt MLCurvature()} appear in \textcolor{goodgreen}{green}, while those in favor of the conventional numerical approach are given in \textcolor{darkcoral}{red}.  (Color online.)}
	\label{tbl:results.steep.stats}
\end{table}

Now, we evaluate our error-correcting models on $\phi_{rose}(\vv{x})$ by integrating them into the {\tt MLCurvature()} routine.  Our tests require computing curvature for steep flower-shaped interfaces, where $3/5 < \max |h\kappa^*| < 2/3$.  The shape parameters for \cref{eq:PolarRoseLevelSetFunction,eq:PolarRoseInterface} and the $h\kappa$- and $h\kappa^\star$-error metrics for each $\eta$ appear in \cref{tbl:results.steep.stats}.  The left column in \cref{fig:results.steep.charts} displays the analyzed free boundary for every $h$.

\Cref{tbl:results.steep.stats} includes auxiliary MAE- and MaxAE-reduction factors, too, to underline the advantages of our hybrid strategy.  For reference, we also provide running wall-times and baseline results collected for twenty redistancing steps.  As discussed below, these data help us gauge the tradeoff between {\tt MLCurvature()} and the standard numerical approach with twice the reinitialization cost.  

First, \cref{tbl:results.steep.stats} shows that \Cref{alg:MLCurvature} leads to better approximations in the $L^1$ norm for relatively coarse grids.  And the same is true in most scenarios when $\nu = 20$.  More specifically, as $h \to 0$, both $h\kappa$ and $h\kappa^\star$ converge to the same quantity for $\nu = 10$.  When $\nu = 20$, however, the situation completely turns in favor of the baseline.  As expected, this behavior results from the increasing number of well-resolved samples, which bias the average error in the benefit of the traditional method.  On the contrary, the MaxAE columns emphasize the strengths of our error-correcting approach.  In this case, the {\tt MLCurvature()}'s MaxAE remains consistently below 24\% of the error incurred by conventional means.  At times, these accuracy factors are as high as one order of magnitude, corroborating the $L^\infty$ error ratio reported in \cref{tbl:results.train.stats}.

The last pair of columns in \cref{tbl:results.steep.stats} summarizes the performance statistics for all the case studies.  We collected these figures from a single-process OpenBLAS-batch-based implementation of \Cref{alg:MLCurvature} in a 2.2GHz Intel Core i7 laptop with 16GB RAM.  \Cref{tbl:results.steep.stats} reports the shortest time out of ten trials for each configuration.  And the cost column, in particular, shows how much more (+) or less (-) time the hybrid procedure needed to complete the task (relative to the baseline).  Clearly, there is always a non-zero cost associated with $\mathcal{F}_\kappa(\cdot)$; however, the potential MaxAE reduction and the low coarse-grid MAE largely out-weight the ML cost in the benefit of an FBP solver.  Notice, too, that the neural system built on $\nu = 10$ is consistently superior in the $L^\infty$ error norm to the baseline with twice the number of redistancing steps.  The empirical savings stand above 30\%, especially for large mesh sizes and in the presence of under-resolution.

\Cref{fig:results.steep.charts} concludes this section with a handful of plots per mesh size for $\nu = 10$.  The first-column diagrams depict the five-petaled interfaces with the nodes where we estimated curvature.  Also, note we have deepened the concavities along the edge to stress $\mathcal{F}_\kappa(\cdot)$'s capabilities in every $h$.  Similar to \cref{fig:results.train.charts}, the second and third columns illustrate the fit quality for the numerical baseline and the ML-corrected approximations.  These charts demonstrate that our {\tt MLCurvature()} routine vastly improves $h\kappa$ as $|h\kappa^*|$ increases.  In like manner, the higher correlation and better regression lines across the third column offer unmistakable evidence that our strategy works well around under-resolved regions.

\begin{figure}[!t]
	\centering
	\begin{subfigure}[!t]{0.30\textwidth}
		\includegraphics[width=\textwidth]{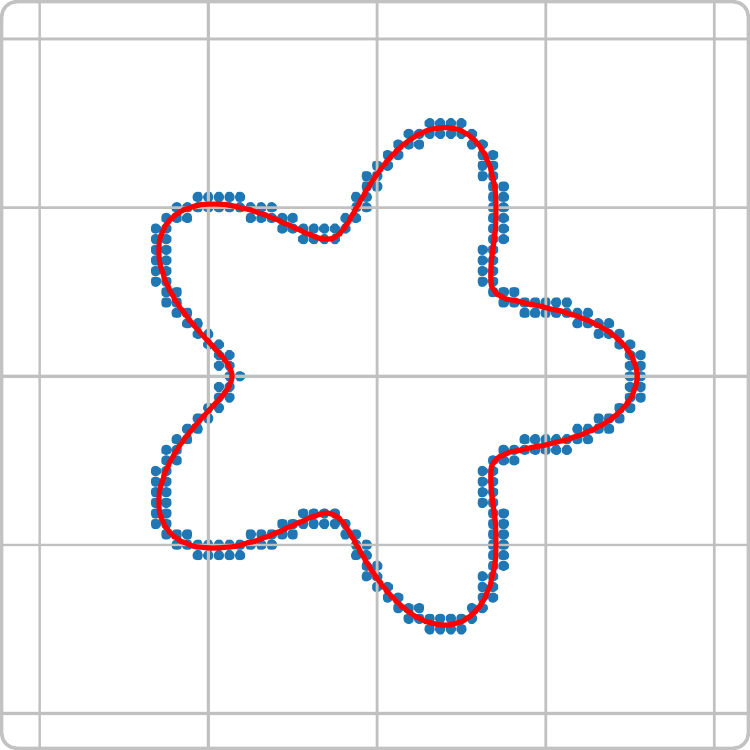}
	\end{subfigure}
	~
	\begin{subfigure}[!t]{0.33\textwidth}
		\includegraphics[width=\textwidth]{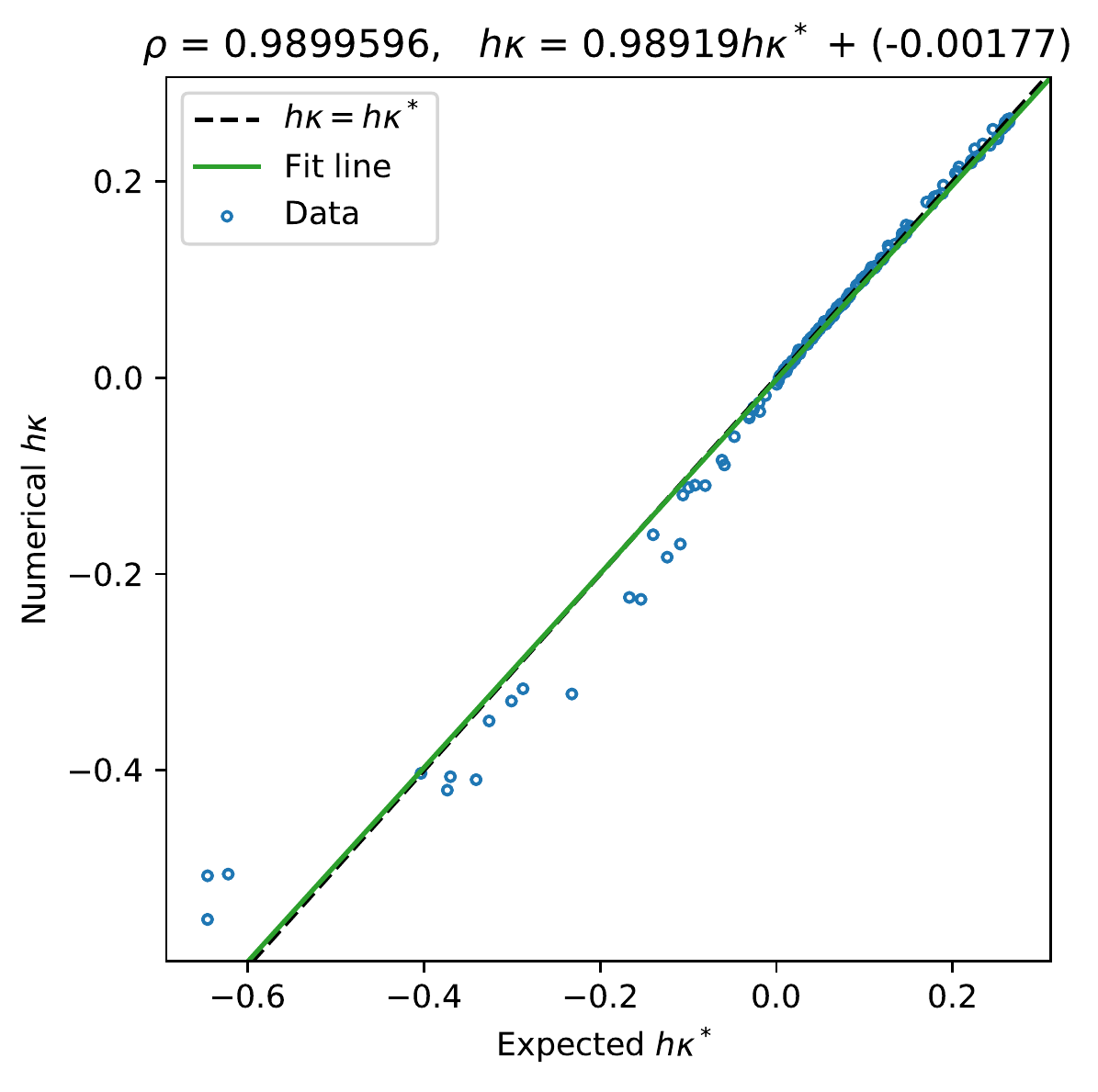}
	\end{subfigure}
	~
	\begin{subfigure}[!t]{0.33\textwidth}
		\includegraphics[width=\textwidth]{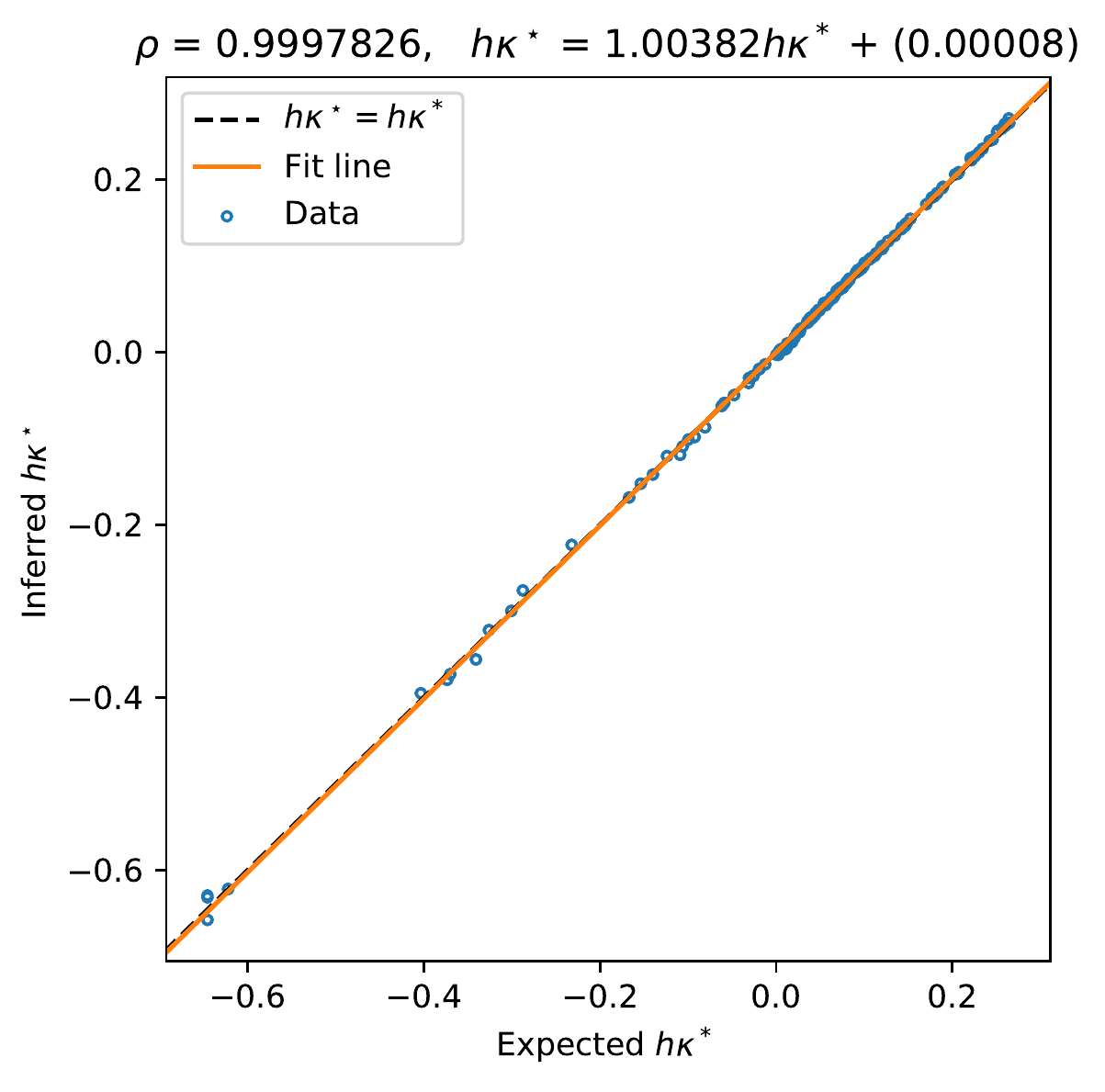}
	\end{subfigure}
	\\
	\begin{subfigure}[!t]{\textwidth}
		\caption{\footnotesize $h = 2^{-6}$}
		\label{fig:results.steep.charts.6}
	\end{subfigure}
	\\
	
	\begin{subfigure}[!t]{0.30\textwidth}
		\includegraphics[width=\textwidth]{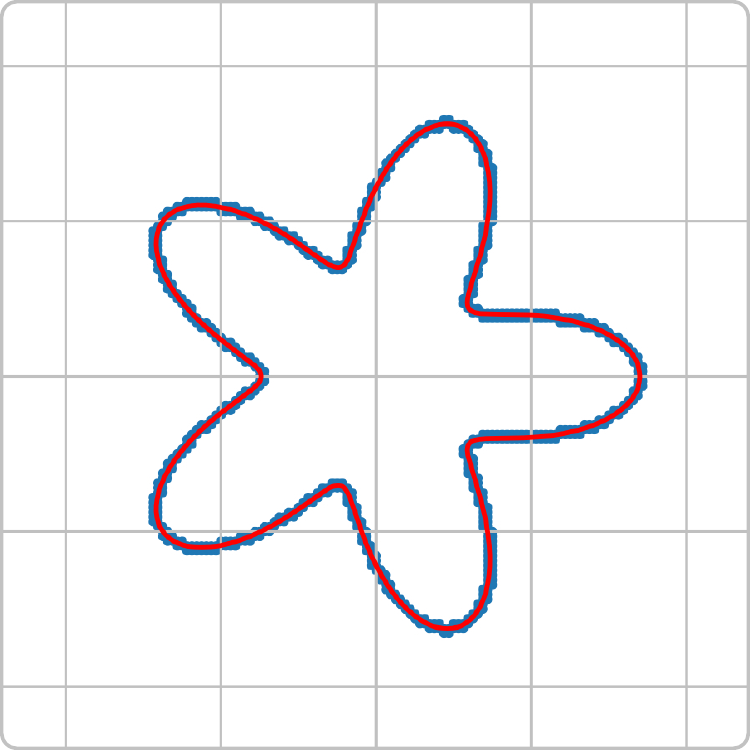}
	\end{subfigure}
	~
	\begin{subfigure}[!t]{0.33\textwidth}
		\includegraphics[width=\textwidth]{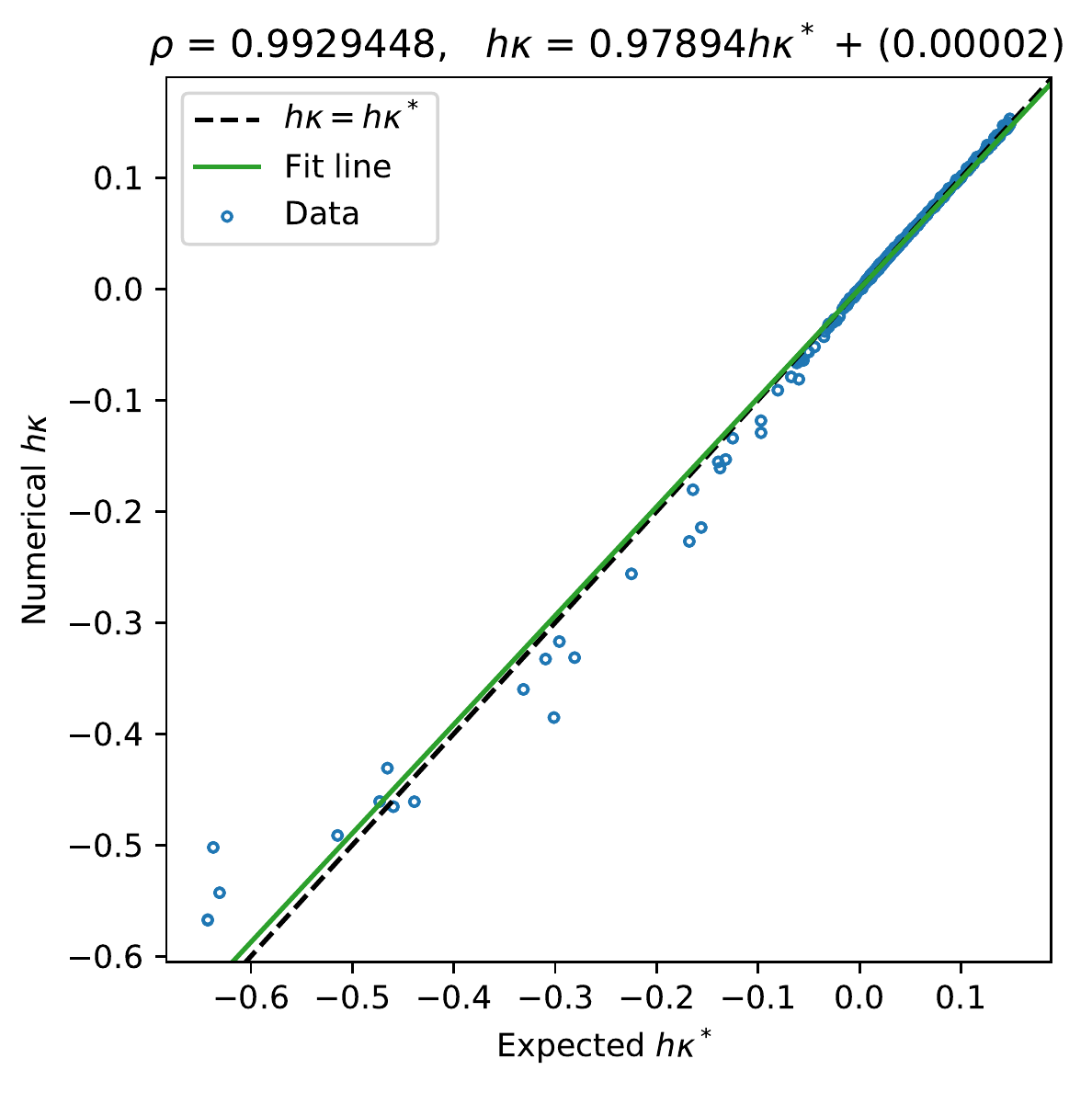}
	\end{subfigure}
	~
	\begin{subfigure}[!t]{0.33\textwidth}
		\includegraphics[width=\textwidth]{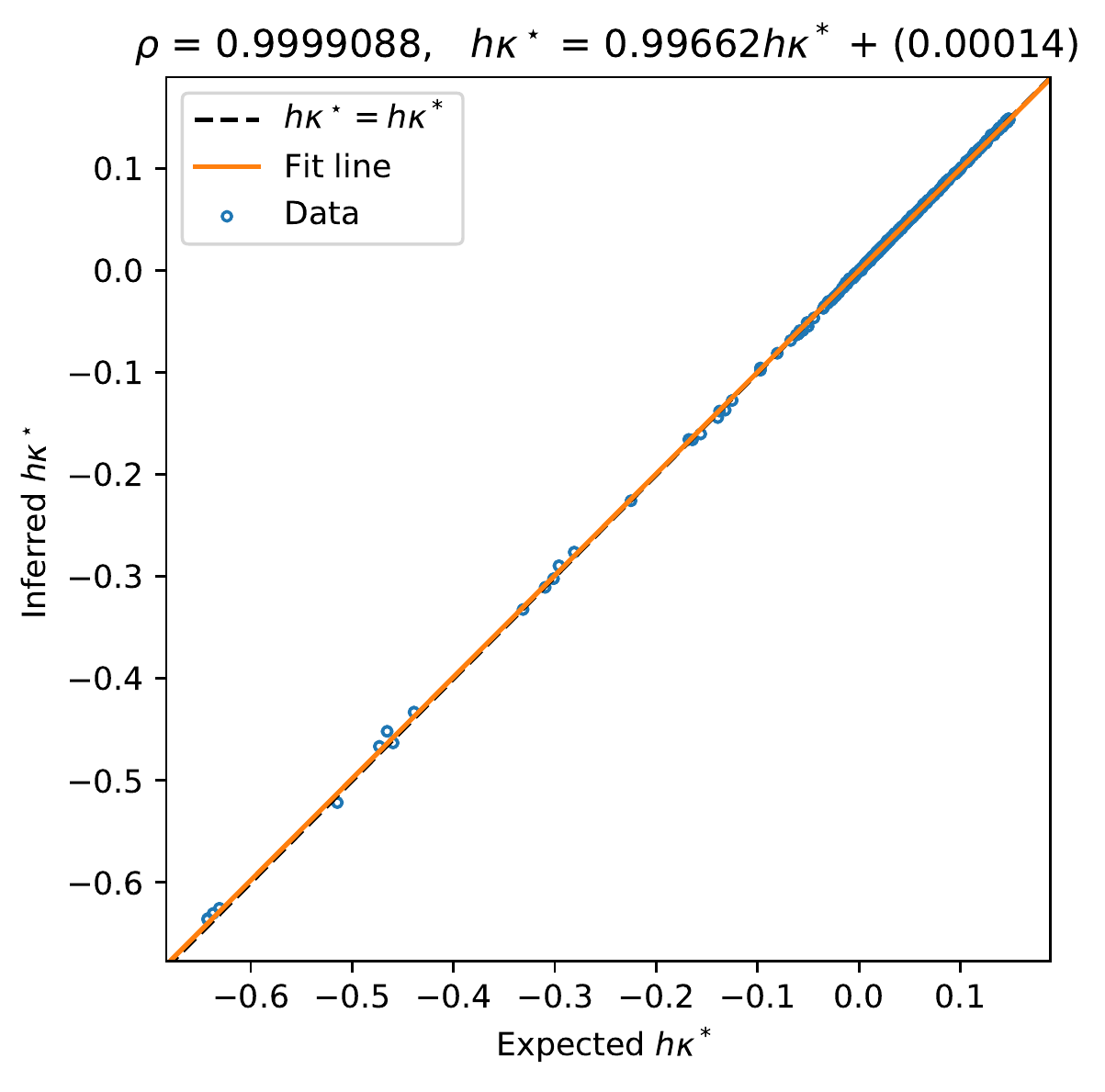}
	\end{subfigure}
	\\
	\begin{subfigure}[!t]{\textwidth}
		\caption{\footnotesize $h = 2^{-7}$}
		\label{fig:results.steep.charts.7}
	\end{subfigure}
	\\
    
	\begin{subfigure}[!t]{0.30\textwidth}
		\includegraphics[width=\textwidth]{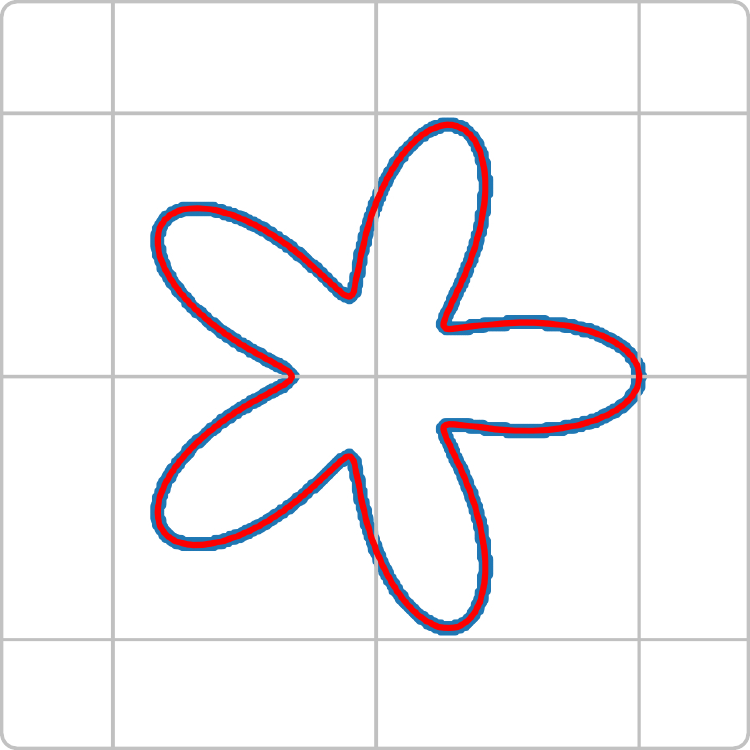}
	\end{subfigure}
	~
	\begin{subfigure}[!t]{0.33\textwidth}
		\includegraphics[width=\textwidth]{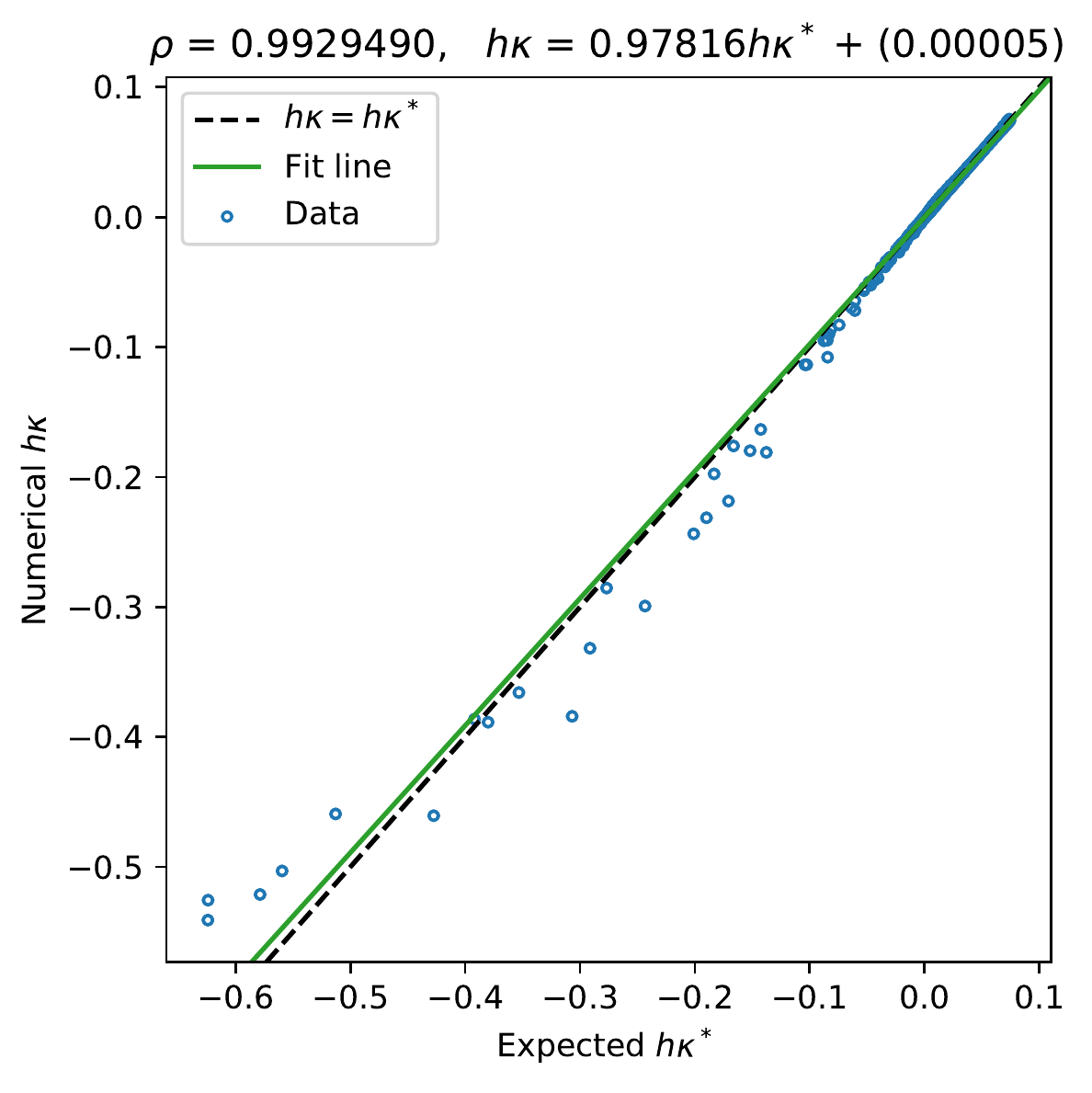}
	\end{subfigure}
	~
	\begin{subfigure}[!t]{0.33\textwidth}
		\includegraphics[width=\textwidth]{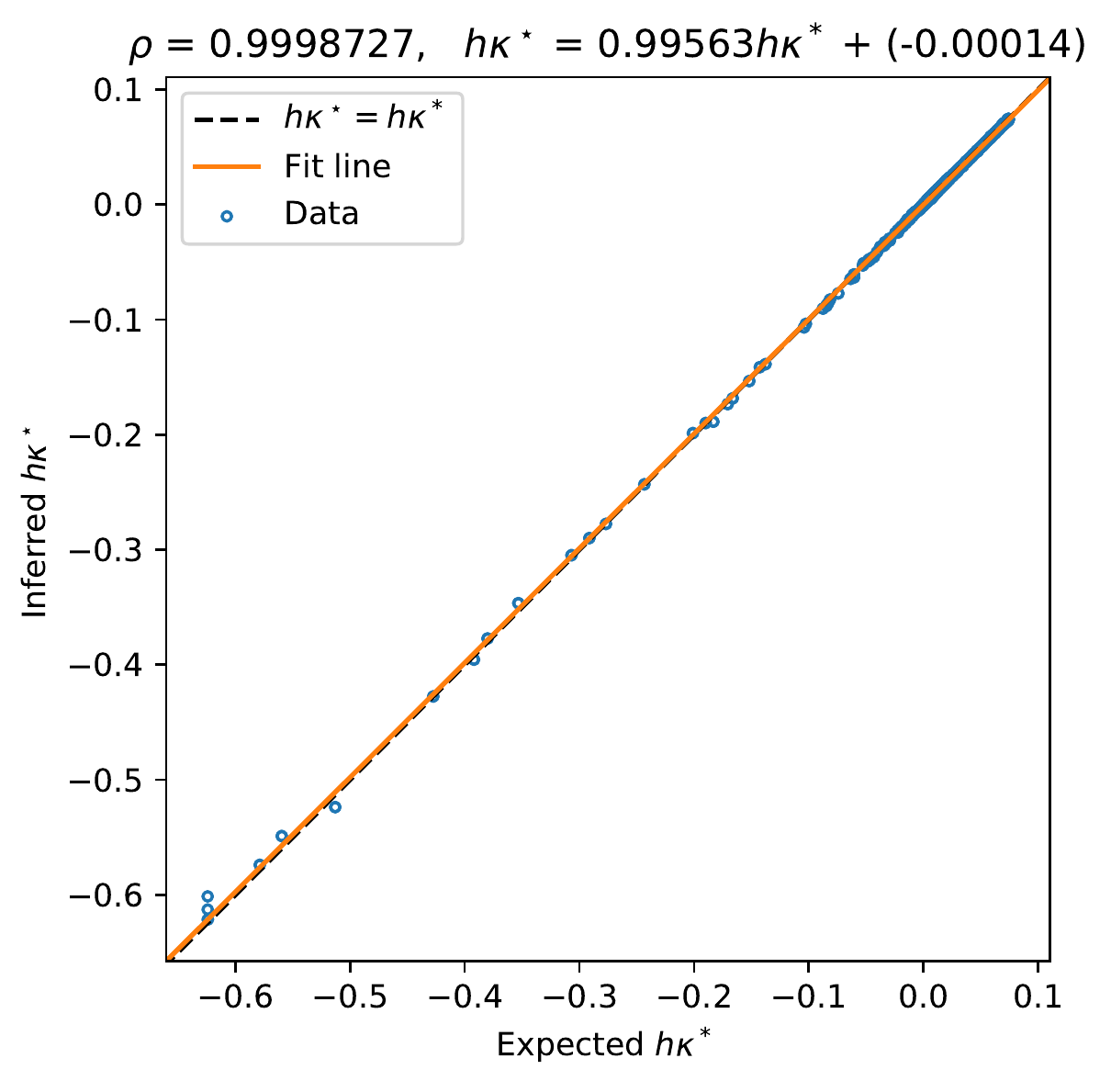}
	\end{subfigure}
	\\
	\begin{subfigure}[!t]{\textwidth}
		\caption{\footnotesize $h = 2^{-8}$}
		\label{fig:results.steep.charts.8}
	\end{subfigure}
	   
	\caption{Steep-curvature flower-shaped interfaces and their correlation plots.  The first column shows the tested interface with concavities varying according to the grid resolution.  We have contoured $\gamma(\theta)$ in red and denoted the sampled grid points next to $\Gamma$ in blue.  The second and third columns depict the fit quality for the numerical baseline (in green) and {\tt MLCurvature()} (in orange) for $\nu = 10$.  Row (a) belongs to $\eta = 6$, row (b) corresponds to $\eta = 7$, and row (c) is associated to $\eta = 8$.  (Color online.)}
	\label{fig:results.steep.charts}
\end{figure}

\begin{figure}[!t]
	\ContinuedFloat
	\centering
	\begin{subfigure}[!t]{0.30\textwidth}
		\includegraphics[width=\textwidth]{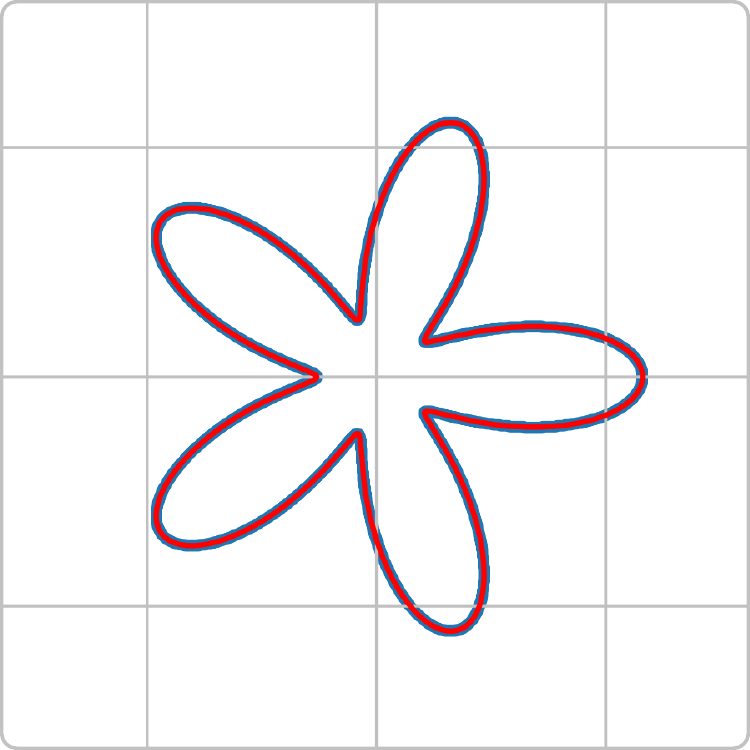}
	\end{subfigure}
	~
	\begin{subfigure}[!t]{0.33\textwidth}
		\includegraphics[width=\textwidth]{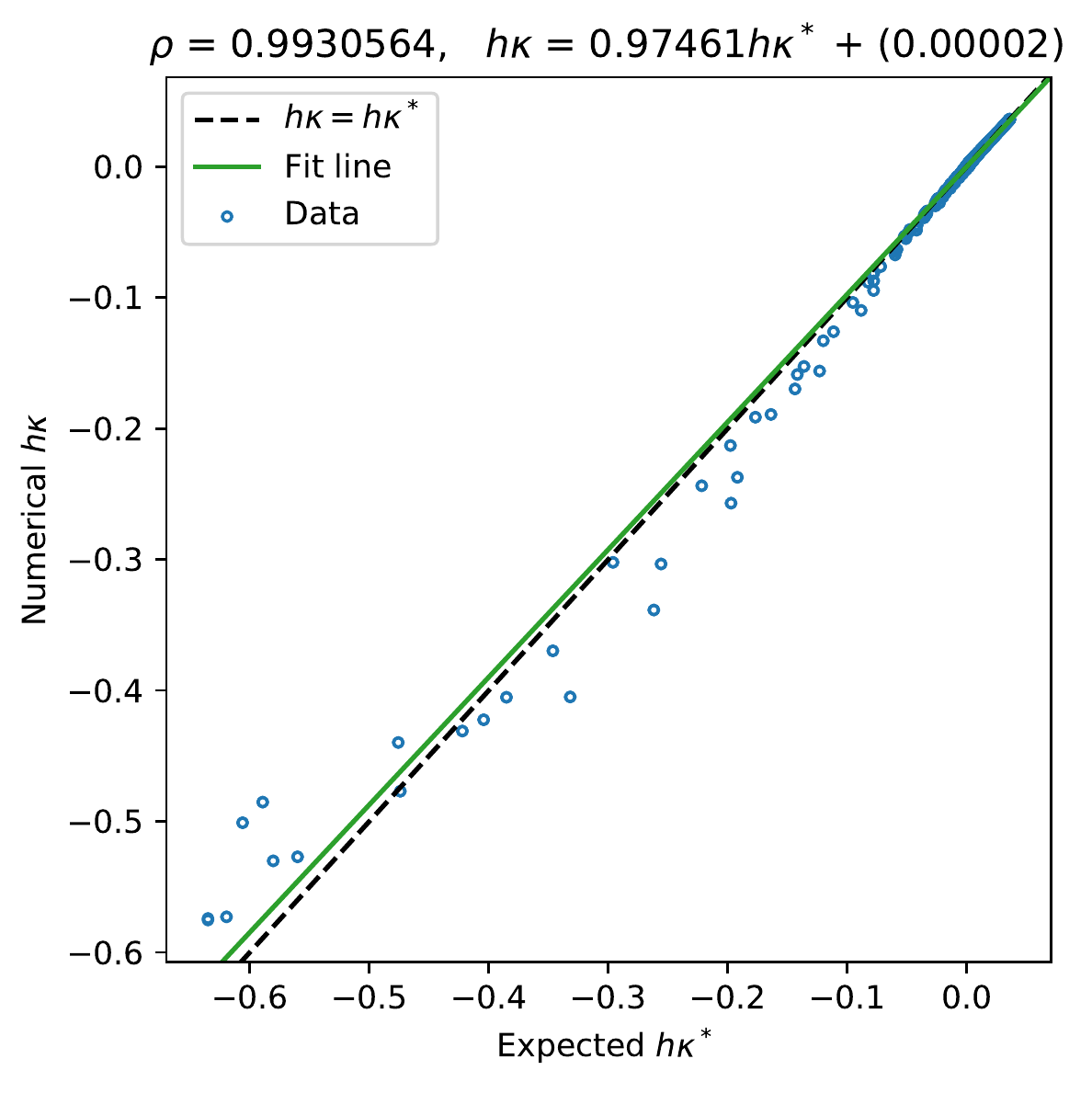}
	\end{subfigure}
	~
	\begin{subfigure}[!t]{0.33\textwidth}
		\includegraphics[width=\textwidth]{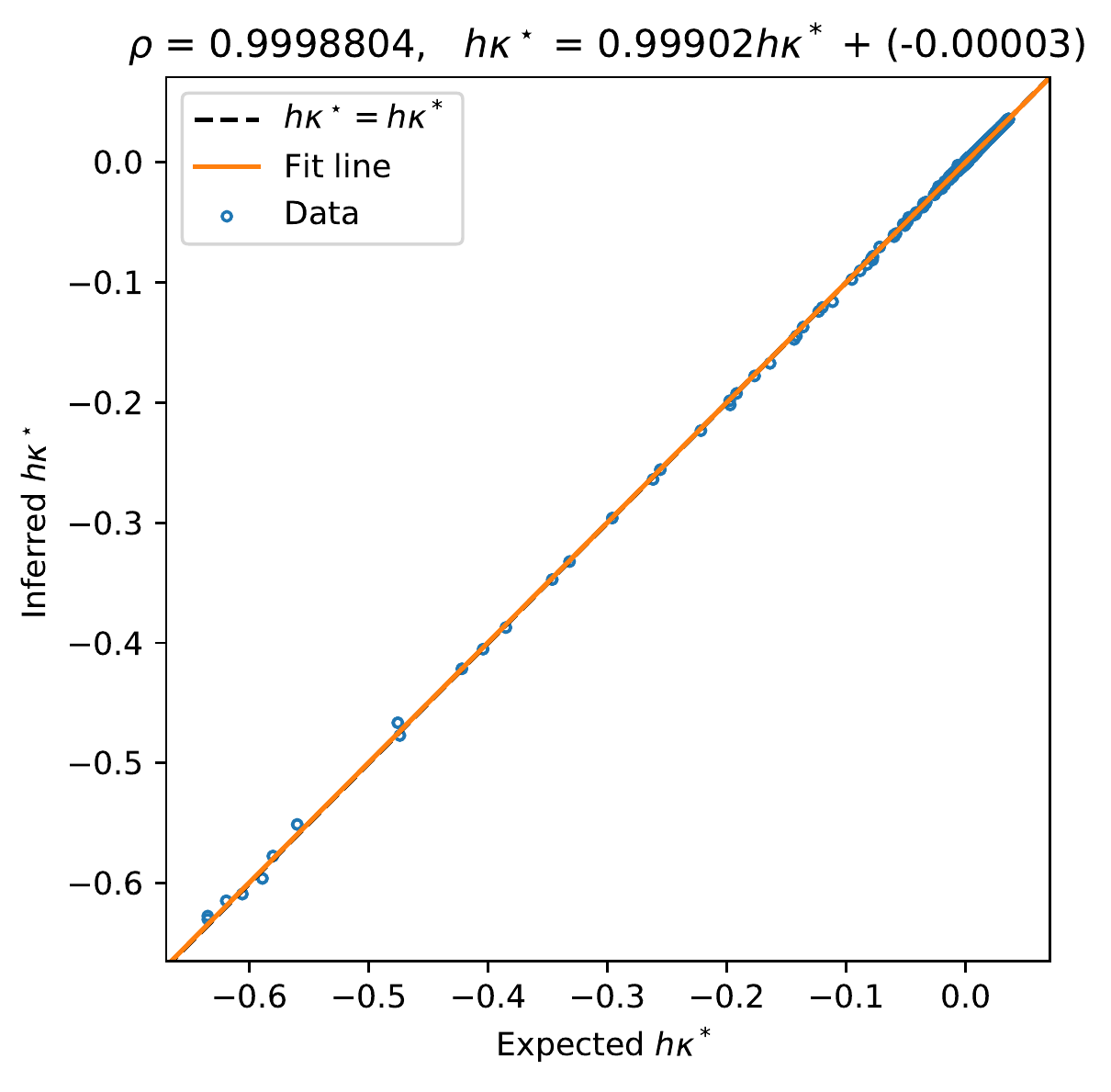}
    \end{subfigure}
	\\
	\begin{subfigure}[!t]{\textwidth}
		\caption{\footnotesize $h = 2^{-9}$}
		\label{fig:results.steep.charts.9}
	\end{subfigure}
	\\
    
	\begin{subfigure}[!t]{0.30\textwidth}
		\includegraphics[width=\textwidth]{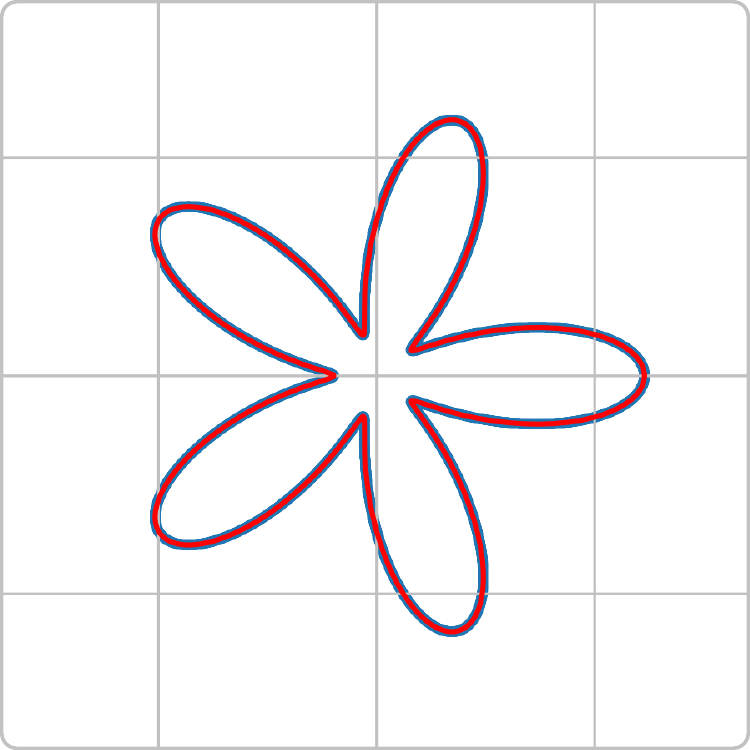}
	\end{subfigure}
	~
	\begin{subfigure}[!t]{0.33\textwidth}
		\includegraphics[width=\textwidth]{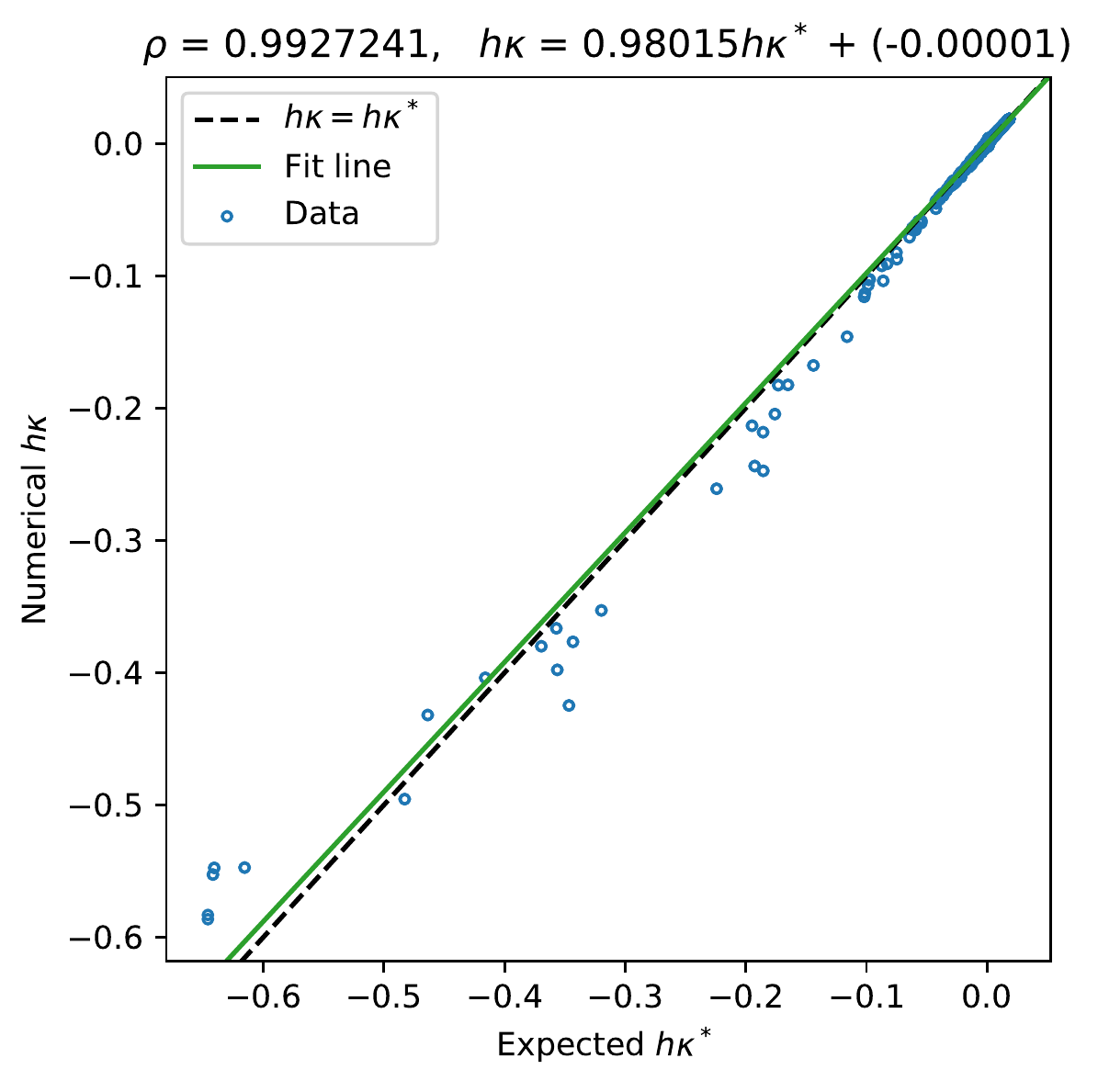}
	\end{subfigure}
	~
	\begin{subfigure}[!t]{0.33\textwidth}
		\includegraphics[width=\textwidth]{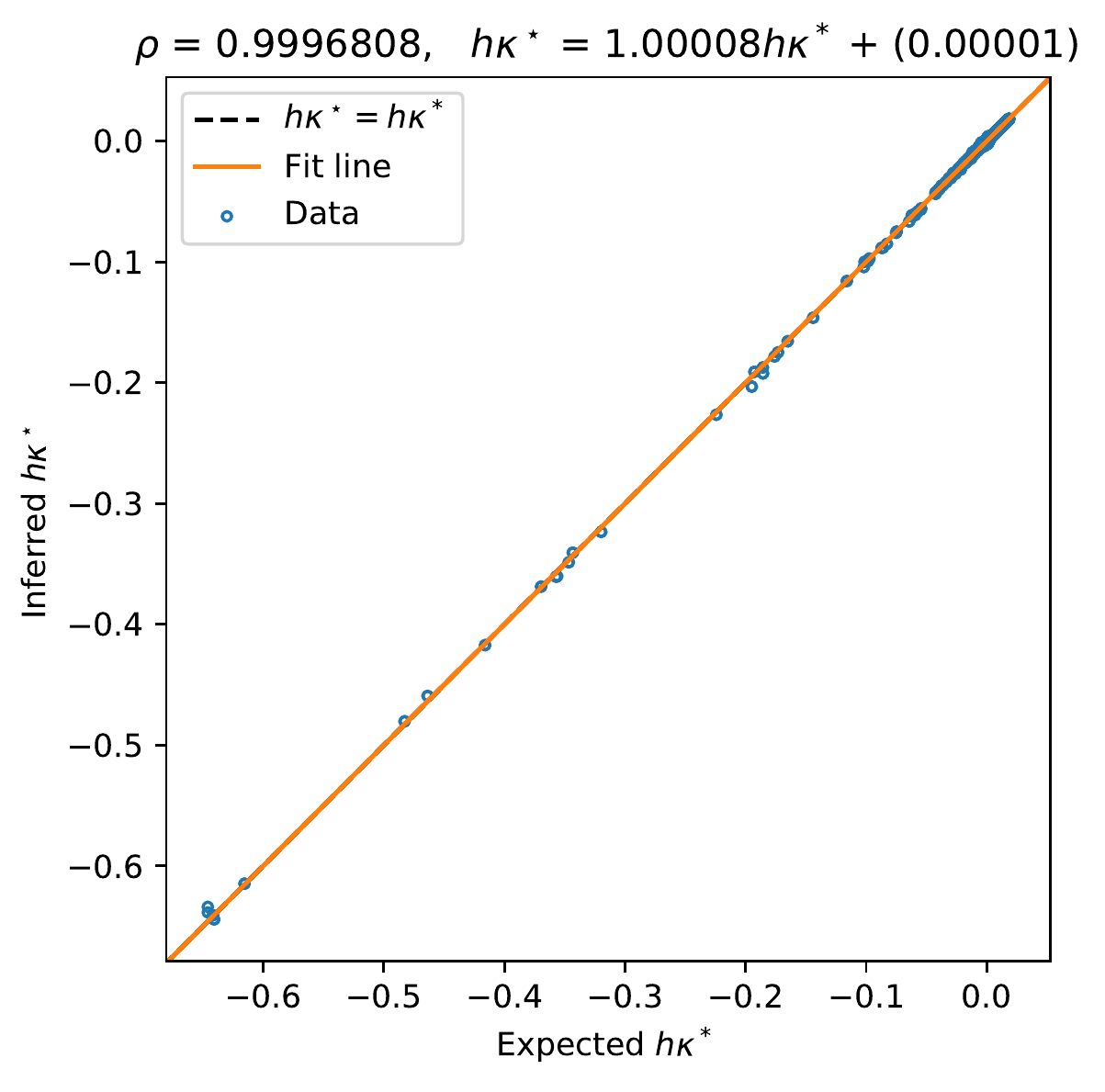}
	\end{subfigure}
	\\
	\begin{subfigure}[!t]{\textwidth}
		\caption{\footnotesize $h = 2^{-10}$}
		\label{fig:results.steep.charts.10}
	\end{subfigure}
	\\
    
	\begin{subfigure}[!t]{0.30\textwidth}
		\includegraphics[width=\textwidth]{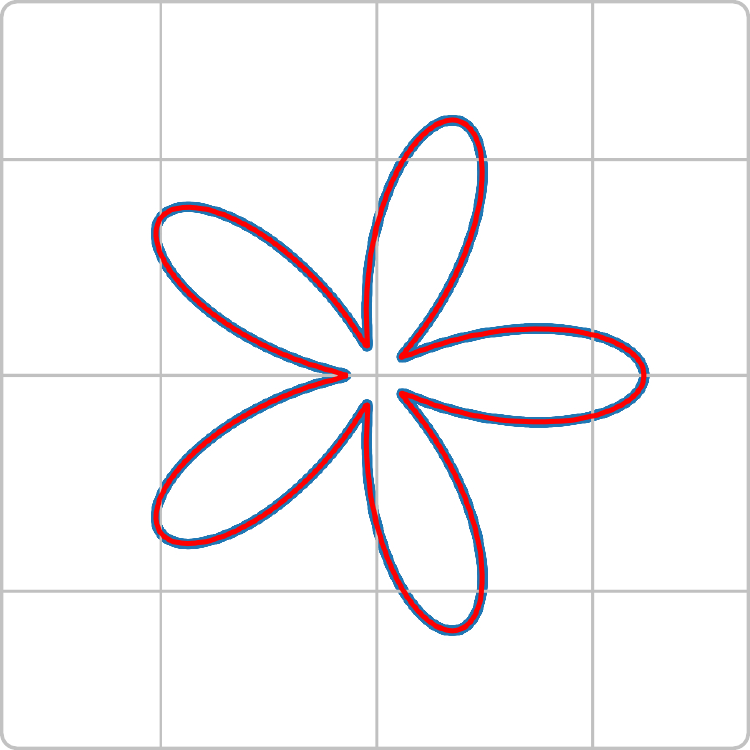}
	\end{subfigure}
	~
	\begin{subfigure}[!t]{0.33\textwidth}
		\includegraphics[width=\textwidth]{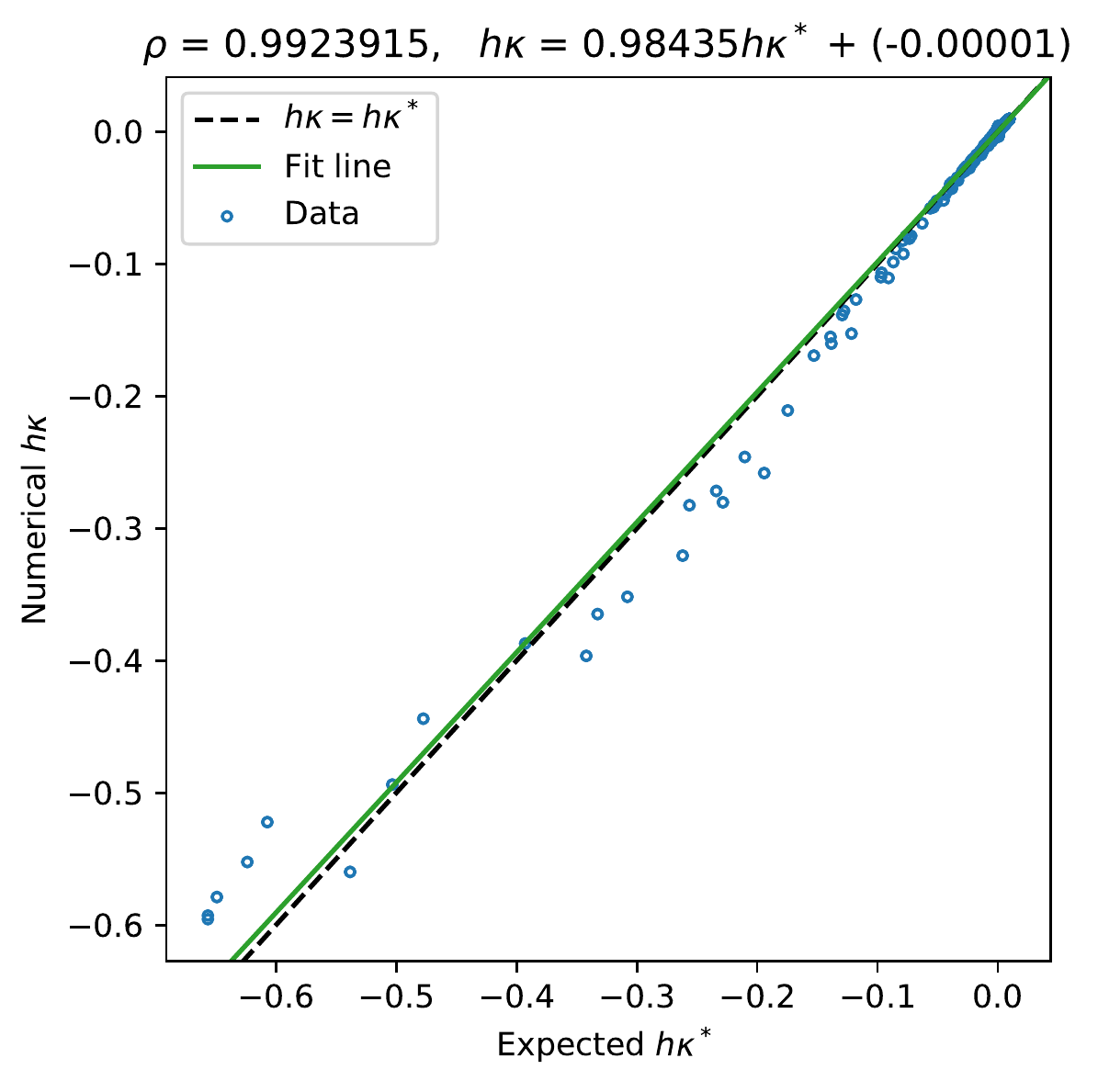}
	\end{subfigure}
	~
	\begin{subfigure}[!t]{0.33\textwidth}
		\includegraphics[width=\textwidth]{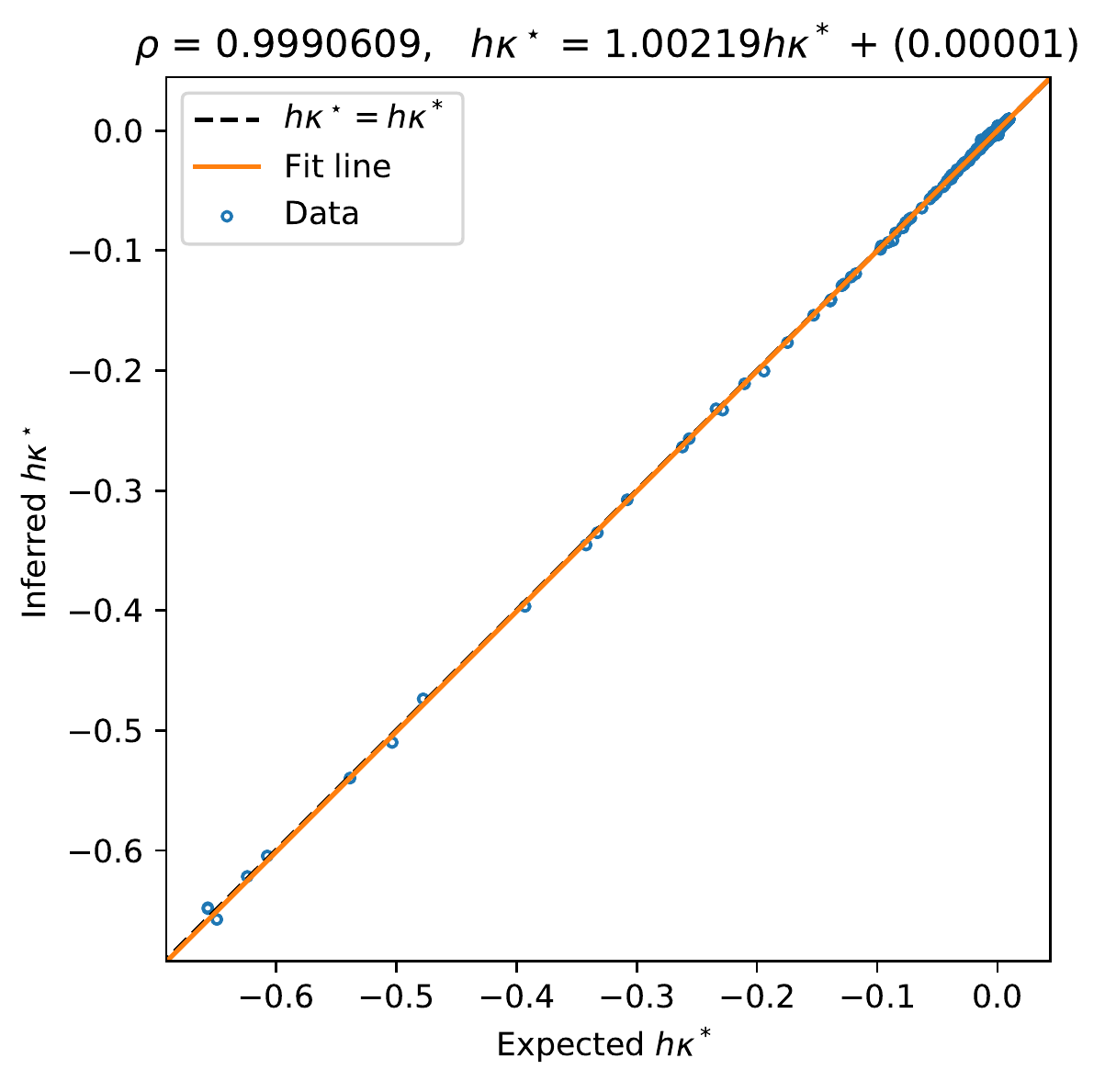}
	\end{subfigure}
	\\
	\begin{subfigure}[!t]{\textwidth}
		\caption{\footnotesize $h = 2^{-11}$}
		\label{fig:results.steep.charts.11}
	\end{subfigure}
	\\
	   
	\caption{(Cont'd.)  Steep-curvature flower-shaped interfaces and their correlation plots.  The first column shows the tested interface with concavities varying according to the grid resolution.  We have contoured $\gamma(\theta)$ in red and denoted the sampled grid points next to $\Gamma$ in blue.  The second and third columns depict fit quality for the numerical baseline (in green) and {\tt MLCurvature()} (in orange) for $\nu = 10$.  Row (d) belongs to $\eta = 9$, row (e) corresponds to $\eta = 10$, and row (f) is associated to $\eta = 11$.  (Color online.)}
	\label{fig:results.steep.charts}
\end{figure}


\colorsubsection{Convergence analysis}
\label{subsec:ConvergenceAnalysis}

\begin{figure}[!t]
	\centering
	\begin{subfigure}[b]{0.36\textwidth}
		\includegraphics[width=\textwidth]{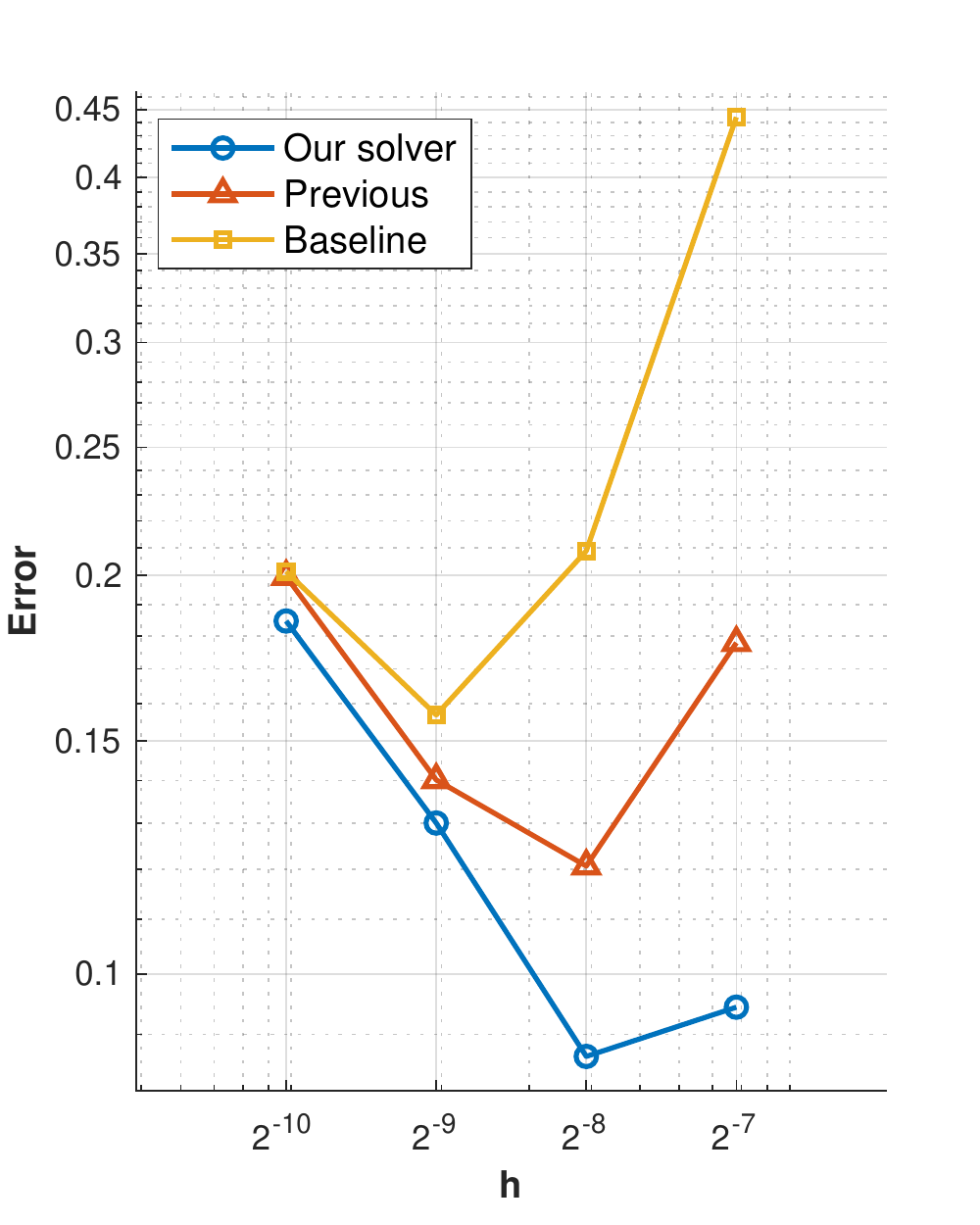}
        \caption{\footnotesize MAE}
        \label{fig:results.convergence.accuracy.mae}
    \end{subfigure}
    ~
	\begin{subfigure}[b]{0.36\textwidth}
		\includegraphics[width=\textwidth]{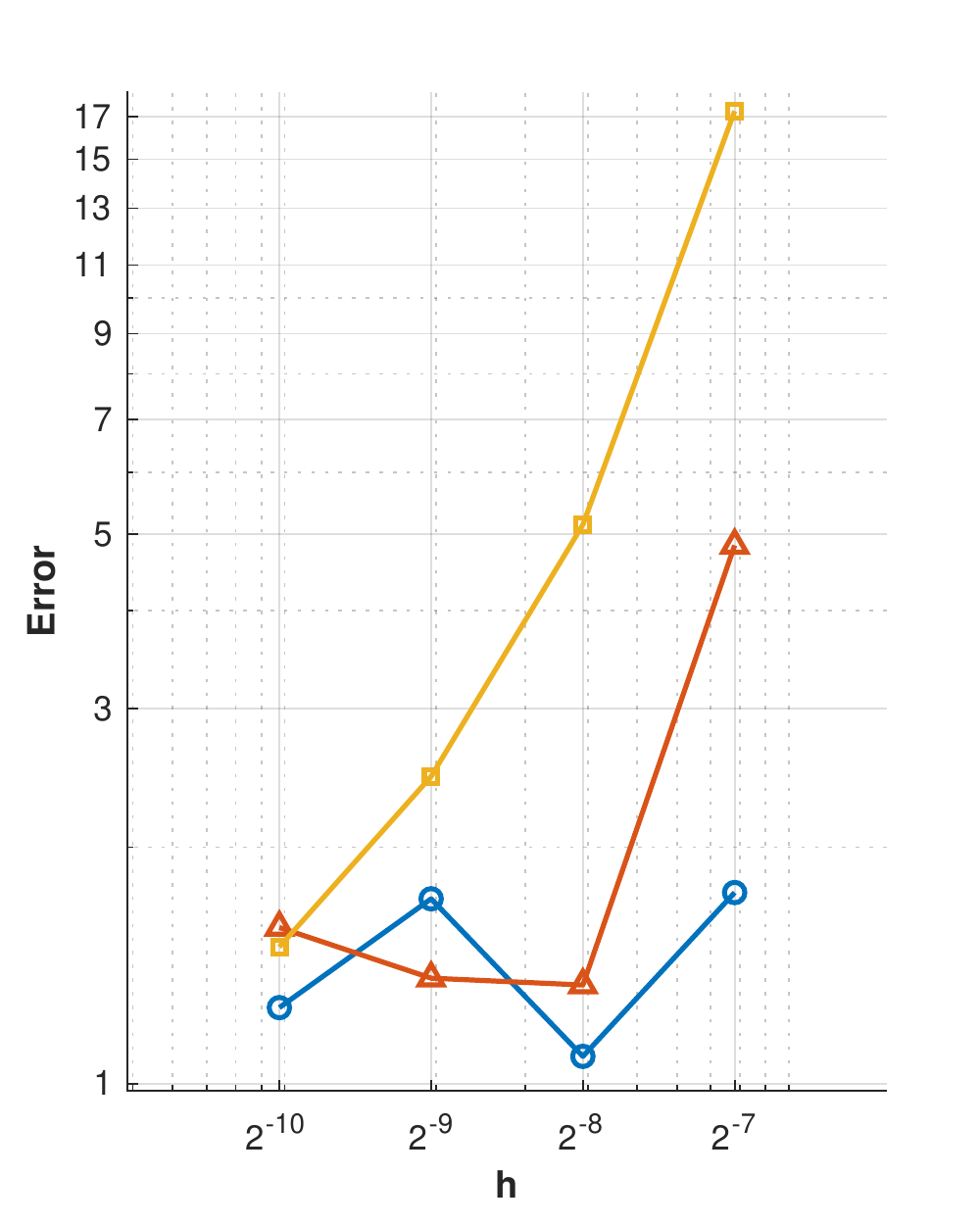}
		\caption{\footnotesize MaxAE}
		\label{fig:results.convergence.accuracy.maxae}
	\end{subfigure}
	   
	\caption{Curvature estimation accuracy in the (a) $L^1$ and (b) $L^\infty$ norms using our proposed system based on the {\tt MLCurvature()} routine, the previous hybrid approach in \cite{Larios;Gibou;HybridCurvature;2021}, and the numerical baseline for $\phi_{rose}(\vv{x})$ with $a = 0.120$, $b = 0.305$, and $\nu = 10$.  (Color online.)}
	\label{fig:results.convergence.accuracy}
\end{figure}

\begin{table}[!t]
	\centering
	\small
	\bgroup
	\def\arraystretch{1.1}%
	\begin{tabular}{rcrcrcr}
		\hline
		\rowcolor{cloud1}
		$\eta$ & {\tt MLCurvature()} & Order & Previous & Order & Baseline & Order \\
		\hline \hline
		7      &            0.094355 &     ~ & 0.177896 &     ~ & 0.444069 &     ~ \\
		\hline
		8      &            0.086604 &  0.12 & 0.120624 &  0.56 & 0.208758 &  1.09 \\
		\hline
		9      &            0.130005 & -0.59 & 0.140138 & -0.22 & 0.156913 &  0.41 \\
		\hline
		10     &            0.184768 & -0.51 & 0.199658 & -0.51 & 0.201407 & -0.36 \\
		\hline
	\end{tabular}
	\egroup
	\caption{Convergence analysis for $\kappa$'s MAE.  Evaluated interface corresponds to $\phi_{rose}(\vv{x})$'s zero-isocontour with $a = 0.120$, $b = 0.305$, and $\nu = 10$.}
	\label{tbl:results.convergence.mae}
\end{table}

\begin{table}[!t]
	\centering
	\small
	\bgroup
	\def\arraystretch{1.1}%
	\begin{tabular}{rcrcrcr}
		\hline
		\rowcolor{cloud1}
		$\eta$ & {\tt MLCurvature()} & Order & Previous & Order & Baseline & Order \\
		\hline \hline
		7      &            1.750567 &     ~ & 4.841265 &     ~ & 17.27762 &    ~ \\
		\hline
		8      &            1.083348 &  0.69 & 1.335203 &  1.86 & 5.143841 & 1.75 \\
		\hline
		9      &            1.718821 & -0.67 & 1.503125 & -0.29 & 2.459796 & 1.06 \\
		\hline
		10     &            1.249617 &  0.46 & 1.580569 & -0.21 & 1.492718 & 0.72 \\
		\hline
	\end{tabular}
	\egroup
	\caption{Convergence analysis for $\kappa$'s MaxAE.  Evaluated interface corresponds to $\phi_{rose}(\vv{x})$'s zero-isocontour with $a = 0.120$, $b = 0.305$, and $\nu = 10$.}
	\label{tbl:results.convergence.maxae}
\end{table}

Our last assessment replicates the convergence analysis conducted in \cite{Larios;Gibou;HybridCurvature;2021}.  The goal is to evaluate some fixed $\phi_{rose}(\vv{x})$ level-set function across several grid resolutions and establish the ratio between the error and $\eta$.  Here, we also use this test to compare {\tt MLCurvature()}'s accuracy with the hybrid inference system in \cite{Larios;Gibou;HybridCurvature;2021}.  For compatibility with the latter, we have performed this assessment for $\eta = 7, 8, 9, 10$.  In particular, the five-armed polar rose's shape parameters in \cref{eq:PolarRoseLevelSetFunction,eq:PolarRoseInterface} are $a = 0.120$ and $b = 0.305$.  These values characterize the steepest contour we can handle for $\eta = 7$.  \Cref{fig:results.convergence.accuracy} plots the curvature precision in the $L^1$ and $L^\infty$ norms for each of the three methods: our solver (i.e., \Cref{alg:MLCurvature}), the previous inference system in \cite{Larios;Gibou;HybridCurvature;2021}, and the numerical baseline.  For clarity, we have restricted level-set reinitialization to only $\nu = 10$ iterations and omitted results for $\nu = 20$.

\Cref{fig:results.convergence.accuracy} shows the lack of MAE and MaxAE convergence in our ML solution, thus confirming the findings in \cite{CurvatureML19} and \cite{Larios;Gibou;HybridCurvature;2021}.  However, it validates that {\tt MLCurvature()}'s precision is predominantly superior to the other methods, with noteworthy improvements for small $\eta$.  \Cref{tbl:results.convergence.mae,tbl:results.convergence.maxae} summarize these errors besides their convergence rates.  

When comparing our routine with the solver in \cite{Larios;Gibou;HybridCurvature;2021}, we can observe a consistent MAE reduction, ranging from 7\% to 47\%.  In this case, the highest gain in accuracy has occurred for $\eta = 7$.  The MaxAE, on the other hand, exhibits a less clear behavior.  Although this metric has dropped below the previous system's error in most cases, it increased by 15\% for $\eta = 9$.  This matter deserves further research, and we plan to investigate whether we can employ reasonably larger data sets in fine grids.  Nevertheless, our strategy has outperformed the numerical baseline in all scenarios.  Given these results, we believe our framework embodies a much more promising and structured method than the earlier study in \cite{Larios;Gibou;HybridCurvature;2021}.


\colorsection{Conclusions}
\label{sec:Conclusions}

We have presented a curvature hybrid solver for two-dimensional, smooth interfaces based on error neural modeling and quantification.  The core of our system is a multilayer perceptron, $\mathcal{F}_\kappa(\cdot)$, designed with the methodology outlined in \cite{Larios;Gibou;ECNetSemiLagrangian;2021}.  $\mathcal{F}_\kappa(\cdot)$'s task is to ingest preprocessed level-set, gradient, and curvature information to fix the numerical curvature estimation at the interface.  Through several case studies on irregular free boundaries, we have established the effectiveness of $\mathcal{F}_\kappa(\cdot)$ as the main constituent of our {\tt MLCurvature()} algorithm.  In particular, including normal vector components and regularization has increased the capacity and generalization in $\mathcal{F}_\kappa(\cdot)$.  Similarly, sample reorientation and symmetry preservation have enhanced neural optimization without sacrificing simplicity.  At the same time, these features have made it possible for {\tt MLCurvature()} to outperform previous inference systems \cite{LALariosFGibou;LSCurvatureML;2021, Larios;Gibou;HybridCurvature;2021}, especially around under-resolved regions.

Our experiments have shown the flexibility of the proposed strategy, irrespective of the mesh size.  Compared to \cite{Larios;Gibou;HybridCurvature;2021}, for example, we have minimized outlying effects by incorporating layer-wise {\tt L2} regularization and exploiting curvature symmetry invariance.  Likewise, we have made network training scalable by introducing dimensionless parametrization and constraining data augmentation to a single reflection.  Alongside probabilistic subsampling, these measures have helped reduce data-set size.  However, it is always possible to adjust a few input constants in \Cref{alg:GenerateCircularDataSet,alg:GenerateSinusoidalDataSet} to produce more learning tuples in high-resolution grids.  As noted in \Cref{subsec:ConvergenceAnalysis}, we plan to investigate whether carefully tuning those constants suffices to further curvature accuracy and amend the situation in \cref{fig:results.convergence.accuracy.maxae}.

We believe our framework provides the basis for promising extensions and enhancements.  For instance, since we have stated the numerical curvature problem from an error perspective, one could first build an estimator to compensate for inaccuracies at the interface nodes.  Then, we could train a model to infer optimal combination coefficients to interpolate $\kappa$ at $\Gamma$.  These coefficients could be extracted from a restricted system, similar to the learned discretizations of \cite{Zhuang;etal;LrndDiscForPassSclrAdvctn2D;2021}.  Another possibility is the introduction of an under-resolution classifier to toggle the ML component as needed.  This idea \cite{TroubledCellIndicator18, ShockDetector20, Buhendwa;Bezgin;Adams;IRinLSwithML;2021} could free us from hand-tuning $h\kappa_{\min}^*$ in \Crefrange{alg:MLCurvature}{alg:GenerateSinusoidalDataSet} and possibly $h\kappa_{\min}^{up}$ in the {\tt MLCurvature()} function.  Lastly, we could include data from neighboring interface nodes, analogous to \cite{Franca;Oishi;MLCurvatureFrontTracking;2022}, but as analytical constraints \cite{Beucler;etal;EnforcingAnalyticConstNnets;2021}.  This way, we could enforce some smoothness and reduce sharp curvature variations between ``successive'' nodes.

Currently, we are working on the three-dimensional extension of our error-correcting approach.  And, to dampen the problem's susceptibility to data explosion, we are considering the sampling procedure delineated in \cite{VOFCurvature3DML19}.  We also leave to future endeavors the integration of {\tt MLCurvature()} into the semi-Lagrangian advection scheme of \cite{Larios;Gibou;ECNetSemiLagrangian;2021}.  There, we would use {\tt MLCurvature()}'s output to train and deploy the neural network within the transport function.  In like manner, it remains to show the interplay of our curvature solver and a full-fledged FBP application.  The latter would help us better understand the effects of \Cref{alg:MLCurvature} on the dynamics of moving fronts.



{\footnotesize
\biboptions{sort&compress}
\bibliographystyle{unsrt}
\bibliography{references}}

\begin{thebibliography}{10}

\bibitem{Friedman10}
A.~Friedman.
\newblock {\em {Variational Principles of Free-Boundary Problems}}.
\newblock Dover Publications, 2010.

\bibitem{Osher1988}
S.~Osher and J.~A. Sethian.
\newblock {Fronts propagating with curvature-dependent speed: Algorithms based
  on Hamilton--Jacobi formulations}.
\newblock {\em J. Comput. Phys.}, 79(1):12--49, November 1988.

\bibitem{Popinet;NumModelsOfSurfTension;18}
S.~Popinet.
\newblock {Numerical models of surface tension}.
\newblock {\em Annu. Rev. Fluid Mech.}, 50(1):49--75, January 2018.

\bibitem{Sussman;Smereka;Osher:94:A-Level-Set-Approach}
M.~Sussman, P.~Smereka, and S.~Osher.
\newblock {A level set approach for computing solutions to incompressible
  two-phase flow}.
\newblock {\em J. Comput. Phys.}, 114(1):146--159, September 1994.

\bibitem{Sussman;Fatemi;Smereka;etal:98:An-Improved-Level-Se}
M.~Sussman, E.~Fatemi, P.~Smereka, and S.~Osher.
\newblock {An improved level set method for incompressible two-phase flows}.
\newblock {\em {Comput. \& Fluids}}, 27(5-6):663--680, June 1998.

\bibitem{Gibou;Chen;Nguyen;etal:07:A-level-set-based-sh}
F.~Gibou, L.~Chen, D.~Nguyen, and S.~Banerjee.
\newblock {A level set based sharp interface method for the multiphase
  incompressible Navier--Stokes equations with phase change}.
\newblock {\em J. Comput. Phys.}, 222(2):536--555, March 2007.

\bibitem{Theillard:2019aa}
M.~Theillard, F.~Gibou, and D.~Saintillan.
\newblock {Sharp numerical simulation of incompressible two-phase flows}.
\newblock {\em J. Comput. Phys.}, 2019.

\bibitem{Losasso;Gibou;Fedkiw:04:Simulating-Water-and}
F.~Losasso, F.~Gibou, and R.~Fedkiw.
\newblock {Simulating water and smoke with an octree data structure}.
\newblock {\em ACM Trans. Graph. (SIGGRAPH Proc.)}, 23(3):457--462, August
  2004.

\bibitem{Losasso:2006aa}
F.~Losasso, T.~Shinar, A.~Selle, and R.~Fedkiw.
\newblock {Multiple interacting liquids}.
\newblock {\em SIGGRAPH ACM TOG}, 25(3):812--819, July 2006.

\bibitem{Gibou:2019aa}
F.~Gibou, D.~Hyde, and R.~Fedkiw.
\newblock {Sharp interface approaches and deep learning techniques for
  multiphase flows}.
\newblock {\em J. Comput. Phys.}, 380:442--463, 2019.

\bibitem{Egan;etal;DirNumSimIncompFlowsOctree;2021}
R.~Egan, A.~Guittet, F.~Temprano-Coleto, T.~Isaac, F.~J. Peaudecerf, J.~R.
  Landel, P.~Luzzatto-Fegiz, C.~Burstedde, and F.~Gibou.
\newblock {Direct numerical simulation of incompressible flows on parallel
  octree grids}.
\newblock {\em J. Comput. Phys.}, 428:110084, March 2021.

\bibitem{Chen;Min;Gibou:09:A-numerical-scheme-f}
H.~Chen, C.~Min, and F.~Gibou.
\newblock {A numerical scheme for the Stefan problem on adaptive Cartesian
  grids with supralinear convergence rate}.
\newblock {\em J. Comput. Phys.}, 228(16):5803--5818, September 2009.

\bibitem{Papac;Gibou;Ratsch:10:Efficient-symmetric-}
J.~Papac, F.~Gibou, and C.~Ratsch.
\newblock {Efficient symmetric discretization for the {P}oisson, heat and
  {S}tefan-type problems with {R}obin boundary conditions}.
\newblock {\em J. Comput. Phys.}, 229(3):875--889, 2010.

\bibitem{Papac;Helgadottir;Ratsch;etal:13:A-level-set-approach}
J.~Papac, A.~Helgadottir, C.~Ratsch, and F.~Gibou.
\newblock {A level set approach for diffusion and {S}tefan-type problems with
  {R}obin boundary conditions on quadtree/octree adaptive {C}artesian grids}.
\newblock {\em J. Comput. Phys.}, 233:241, 2013.

\bibitem{Mirzadeh;Gibou:14:A-conservative-discr}
M.~Mirzadeh and F.~Gibou.
\newblock {A conservative discretization of the {P}oisson--{N}ernst--{P}lanck
  equations on adaptive {C}artesian grids}.
\newblock {\em J. Comput. Phys.}, 274:633--653, 2014.

\bibitem{Theillard;Gibou;Pollock:14:A-Sharp-Computationa}
M.~Theillard, F.~Gibou, and T.~Pollock.
\newblock {A sharp computational method for the simulation of the
  solidification of binary alloys}.
\newblock {\em J. Sci. Comput.}, 63:330--354, 2015.

\bibitem{Boudon;etal;3DPlantMorphogenesis;2015}
F.~Boudon, J.~Chopard, O.~Ali, B.~Gilles, O.~Hamant, A.~Boudaoud, J.~Traas, and
  C.~Godin.
\newblock {A computational framework for 3D mechanical modeling of plant
  morphogenesis with cellular resolution}.
\newblock {\em PLoS Comput. Biol.}, 11(1):e1003950, January 2015.

\bibitem{Ocko;Heyde;Mahadevan;MorphTermiteMounds;2019}
S.~A. Ocko, A.~Heyde, and L.~Mahadevan.
\newblock {Morphogenesis of termite mounds}.
\newblock {\em Proc. Natl. Acad. Sci. USA}, 116(9):3379--3384, February 2019.

\bibitem{AliasBuenzli20}
M.~A. Alias and P.~R. Buenzli.
\newblock {A level-set method for the evolution of cells and tissue during
  curvature-controlled growth}.
\newblock {\em Int. J. Numer. Methods Biomed. Eng.}, 36(1):e3279, January 2020.

\bibitem{Lervag;CalcCurvatureLSM;2014}
K.~Y. Lerv{\r{a}}g.
\newblock {Calculation of interface curvature with the level-set method}.
\newblock \url{https://arxiv.org/abs/1407.7340}, July 2014.

\bibitem{Sethian:99:Level-set-methods-an}
J.~A. Sethian.
\newblock {\em {Level Set Methods and Fast Marching Methods}}.
\newblock Cambridge Monogr. Appl. Comput. Math. Cambridge University Press,
  Cambridge, UK, 2nd edition, 1999.

\bibitem{Osher;Fedkiw:02:Level-Set-Methods-an}
S.~Osher and R.~Fedkiw.
\newblock {\em {Level Set Methods and Dynamic Implicit Surfaces}}.
\newblock Appl. Math. Sci. 153. Springer, Cham, 2002.

\bibitem{GFO18}
F.~Gibou, R.~Fedkiw, and S.~Osher.
\newblock {A review of level-set methods and some recent applications}.
\newblock {\em J. Comput. Phys.}, 353:82--109, January 2018.

\bibitem{Hirt;Nichols:81:Volume-of-Fluid-VOF-}
C.~W. Hirt and B.~D. Nichols.
\newblock {Volume of fluid (VOF) method for the dynamics of free boundaries}.
\newblock {\em J. Comput. Phys.}, 39:201--225, 1981.

\bibitem{QB10}
R.~S. Qin and H.~K. Bhadeshia.
\newblock {Phase field method}.
\newblock {\em Materials Sci. Tech.}, 26(7):803--811, 2010.

\bibitem{Tryggvason;Bunner;Esmaeeli;etal:01:A-Front-Tracking-Met}
G.~Tryggvason, B.~Bunner, A.~Esmaeeli, D.~Juric, N.~Al-Rawahi, W.~Tauber,
  J.~Han, S.~Nas, and Y.-J. Jan.
\newblock {A front-tracking method for the computations of multiphase flow}.
\newblock {\em J. Comput. Phys.}, 169(2):708--759, May 2001.

\bibitem{Chene;Min;Gibou:08:Second-order-accurat}
A.~du~Ch{\'e}n{\'e}, C.~Min, and F.~Gibou.
\newblock {Second-order accurate computation of curvatures in a level set
  framework using novel high-order reinitialization schemes}.
\newblock {\em J. Sci. Comput.}, 35:114--131, June 2008.

\bibitem{Zhao:04:A-Fast-Sweeping-Meth}
H.~Zhao.
\newblock {A fast sweeping method for eikonal equations}.
\newblock {\em Math. Comp.}, 74:603--627, 2005.

\bibitem{Detrixhe;Gibou;Min:13:A-parallel-fast-swee}
M.~Detrixhe, F.~Gibou, and C.~Min.
\newblock {A parallel fast sweeping method for the eikonal equation}.
\newblock {\em J. Comput. Phys.}, 237:46--55, March 2013.

\bibitem{Macklin;Lowengrub;ImprovedCurvatureAppTumorGrowth;2006}
P.~Macklin and J.~Lowengrub.
\newblock {An improved geometry-aware curvature discretization for level set
  methods: Application to tumor growth}.
\newblock {\em J. Comput. Phys.}, 215(2):392--401, July 2006.

\bibitem{LALariosFGibou;LSCurvatureML;2021}
L.~{\'A}. Larios-C{\'a}rdenas and F.~Gibou.
\newblock {A deep learning approach for the computation of curvature in the
  level-set method}.
\newblock {\em SIAM J. Sci. Comput.}, 43(3):A1754--A1779, January 2021.

\bibitem{Larios;Gibou;HybridCurvature;2021}
L.~{\'A}. Larios-C{\'a}rdenas and F.~Gibou.
\newblock {A hybrid inference system for improved curvature estimation in the
  level-set method using machine learning}.
\newblock {\em J. Comput. Phys.}, page 111291, May 2022.

\bibitem{CurvatureML19}
Y.~Qi, J.~Lu, R.~Scardovelli, S.~Zaleski, and G.~Tryggvason.
\newblock {Computing curvature for volume of fluid methods using machine
  learning}.
\newblock {\em J. Comput. Phys.}, 377:155--161, 2019.

\bibitem{A18}
C.~C. Aggarwal.
\newblock {\em {Neural Networks and Deep Learning -- A Textbook}}.
\newblock Springer, Cham, 2018.

\bibitem{Mehta19}
P.~Mehta, M.~Bukov, C.~Wang, A.~G.~R. Day, C.~Richardson, C.~K. Fisher, and
  D.~J. Schwabd.
\newblock {A high-bias, low-variance introduction to machine learning for
  physicists}.
\newblock {\em Phys. Rep.}, 810:1--124, May 2019.

\bibitem{VOFCurvature3DML19}
H.~V. Patel, A.~Panda, J.~A.~M. Kuipers, and E.~A. J.~F. Peters.
\newblock {Computing interface curvature from volume fractions: A machine
  learning approach}.
\newblock {\em {Comput. \& Fluids}}, 193:104263, October 2019.

\bibitem{DespresJourdren;MLDesignOfVOF;20}
B.~Despr{\'e}s and H.~Jourdren.
\newblock {Machine learning design of volume of fluid schemes for compressible
  flows}.
\newblock {\em J. Comput. Phys.}, 408(1):109275, May 2020.

\bibitem{NPLIC20}
M.~Ataei, M.~Bussmann, V.~Shaayegan, F.~Costa, S.~Han, and C.~B. Park.
\newblock {NPLIC: A machine learning approach to piecewise linear interface
  construction}.
\newblock {\url{https://arxiv.org/abs/2007.04244}}, January 2021.

\bibitem{Buhendwa;Bezgin;Adams;IRinLSwithML;2021}
A.~B. Buhendwa, D.~A. Bezgin, and N.~Adams.
\newblock {Consistent and symmetry preserving data-driven interface
  reconstruction for the level-set method}.
\newblock \url{https://arxiv.org/abs/2104.11578}, April 2021.

\bibitem{Franca;Oishi;MLCurvatureFrontTracking;2022}
H.~L. Fran{\c{c}}a and C.~M. Oishi.
\newblock {A machine learning strategy for computing interface curvature in
  front-tracking methods}.
\newblock {\em J. Comput. Phys.}, 450:110860, February 2022.

\bibitem{Larios;Gibou;ECNetSemiLagrangian;2021}
L.~{\'A}. Larios-C{\'a}rdenas and F.~Gibou.
\newblock {Error-correcting neural networks for semi-Lagrangian advection in
  the level-set method}.
\newblock \url{https://arxiv.org/abs/2110.11611}, October 2021.

\bibitem{Pathak;etal;MLToAugCoarseGridCFD;2020}
J.~Pathak, M.~Mustafa, K.~Kashinath, E.~Motheau, T.~Kurth, and M.~Day.
\newblock {Using machine learning to augment coarse-grid computational fluid
  dynamics simulations}.
\newblock \url{https://arxiv.org/abs/2010.00072}, 2020.

\bibitem{Dong;Loy;He;SuperResolution;2014}
C.~Dong, C.~C. Loy, K.~He, and X.~Tang.
\newblock {Learning a deep convolutional network for image super-resolution}.
\newblock In D.~Fleet, T.~Pajdla, B.~Schiele, and T.~Tuytelaars, editors, {\em
  Computer Vision -- ECCV 2014}, pages 184--199, Cham, 2014. Springer
  International Publishing.

\bibitem{Min;Gibou:07:A-second-order-accur}
C.~Min and F.~Gibou.
\newblock {A second order accurate level set method on non-graded adaptive
  Cartesian grids}.
\newblock {\em J. Comput. Phys.}, 225(1):300--321, July 2007.

\bibitem{Shu;Osher:89:Efficient-Implementa}
C.-W. Shu and S.~Osher.
\newblock {Efficient implementation of essentially non-oscillatory shock
  capturing schemes, II}.
\newblock {\em J. Comput. Phys.}, 83(1):32--78, July 1989.

\bibitem{Jiang;Peng:00:Weighted-ENO-Schemes}
G.-S. Jiang and D.~Peng.
\newblock {Weighted ENO schemes for Hamilton--Jacobi equations}.
\newblock {\em SIAM J. Sci. Comput.}, 21(6):2126--2143, 2000.

\bibitem{Mirzadeh;etal:16:Parallel-level-set}
M.~Mirzadeh, A.~Guittet, C.~Burstedde, and F.~Gibou.
\newblock {Parallel level-set methods on adaptive tree-based grids}.
\newblock {\em J. Comput. Phys.}, 322:345--364, October 2016.

\bibitem{Burstedde;Wilcox;Ghattas:11:p4est:-Scalable-Algo}
C.~Burstedde, L.~C. Wilcox, and O.~Ghattas.
\newblock {\texttt{p4est}: Scalable algorithms for parallel adaptive mesh
  refinement on forests of octrees}.
\newblock {\em SIAM J. Sci. Comput.}, 33(3):1103--1133, May 2011.

\bibitem{Strain1999}
J.~Strain.
\newblock {Tree methods for moving interfaces}.
\newblock {\em J. Comput. Phys.}, 151(2):616--648, May 1999.

\bibitem{scikit-learn11}
F.~Pedregosa, G.~Varoquaux, A.~Gramfort, V.~Michel, B.~Thirion, O.~Grisel,
  M.~Blondel, P.~Prettenhofer, R.~Weiss, V.~Dubourg, J.~Vanderplas, A.~Passos,
  D.~Cournapeau, M.~Brucher, M.~Perrot, and E.~Duchesnay.
\newblock {Scikit-learn: Machine learning in Python}.
\newblock {\em J. Mach. Learn. Res.}, 12:2825--2830, 2011.

\bibitem{Turk;Pentland;Eigenfaces;1991}
M.~Turk and A.~Pentland.
\newblock {Eigenfaces for recognition}.
\newblock {\em J. Cogn. Neurosci.}, 3(1):71--86, 1991.

\bibitem{Parker;CS170A;2016}
D.~S. Parker.
\newblock {Exploring the Matrix -- Adventures in Modeling with Matlab}.
\newblock UCLA Course Reader Solutions, January 2016.

\bibitem{BKOS00}
M.~de~Berg, M.~van Kreveld, M.~Overmars, and O.~Schwarzkopf.
\newblock {\em {Computational Geometry - Algorithms and Applications}}.
\newblock Springer, Cham, 2nd edition, 2000.

\bibitem{Min2004}
C.~Min.
\newblock {Local level set method in high dimension and codimension}.
\newblock {\em J. Comput. Phys.}, 200(1):368--382, October 2004.

\bibitem{Swokowski88}
E.~W. Swokowski.
\newblock {\em {Calculus with Analytic Geometry}}.
\newblock PWS Publishers, 4th edition, 1988.

\bibitem{CGUsingOpenGL01}
Jr. F.~S.~Hill.
\newblock {\em {Computer Graphics Using OpenGL}}.
\newblock Prentice-Hall Inc., 2nd edition, 2001.

\bibitem{Heath;SciComput;2018}
M.~T. Heath.
\newblock {\em {Scientific Computing: An Introductory Survey}}.
\newblock SIAM, Philadelphia, revised 2nd edition, 2018.

\bibitem{ComputerAnimation08}
R.~Parent.
\newblock {\em {Computer Animation: Algorithms and Techniques}}.
\newblock Morgan Kaufmann, 2nd edition, 2008.

\bibitem{Boost;2019}
{The Boost Community}.
\newblock {Boost C++ libraries}.
\newblock \url{https://www.boost.org}, August 2019.
\newblock v1.71.0.

\bibitem{Tensorflow15}
M.~Abadi, A.~Agarwal, P.~Barham, E.~Brevdo, Z.~Chen, C.~Citro, G.~S. Corrado,
  A.~Davis, J.~Dean, M.~Devin, S.~Ghemawat, I.~Goodfellow, A.~Harp, G.~Irving,
  M.~Isard, R.~Jozefowicz, Y.~Jia, L.~Kaiser, M.~Kudlur, J.~Levenberg,
  D.~Man{\'e}, M.~Schuster, R.~Monga, S.~Moore, D.~Murray, C.~Olah, J.~Shlens,
  B.~Steiner, I.~Sutskever, K.~Talwar, P.~Tucker, V.~Vanhoucke, V.~Vasudevan,
  F.~Vi{\'e}gas, O.~Vinyals, P.~Warden, M.~Wattenberg, M.~Wicke, Y.~Yu, and
  X.~Zheng.
\newblock {TensorFlow: Large-Scale Machine Learning on Heterogeneous Systems}.
\newblock \url{https://www.tensorflow.org}, 2015.

\bibitem{Keras15}
F.~Chollet et~al.
\newblock {Keras}.
\newblock \url{https://keras.io}, 2015.

\bibitem{Pandas2010}
W.~McKinney.
\newblock {Data Structures for Statistical Computing in Python}.
\newblock In {S. van der Walt and J. Millman}, editor, {\em {Proceedings of the
  9th Python in Science Conference}}, pages 56--61, 2010.
\newblock \url{https://pandas.pydata.org}.

\bibitem{LeCun;EfficientBackProp;98}
Y.~A. LeCun, L.~Bottou, G.~B. Orr, and K.-R. M{\"u}ller.
\newblock {\em {Efficient BackProp}}, volume 7700 of {\em Lecture Notes in
  Comput. Sci.}, pages 9--48.
\newblock Springer Berlin Heidelberg, Berlin, Heidelberg, 2012.

\bibitem{json;2021}
N.~Lohmann.
\newblock {JSON for modern C++}.
\newblock \url{https://github.com/nlohmann/json}, August 2020.
\newblock v3.9.1.

\bibitem{frugally-deep;2021}
T.~Hermann.
\newblock {Frugally-deep}.
\newblock \url{https://github.com/Dobiasd/frugally-deep}, February 2021.
\newblock v0.15.2.

\bibitem{openblas;2021}
Z.~Xianyi and M.~Kroeker.
\newblock {OpenBLAS: An optimized BLAS library}.
\newblock \url{https://github.com/xianyi/OpenBLAS}, July 2021.
\newblock v0.3.17.

\bibitem{Zhuang;etal;LrndDiscForPassSclrAdvctn2D;2021}
J.~Zhuang, D.~Kochkov, Y.~Bar-Sinai, M.~P. Brenner, and S.~Hoyer.
\newblock {Learned discretizations for passive scalar advection in a
  two-dimensional turbulent flow}.
\newblock {\em Phys. Rev. Fluids}, 6(6):064605, June 2021.
\newblock \url{https://github.com/google-research/data-driven-advection}.

\bibitem{TroubledCellIndicator18}
D.~Ray and J.~S. Hesthaven.
\newblock {An artificial neural network as a troubled-cell indicator}.
\newblock {\em J. Comp. Phys.}, 367:166--191, April 2018.

\bibitem{ShockDetector20}
N.~R. Morgan, S.~Tokareva, X.~Liu, and A.~D. Morgan.
\newblock {A machine learning approach for detecting shocks with high-order
  hydrodynamic methods}.
\newblock {\em AIAA SciTech Forum}, January 2020.

\bibitem{Beucler;etal;EnforcingAnalyticConstNnets;2021}
T.~Beucler, M.~Pritchard, S.~Rasp, J.~Ott, P.~Baldi, and P.~Gentine.
\newblock {Enforcing analytic constraints in neural networks emulating physical
  systems}.
\newblock {\em Phys. Rev. Lett.}, 126(9):098302, March 2021.

\end{thebibliography}


\newpage
\color{navy}\sffamily\section*{Statements and declarations}\rmfamily\color{black}

\color{navy}\subsection*{Funding}\color{black}

The authors declare that no funds, grants, or other support were received during the preparation of this manuscript.


\color{navy}\subsection*{Competing interests}\color{black}

The authors have no relevant financial or non-financial interests to disclose.


\color{navy}\subsection*{Author contributions}\color{black}

\textbf{Luis Ángel Larios-Cárdenas}: Conceptualization, Methodology, Software, Validation, Formal analysis, Investigation, Resources, Data curation, Writing – Original draft, Visualization.  \textbf{Frédéric Gibou}: Conceptualization, Methodology, Resources, Writing – Review \& editing, Supervision, Project administration, Funding acquisition.


\color{navy}\subsection*{Data availability}\color{black}

The training datasets generated during the current study are not publicly available due to space limitations but are available from the corresponding author on reasonable request.  Our neural networks, their preprocessing objects, and the testing flower-interface data sets we used, however, are publicly available at \url{https://github.com/UCSB-CASL/Curvature_ECNet_2D}.

\end{document}